\documentclass[mnsc,nonblindrev]{informs3a}

\OneAndAHalfSpacedXI 

\usepackage{natbib}
 \bibpunct[, ]{(}{)}{,}{a}{}{,}%
 \def\bibfont{\small}%
 \usepackage{amssymb}
 \usepackage{pifont}
\usepackage{makecell}

\renewcommand{\bibfont}{\footnotesize\linespread{1}\selectfont}

\renewcommand{\bar}{\widebar}

\newcommand{\R}{ \mathbb{R} }
\renewcommand{\P}{ \mathbb{P}}

\newcommand{\N}{ \mathbb{N} }

\newcommand{\E}{ \mathbb{E} }

\newcommand{\I}{ \mathbb{I} }

\newcommand{\cU}{\mathcal{U}}
\newcommand{\cX}{\mathcal{X}}
\newcommand{\cI}{\mathcal{I}}
\newcommand{\cE}{\mathcal{E}}
\newcommand{\X}{ \mathcal{X} }
\newcommand{\G}{ \mathcal{G} }

\newcommand{\mcR}{ \mathcal{R} }

\newcommand{\bD}{{\bm D}}

\usepackage{parskip} 

\newcommand{\intd}{\mathrm{d}}

\usepackage{graphicx}   
\usepackage{subcaption} 

\usepackage{array} 
\usepackage{cases}
\usepackage{pifont}
\usepackage{braket}
\usepackage{algorithm}
\usepackage{algorithmicx}
\usepackage[noend]{algpseudocode}
\def\1{(\mathrm{\uppercase\expandafter{\romannumeral1}})}
\def\2{(\mathrm{\uppercase\expandafter{\romannumeral2}})}
\def\O{\mathcal O}
\newcommand{\mh}{\mathcal{H}}
\RequirePackage[colorlinks,citecolor=blue,linkcolor=blue,urlcolor=blue]{hyperref}
\usepackage{enumitem}
\DeclareMathAlphabet{\mathscr}{OT1}{pzc}{m}{it}

\usepackage{enumitem}
\usepackage{tablefootnote}

\def\O{\mathcal O}

\def\X{\mathcal X}

\usepackage{bm}
\usepackage{hyperref}

\usepackage{xurl}

\TheoremsNumberedThrough     
\ECRepeatTheorems

\EquationsNumberedThrough    

\usepackage{mathabx}
\renewcommand{\hat}{\widehat}
\renewcommand{\tilde}{\widetilde}

\begin{document}



\RUNAUTHOR{Chen, Lyu, Yuan, Zhou}

\RUNTITLE{Nonstationary Newsvendor from Uncensored to Censored Demand}

\TITLE{Learning When to Restart: Nonstationary Newsvendor from Uncensored to Censored Demand}

\ARTICLEAUTHORS{%
 \AUTHOR{Xin Chen\footnotemark[1]}
\AFF{H. Milton Stewart School of Industrial and Systems Engineering, Georgia Institute of Technology, Atlanta, Georgia 30332, USA, \EMAIL{xin.chen@isye.gatech.edu}}
 \AUTHOR{Jiameng Lyu\footnotemark[1]}
 \AFF{Department of Management Science, School of Management, Fudan  University, Shanghai, 200433, China,  \EMAIL{jiamenglyu@fudan.edu.cn}}
 \AUTHOR{Shilin Yuan\footnotemark[1]}
\AFF{Department of Decisions, Operations and Technology, Chinese University of Hong Kong, Shatin, Hong Kong, \EMAIL{shilinyuan@cuhk.edu.hk}}
 \AUTHOR{Yuan Zhou\footnotemark[1]}
 \AFF{Yau Mathematical Sciences Center \& Department of Mathematical  Sciences, Tsinghua University, Beijing 100084, China, \EMAIL{yuan-zhou@tsinghua.edu.cn}}
 } 
\renewcommand{\thefootnote}{\fnsymbol{footnote}}
\footnotetext[1]{Author names listed in alphabetical order.}
\renewcommand{\thefootnote}{\arabic{footnote}}

\ABSTRACT{
{We study nonstationary newsvendor problems under nonparametric demand models and general distributional measures of nonstationarity, addressing the practical challenges of unknown degree of nonstationarity and demand censoring. 
 We propose a novel distributional-detection-and-restart framework for learning in nonstationary environments, and instantiate it through two efficient algorithms for the uncensored and censored demand settings. 
 The algorithms are fully adaptive, requiring no prior knowledge of the degree and type of nonstationarity, and offer a flexible yet powerful approach to handling both abrupt and gradual changes in nonstationary environments.
We establish a comprehensive optimality theory for our algorithms by deriving matching regret upper and lower bounds under both general and refined structural conditions with nontrivial proof techniques that are of independent interest.
Numerical experiments using real-world datasets, including nurse staffing data for emergency departments and COVID-19 test demand data, showcase the algorithms' superior and robust empirical performance.
While motivated by the newsvendor problem, the distributional-detection-and-restart framework applies broadly to a wide class of nonstationary stochastic optimization problems. Managerially, our framework provides a practical, easy-to-deploy, and theoretically grounded solution for decision-making under nonstationarity.}
}

\maketitle

\section{Introduction}

The newsvendor model is a fundamental framework in operations management that captures the essence of inventory decisions under demand uncertainty. In this model, the decision-maker must determine the stock level in advance to balance the trade-off between underage cost, arising from insufficient inventory that leads to lost sales or customer dissatisfaction, and overage cost, resulting from excess inventory that incurs holding or disposal costs. This simple yet powerful formulation has been widely applied in various contexts, including retail inventory management for perishable goods and healthcare staffing decisions. While the classic model assumes that the demand distribution is known, real-world applications often face significant uncertainty about this distribution. To address this challenge, recent research has focused on the \textit{data-driven newsvendor model}, which leverages historical demand data to guide inventory decisions without requiring full knowledge of the underlying demand distribution.

A common assumption in most existing literature of  {data-driven newsvendor problems} is stationarity:   demand data \(D_1, \cdots, D_T\) are drawn from the \textit{identical} distribution. 
However, nonstationarity is an intrinsic feature of the real-world operating environment. 
Stationarity, although analytically convenient, may fail to capture important temporal dynamics, resulting in outdated and suboptimal decisions, as well as substantial operational inefficiencies. In practice,
nonstationarity can generally be categorized into two types:  abrupt changes and gradual changes. Abrupt changes may result from sudden events such as regulatory shifts or viral marketing campaigns, whereas gradual changes can arise from seasonal trends, market evolution, or long-term economic shifts. {
For instance, in perishable retail (e.g., supermarket fresh produce), seasonal and climatic trends gradually change the underlying demand pattern over time, and a flash promotion can trigger sudden demand shocks. Similarly, the demand for emergency departments also exhibits nonstationarity.}
{Besides, the timing and degree of these changes are typically unknown a priori.
Thus, rather than committing to a fixed adjustment strategy, we should adaptively adjust the data-driven policies online according to real-time environments to mitigate the adverse effects of nonstationarity.}

A further layer of complexity in data-driven newsvendor problems arises from demand censoring, a pervasive phenomenon in practice. In many settings, when inventories are insufficient, actual demand remains unknown beyond recorded sales, and only censored demand is observable, such as retail stockouts, limited-capacity flash sales, or perishable inventory systems. 
{In nonstationary environments,  demand censoring further obscures the true demand pattern, exacerbating the difficulty of demand learning and inventory decision-making. }

Motivated by the above practical needs, this paper makes several contributions. We introduce a direct and distribution-level characterization of nonstationarity, operating in a fully nonparametric demand model. Based on this characterization, we propose a flexible algorithm design framework tailored to nonstationary environments and develop efficient algorithms for both censored and uncensored demand settings. 
Moreover, we provide a comprehensive regret analysis that establishes the optimal dynamic regret of the proposed algorithms, under various nonstationarity measures and across different distribution classes for both uncensored and censored demand settings.

\subsection{Our Contributions}
In the following,  we present a more detailed discussion of our main contributions.

\textbf{Contribution I: Fully nonparametric modeling and direct characterization of nonstationarity.}
We introduce a direct, distribution-level characterization of nonstationarity within a fully nonparametric demand model. Existing work on nonstationary newsvendor problems often adopts a parametric perspective, for example, assuming a drifting parameter \(\{\mu_t\}_{t=1}^T\) plus i.i.d.\ noise. In contrast, we impose no parametric assumptions on the demand distributions. We treat the observed sequence \(\{D_t\}_{t=1}^T\)  as draws from an arbitrary sequence of distributions and model nonstationarity directly as distributional change. This approach eliminates restrictive structural assumptions on the noise or the mechanism through which nonstationarity enters the model.

We propose two distribution-level measures of nonstationarity on the sequence of distributions that capture both abrupt and gradual changes commonly encountered in real-world settings:
\begin{enumerate}
    \item The \emph{switch number budget}:
      \begin{equation*}
    \sum_{t=2}^T \mathbb{I}\bigl[D_t \overset{d}{\neq} D_{t-1} \bigr]  \leq S,
    \end{equation*}
    which bounds the number of distribution switches.\footnote{$D_t \overset{d}{\neq} D_{t-1} $ refers to that $D_t$ and $D_{t-1}$ have different distributions.  $\mathbb I[\cdot]$ is the indicator function.}
    \item The \emph{total variation budget}:
    \begin{equation*}
    \sum_{t=2}^T \|D_t - D_{t-1}\|_{\mathrm{TV}} \leq V,
    \end{equation*}
    which controls the total variation of the distributional shift.\footnote{The total variation between two distributions, defined as \(\|D_t-D_{t-1}\|_{\mathrm{TV}} = \frac{1}{2} \int |dD_t - dD_{t-1}|\), is one of the most common metrics for measuring the distance between distributions.}
\end{enumerate}
The variation budget \(V\)  is well suited to environments with gradual change, while the switch budget \(S\)  is better for settings with a small number of abrupt changes (where \(V\) may be large even when \(S\) is small). {In contrast, prior work has largely focused on parameter-level variation in sequences \(\{\mu_t\}_{t=1}^T\) and has not explicitly accounted for the number of distribution switches. Moreover, our distribution-level measures are more general and natural than parameter-level ones, and are capable of capturing a broader class of nonstationary behaviors.}

\textbf{Contribution II: A powerful framework and two efficient algorithms.}
 We propose the new \emph{distributional-detection-and-restart} framework for learning in nonstationary environments without knowledge of the degree of nonstationarity, and instantiate it through two novel algorithms for the nonstationary newsvendor problems.  The core idea of our distributional-detection-and-restart framework is to monitor the stationarity of the distributions and search for the optimal solution based on the data detected to be stationary.
 This monitoring is carried out through a distributional nonstationarity detection layer.  {Unlike prior works that detect changes using function values, our method performs distributional detection by comparing empirical distributions across different time intervals.}  
Guided by this framework, we develop two algorithms tailored to different demand observation settings: uncensored and censored demand. Specifically,
\begin{enumerate}
 \item \textit{For uncensored demand:} Algorithm~\ref{alg:saa-uncensored-nvp} combines the distributional-detection-and-restart mechanism with Sample Average Approximation (SAA). It detects changes in the demand distribution and restarts the SAA solver using only the most recent stationary data segment, ensuring continued adaptation without prior knowledge of the degree of nonstationarity.
 \item \textit{For censored demand:} Algorithm~\ref{alg:saacensored} integrates a similar detection layer with a specially designed elimination-based simulation optimization method. It successfully addresses the challenges due to the partial observability inherent in the censored demand setting.
\end{enumerate}
The algorithms are fully adaptive, requiring no prior knowledge of the degree and type of nonstationarity, and offer a flexible yet powerful approach to handling both abrupt and gradual changes in nonstationary environments. 
 Numerical experiments on two real‑world datasets demonstrate that our methods consistently outperform existing approaches, 
showcasing our algorithms’ superior and robust empirical performance.
{Additionally, while we use the newsvendor problem as our primary motivation, the distributional-detection-and-restart framework we develop is broadly applicable to a general class of nonstationary stochastic optimization problems found throughout operations management and machine learning (see Section~\ref{sec:extension}).}

\textbf{Contribution III: Comprehensive analysis of optimality under different conditions.}
We establish a comprehensive optimality theory for our algorithms by deriving matching regret upper and lower bounds under both general and refined structural conditions.
{Our epoch-based analysis involves a carefully designed distributional-detection layer that balances detection power against restart frequency; the former controls the regret incurred within an epoch, while the latter determines the number of epochs. }
 This analytic framework is novel and broadly applicable, combining nontrivial proof techniques that should be of independent interest for other nonstationary decision‑making problems.

Our results demonstrate that our algorithms attain minimax optimal dynamic regret in $T$, $S$, and $V$  across different distribution classes, for both censored and uncensored demand settings (see Table~\ref{tab:theoretical guarantee of our algorithms} for a summary). This comprehensive characterization covers the full spectrum of nonstationary environments, from abruptly changing to gradually drifting environments.
First, we prove the regret upper bound for general demand distributions. We then improve the analysis for demand distributions that satisfy the global minimal separation condition\footnote{The condition refers to the fact that the demand density is lower bounded by a positive constant globally. See Definition~\ref{def:global-separation} for a detailed definition.}, which implies strong convexity of the newsvendor cost function.
Furthermore, we extend this analysis to objective functions satisfying the Polyak-Lojasiewicz (PL) condition, encompassing a class of nonconvex functions in Operations Management (OM) contexts (see Section~\ref{sec:extension}).
For the lower bound, we construct new hard instances specifically designed for the newsvendor problem that capture the intrinsic difficulty of learning in nonstationary environments.

\begin{table}[h!]
\caption{Theoretical guarantees for our algorithms under both uncensored and censored demand settings.} 
\label{tab:theoretical guarantee of our algorithms}
\renewcommand{\arraystretch}{1.5} 
\begin{center}
\scalebox{1}{%
\begin{tabular}{l|c|c|c|c}
\hline
     Setting   & \makecell[l]{Upper Bound:$V$} & \makecell[l]{Lower Bound:$V$} & \makecell[l]{Upper Bound:$S$} & \makecell[l]{Lower Bound:$S$}\\ \hline
{General Distribution}     & $\mathcal O(V^{1/3}T^{2/3}\sqrt{\log T})$ & $ \Omega(V^{1/3}T^{2/3})$     & $\mathcal O(\sqrt{ST\log T})$ & $ \Omega(\sqrt{ST})$  \\ \hline
\makecell[l] {Global Minimal\\ Separation Condition}     & $\mathcal O(V^{2/3}T^{1/3}\log T)$ & $ \Omega(V^{2/3}T^{1/3})$  & $\mathcal O(S\log(T/S)\log T)$ & $ \Omega(S\log(T/S))$    \\ \hline
\end{tabular}%
}
\end{center}
\end{table}

\subsection{Related Works}
In this section, we will carefully discuss three streams of work that are related to our work.

\textbf{Nonstationary Newsvendor Problems.} 
Two recent papers studied nonstationary newsvendor problems and provide regret analyses \citep{keskin2023nonstationary,an2025nonstationary}.  \cite{keskin2023nonstationary} model each demand distribution as \(D_t=\mu_t+\epsilon_t\), where the mean \(\mu_t\) drifts over time and the noise \(\epsilon_t\) is i.i.d., continuous, and bounded. They adopted a variation measure on the sequence of means \(\boldsymbol{\mu}=(\mu_1,\dots,\mu_T)\), defined as $V_{\boldsymbol{\mu}}=\max_{0\le K\le T-1}\ \max_{(t_0,\dots,t_K)\in\mathcal{P}}\ \sum_{k=1}^K \bigl|\mu_{t_k}-\mu_{t_{k-1}}\bigr|^2,$
where \(\mathcal{P}\) denotes the set of all partitions. They designed policies that attain optimal regret when the degree of nonstationarity is known a priori (for example, to schedule restarts).  \cite{an2025nonstationary} considered a family of sub-Gaussian parametric distributions \(\mathcal{D}=\{D_\mu:\mu\in[\mu_{\min},\mu_{\max}]\}\), adopted the same nonstationarity measure \(V_{\boldsymbol{\mu}}\), and proposed a self-adaptive fixed-window procedure that does not require prior knowledge of the degree of nonstationarity.

Our work makes a significant step for nonstationary newsvendor problems beyond these lines in three key ways. 
{First, the previous studies take a parametric view of nonstationarity, where nonstationarity is represented through a drifting parameter sequence $\{\mu_t\}_{t=1}^T$. Such a formulation imposes strong structural assumptions on how demand evolves over time. In contrast, we make no parametric assumptions on the demand process. We treat the observed demand sequence as samples drawn from an arbitrary sequence of distributions and model nonstationarity directly at the distributional level. This nonparametric formulation avoids restrictive assumptions on the noise structure or the mechanism through which nonstationarity enters the system, thereby offering a more faithful description of real-world demand patterns that are often irregular and complex. From a practical perspective, the ability to flexibly capture distributional shifts is crucial, since real-world demand rarely conforms to simple parametric dynamics. By adopting a nonparametric framework, our approach adapts more effectively to the data and, as a consequence, achieves superior empirical performance. In particular, across two real-world datasets, our algorithm consistently outperforms both the parametric algorithm of \cite{keskin2023nonstationary} and a range of heuristic benchmarks. These results provide strong evidence that nonparametric modeling of nonstationarity is not only theoretically appealing but also practically necessary for achieving robust performance in data-driven inventory control.}

Second, we introduce two distribution-level measures of nonstationarity: a switch budget and a variation budget. This comprehensive characterization covers the full spectrum of nonstationary environments, from abruptly changing to gradually drifting environments. Prior works focused on parameter-level variation of the sequence ${\mu_t}$ and did not explicitly account for the number of switches. {Third, demand censoring is pervasive in practice, but remains unaddressed by existing algorithms in prior works. The primary difficulty arises because it obscures the true demand pattern, thereby exacerbating the challenges in demand learning. Moreover, it significantly complicates reliable parameter estimation, rendering methods that rely on accurate estimates ineffective.} In contrast, our proposed algorithmic framework successfully overcomes the challenges induced by demand censoring. For a summary of comparisons, please refer to Table~\ref{tab:comparison in inventory}.

  \begin{table}[h!]
\caption{Comparisons with \cite{keskin2023nonstationary,an2025nonstationary}.} 
\label{tab:comparison in inventory}
\renewcommand{\arraystretch}{1.5} 
\begin{center}
\scalebox{0.9}{%
\begin{tabular}{l|c|c|c|c|c}
\hline
       & \makecell[l]{Nonstationarity \\ characterization} & \makecell[l]{Total Variation\\ Budget} & \makecell[l]{Switch Number\\ Budget} & \makecell[l]{Assumption on\\ Known Budget} & \makecell[l]{Censored Demand}   
       \\ \hline
\cite{keskin2023nonstationary}   & Parametric\tablefootnote{``Parametric'' refers to that they define nonstationarity based on parameters $\{\mu_t\}_{t=1}^T$. } & \ding{51} & \ding{55}     & Required    &  \ding{55} 
 \\ \hline
\cite{an2025nonstationary}   & Parametric  & \ding{51} & \ding{55}     & Not Required    &  \ding{55} 
\\ \hline
Our Algorithm~\ref{alg:saa-uncensored-nvp}     & \textbf{Nonparametric}   & \ding{51} & \ding{51}    & \textbf{Not Required}    &  \ding{55}  \\ \hline
Our Algorithm~\ref{alg:saacensored}  & \textbf{Nonparametric}     & \ding{51} & \ding{51}    & \textbf{Not Required}    &  \ding{51}  \\ \hline
\end{tabular}%
}
\end{center}
\end{table}

\textbf{Other Related Work on Data-driven Inventory Management.} 
There is a stream of work applying SAA in newsvendor problems focusing on sample complexity  \cite{levi2007nearoptimal, levi2015data,besbes2023big,chen2024survey} or cumulative regret  \cite{besbes2013implications,lin22datadriven}.  Recently,  \cite{lyu2024closing} closed the regret gaps and established
the optimality of SAA assuming both global and local strong convexity of the cost objective function, 
In addition to newsvendor problems, sampling-based methods have been applied to other inventory management problems, see e.g. \cite{cheung2019sampling,zhang2021sampling,qin2022data,qin2023sailing,wang2023hybrid}.
There is also a stream of work focusing on the gradient-descent method for data-driven inventory management problems, see e.g. \cite{huh2009nonparametric,huh2009adaptive,shi2016nonparametric,zhang2018perishable,zhang2020closing,yuan2021marrying,lyu2023minibatch}.
\cite{lyu2023minibatch, chen2023learning,chen2024learning} proposed bandit algorithms for inventory management.
Additionally, \cite{chu2023solving,feng2025contextual} studied the newsvendor problem based on the operational data analytics framework introduced in \cite{feng2023framework}.

\textbf{Nonstationary Stochastic Optimization.} We briefly review studies about nonstationary stochastic optimization and discuss the uniqueness of our work from two perspectives. 

\begin{itemize}
\item \textit{Types of Feedback.} Existing studies often employ a \textit{black-box} approach, where the algorithm observes limited feedback at decision points. Common examples include deterministic first-order feedback \citep{mokhtari2016online}, noisy first-order feedback \citep{besbes2015non,cao2020online,kim2022online}, and noisy zero-order feedback \citep{flaxman2004online,chen2019nonstationary}. However, in the newsvendor problem and other operations management tasks, the cost function is predefined. For instance, inventory management problems use a piecewise linear cost function.

This paper adopts a \textit{white-box} perspective, where the inner-layer loss function is explicitly known, and random variable realizations can be observed. 
This approach enables the extraction of additional information specific to the task, thereby enhancing overall performance. {Therefore, our work contributes to this literature by proposing a fully adaptive methodology tailored to ``white-box'' problems where direct distributional detection offers a more powerful and natural approach than the function-value or gradient-based feedback methods common in ``black-box'' settings.}

\item \textit{Measure of Nonstationarity.} Different studies adopt distinct approaches to characterize nonstationarity, tailored to their contexts or methods (see \cite{pun2024online} for a comprehensive review).  \cite{zinkevich2003online} introduced the cumulative path variation of the optimal solution in the context of online gradient descent to capture dynamic regret, a definition widely used in subsequent works \citep{besbes2015non,mokhtari2016online}. \cite{chen2019nonstationary} proposed a generalized $L_{p,q}$ variation of path variation. Recently, \cite{huang2023stability}  proposed an adaptive moving-window algorithm for statistical learning in nonstationary environments. They established the optimal regret with respect to $T$ and the path variation.
{Unlike their method, which detects changes using function‑value tests, our approach performs distributional detection by comparing empirical distributions across different time intervals.} 
In our paper, we model the nonstationarity of the environment by the switch and variation of distributions directly, which aligns well with the white-box model. 
\end{itemize}

\textbf{Nonstationary Reinforcement Learning.}
A series of works investigates reinforcement learning in the nonstationary setting, including the special cases of the multi-armed bandit (MAB) and linear bandit.
Most of these works assume the degree of nonstationarity is known in advance \citep{besbes2019optimal,cheung2022hedging,russac2019weighted,pmlr-v108-zhao20a,mao2024model,cheung2023nonstationary}.
\cite{besbes2019optimal} studied MAB with nonstationary reward. Their algorithm restarts the EXP3
algorithm \citep{auer2002nonstochastic} by a fixed cycle.  
\cite{pmlr-v108-zhao20a}  generalized the restarting algorithm to the linear bandit problem. \cite{cheung2022hedging} and \cite{russac2019weighted} proposed moving-window and exponential weighting types of learning algorithms for the linear bandit problem, respectively.  In these works, to obtain the optimal regret, the cycle length, window size, and weighting parameter should be fine-tuned with the knowledge of the degree of nonstationarity. \cite{cheung2022hedging,cheung2023nonstationary}  proposed a bandit-over-reinforcement-learning technique to handle the unknown degree of nonstationarity setting, but this technique leads to suboptimal regret.

Another stream of works relaxed the known degree of nonstationarity assumption and still achieved the optimal regret. 
\cite{auer2019adaptively} proposed an adaptive restarting algorithm for MAB and established the optimal regret for the switch number of parameter changes. 
\cite{wei2021non} proposed a sophisticated adaptive restarting framework for general reinforcement learning. They demonstrated that any algorithm satisfying specific properties in a weakly nonstationary environment can be combined with their framework to achieve optimal regret without knowledge of the degree of nonstationarity. 
\cite{wang2023adaptivity} designed a bandit optimization algorithm that satisfies the sufficient properties for applying the adaptive restarting framework in \cite{wei2021non}.
{The key distinction between our approach and the framework in \cite{wei2021non} is that the latter adopts a black-box setting and detects changes via function-value tests, whereas we adopt a white-box setting and perform detection at the distribution level. Our distribution-detection-and-restart framework enables us to develop efficient algorithms, derive improved regret bounds, and address challenges specific to operations management problems.}

\subsection{Organization}
The remainder of the paper is organized as follows.
In Section~\ref{sec: Problem Formulation}, we introduce the nonstationary newsvendor model, notation, admissible policies, the two observation models (uncensored and censored), and the two nonstationarity measures $S$ and $V$ used throughout the paper.
Section~\ref{sec:uncensored} studies the uncensored-demand setting: we present our algorithm with uncensored demand and establish matching upper and lower bounds on its dynamic regret.
In Section~\ref{sec:censored}, we extend the methodology to the censored-demand case, describes the necessary estimator and algorithmic adaptations, and states the corresponding performance guarantees.
In Section~\ref{sec:extension}, we develop a general framework for nonstationary stochastic optimization and present several other applications in operations management.
In Section~\ref{sec:Numerical Experiments}, we present empirical evaluations on two real-world datasets, and compare our algorithm with baseline methods.
In Section~\ref{sec: Conslusion}, we conclude the paper and outlines directions for future work.
All technical proofs, auxiliary lemmas, and supplementary details are deferred to the appendix.

\section{Problem Formulation}\label{sec: Problem Formulation}
We consider the nonstationary newsvendor problem over $T$ periods.
In each period $t\in[T]$, the manager needs to decide order quantity $x_t$ from set $\mathcal{X}$ to minimize the following newsvendor cost
\[ 
\min_{x\in{\mathcal X}}f_t(x) := \mathbb{E}_{D_t \sim G_t}[F(x,D_t)],
\]
where $F(x,d) = h(x-d)^+ + b(d-x)^+$, $h$ and $b$ are unit overage cost and underage cost, respectively. For each $t \in [T]$, $D_t$ is a random demand 
with Cumulative Distribution Function (CDF) $G_t(y)$.

A policy $\pi=(\pi_t)_{t \geq 1}$ is called admissible if for any $t$, the mapping $\pi_t:\mh_t\mapsto x_t $ is measurable, where $\mathcal{H}_t$ is the historical data at the beginning of period $t$. In this paper, we will consider the following two different demand observation settings:
\begin{itemize}
    \item \textbf{Uncensored demand.} After we make the decision $x_t$,  we can observe $D_t$, and $\mathcal{H}_t:= \{D_i \}_{i=1}^{t-1} $.

    \item \textbf{Censored demand.} With order quantity $x_t$, we can only observe the sales data $\hat{D}_t=\min\{x_t, D_t\}$. Therefore, in this case, we have $\mathcal{H}_t:= \{\hat{D}_i \}_{i=1}^{t-1} $.
\end{itemize}
The objective is to design a policy $\pi$ that minimizes the following  dynamic regret 
\begin{equation}
    \label{eq:regret}
    \mathcal{R}^{\pi}(T) = \sum_{t=1}^T  f_t(x^{\pi}_t)  - \sum_{t=1}^T  f_t(x_t^*),
\end{equation}
where $x^{\pi}_t$ is the decision generated by policy $\pi$ at time $t$ and $x_t^*=\argmin_{x\in \mathcal{X}} f_t(x)$.

\begin{remark}[Uniqueness of the Optimal Solution]
Without loss of generality, we assume that the optimal solution to the problem $\argmin_{x\in \mathcal{X}} f_t(x)$ is unique. This is because the main objective of our algorithm design and analysis is to control the performance gap between our solution and the optimal one. Therefore, even if there are multiple optimal solutions, they all yield the same optimal value, and it suffices for us to compare with just one of them.
\end{remark}

Since the optimal solution $x_t^*, t \in [T]$ can change arbitrarily, no algorithm can achieve a sublinear regret if we do not impose any assumptions on the nonstationarity of the random variables $
\{D_t\}_{t=1}^T$.
As mentioned in the introduction, we propose two distribution-level measures of nonstationarity on the sequence of distributions that capture both abrupt and gradual changes as follows.
\begin{enumerate}
\item Upper bound $S$ on the number of distribution switches, i.e., 
\begin{equation}
    \label{def:S}
\sum_{t=2}^{T} \mathbb{I}\bigl[D_t \overset{d}{\neq} D_{t-1} \bigr]  \le S.
\end{equation}
    \item Upper bound $V$ on the total variation of the distributional shift, i.e.,
\begin{equation}
\label{def:V}
\sum_{t=2}^{T}\|D_{t}-D_{t-1}\|_{\mathrm{TV}}\le V.
\end{equation}

\end{enumerate}

Our goal is to design algorithms and provide theoretical guarantees for their dynamic regret, i.e., optimality in terms of the time horizon $T$ and these two nonstationarity measures.

For convenience of regret analysis, given a time interval $\mathcal I:=[t_1,t_2]$, we define $$S_{\mathcal I}:=\sum_{t=t_1+1}^{t_2}\mathbb{I}\bigl[D_t \overset{d}{\neq} D_{t-1} \bigr]\text{ 
 ~~and~~  } V_{\mathcal I}:=\sum_{t=t_1+1}^{t_2}\|D_{t}-D_{t-1}\|_{\mathrm{TV}}.$$ 

Throughout this paper, we assume that $\mathcal{X}=[0,\bar{x}] \subset\mathbb R_+$. We also require that random variables $
\{D_t \mid t \in [T]\}$ are bounded\footnote{The bounded assumption on demand is very common in the data-driven inventory management literature; see, e.g., \cite{besbes2013implications,keskin2023nonstationary}.} and supported on $[0,\bar{x}]$ and independent over time.

\section{Uncensored Demand: Algorithm Design and Analysis}\label{sec:uncensored}
In this section, we present the algorithm design and analysis for the uncensored newsvendor problem. We first introduce the algorithm in Section \ref{subsec:Algorithm Design-uncensored}. Then we present upper bounds for the algorithm in Section \ref{subsec:uncensored-ub}, and provide the matching lower bounds in Section \ref{subsec:uncensored-lb}.
\subsection{Algorithm Design}
\label{subsec:Algorithm Design-uncensored}

Before diving into the details of the algorithm for uncensored newsvendor problems, we first introduce the high-level idea of our novel design framework.

\textbf{Distributional‑detection‑and‑restart framework.} We propose a novel design framework, which is called  \textit{distributional‑detection‑and‑restart}, for learning in nonstationary environments without prior knowledge of the degree of nonstationarity.
The core idea is to continuously monitor whether the demand distribution is stationary and, when stationarity is detected, search for the optimal solution using only the data from the most recent stationary segment. Monitoring is implemented via a distributional nonstationarity detection layer that compares empirical distributions over different time intervals.
This framework offers an intuitive and flexible design.
Guided by this framework, we develop two algorithms tailored to different observation regimes: uncensored and censored demand. In the following, we present the algorithm for the uncensored demand. 

For the uncensored setting, Algorithm~\ref{alg:saa-uncensored-nvp} combines a distributional nonstationarity detection layer into SAA. This layer continuously monitors historical sequential data, and upon detecting nonstationarity, discards outdated data. The algorithm then restarts the SAA process with newly collected data, ensuring adaptability and precision. Thus, we refer to our approach as \textit{Nonstationary SAA with Adaptive Restarts}, NSAA for short.

\begin{algorithm}[h!]
\caption{Nonstationary SAA with Adaptive Restarts for Uncensored Newsvendor Problem}
\label{alg:saa-uncensored-nvp}
\begin{algorithmic}[1]
\State \textbf{Input:} Confidence level $\delta\in(0,1)$.
\State \textbf{Initialization:} Set the epoch number $\tau=1$, and initial time step $t=1$, $l_1=1$.
\For{$\tau =1,2,\cdots$}
\For{$t=l_\tau,l_\tau+1,\cdots$}\label{Line:for}
\State \textbf{If} $t = l_{\tau}$ \textbf{then} set $x_t$ randomly from $\mathcal{X}$.  
\State \textbf{If} $t > l_{\tau}$ \textbf{then} optimize the following problem to obtain the decision $x_t$:
\[
\min_{x\in\mathcal X}\frac{1}{t-l_\tau} \sum_{k=l_\tau}^{t-1} F(x,D_k).
\]
\State  \label{Line:detect} \textbf{If} there exists $s\in [l_\tau,t]$ such that
\[
\max_{y\in \R} \vert\hat{G}_{l_\tau,t-1}(y) - \hat G_{s,t}(y) \vert > 2\sqrt{\frac{\ln (2T^2/\delta)}{t-l_\tau}}+ 2\sqrt{\frac{\ln (2T^2/\delta)}{t-s+1}},
\]
\State \textbf{then} set $\tau=\tau+1$ and $l_\tau = t+1$, and start a new epoch at Line \ref{Line:for}.
\EndFor
\State \textbf{end for}
\EndFor
\State \textbf{end for}
\end{algorithmic}
\end{algorithm}
Our algorithm operates in epochs indexed by $\tau$. After making the decision $x_t$ and observing demand $D_t$, the algorithm checks whether the distribution remains stationary (see Line \ref{Line:detect}). To perform this check of nonstationarity, we construct the empirical CDF based on the data collected from periods $s$ to $t$ as 
\[
\hat{G}_{s,t}(y) = \frac{1}{t-s+1}\sum_{i=s}^t \mathbb{I}[D_{i} \leq y],
\]  
and derive a confidence region for the true CDF, $G_{s,t}(y) = 1/(t-s+1)\sum_{i=s}^t\E[\mathbb{I}[D_{i} \leq y]]$. If there exists a time $s$ such that $[s, t]$, where the confidence region of $\hat G_{s,t}(y)$ is disjoint from that of $\hat{G}_{l_\tau,t-1}(y)$, it indicates a change in the distribution during certain periods. 
Specifically, using a uniform concentration bound for empirical CDFs, we have with high probability
\[
\sup_{y\in\mathbb{R}}|\hat G_{l_\tau,t-1}(y)-G_{l_\tau,t-1}(y)| \le \sqrt{\frac{\ln(2T^2/\delta)}{t-l_\tau}}\text{~~and~~} \sup_{y\in\mathbb{R}}|\hat G_{s,t}(y)-G_{s,t}(y)| \le \sqrt{\frac{\ln(2T^2/\delta)}{t-s+1}}.
\]
Therefore, if
\[
\sup_{y\in\mathbb{R}}|\hat G_{l_\tau,t-1}(y)-\hat G_{s,t}(y)|
>2\sqrt{\frac{\ln(2T^2/\delta)}{t-l_\tau+1}}+2\sqrt{\frac{\ln(2T^2/\delta)}{t-s+1}},
\]
we can conclude that a distributional change has occurred on the interval considered.
Note that the supremum over $y$ only needs to be checked at the observed sample values $y=D_1,\dots,D_t$. Upon detecting such a change, the algorithm initiates a new epoch, $\tau+1$, and sets the starting period for the new epoch as $l_{\tau+1} = t+1$.

\subsection{Regret Analysis: Upper Bound}\label{subsec:uncensored-ub}
In this section, we give upper bound results on the performance of Algorithm \ref{alg:saa-uncensored-nvp}. The distributional-detection-and-restart procedures of Algorithm \ref{alg:saa-uncensored-nvp} can be generalized to a broader class of stochastic optimization problems. Therefore, in Section \ref{sec:extension}, we extend the algorithm design and analysis under a general framework and defer the proof in this section as a special case of the general theory.

\begin{theorem}\label{thm:uncensored-nvp-ub-general}
With probability at least $1-\delta$, the following conclusions hold.
\begin{enumerate}
    \item \textit{Switch number $S$:} the dynamic regret of Algorithm~\ref{alg:saa-uncensored-nvp} can be upper bounded by 
    \[ 
    \mcR^{\pi}(T) \leq 40 (h+b)\bar{x} \sqrt{(S+1)T\ln (2 T^2/\delta)} +2(h+b)\bar{x} (S+1).
    \]
    \item \textit{Total variation $V$:} the dynamic regret of Algorithm~\ref{alg:saa-uncensored-nvp} can be upper bounded by 
\[ 
\mcR^{\pi}(T) \leq
    4(h+b)\bar{x}(V^{2/3}T^{1/3}+1) + 96(h+b)\bar{x} (V^{1/3}+1)T^{2/3}\ln^{1/2}(2T^2/\delta).
\]
\end{enumerate}
Therefore, we have $\E[\mcR^{\pi}(T)] = \O(V^{1/3}T^{2/3}\log^{1/2}T)$ and $\E[\mcR^{\pi}(T)] = \O(\sqrt{ST\log T})$.
\end{theorem}

The core idea behind the proof of the above theorems is to bound the regret according to epochs. Since the data is detected to be near-stationary within an epoch, the detection step (Line \ref{Line:detect}) ensures that we have a good estimation of $f_t(x)$, thereby providing an upper bound on $f_t(x_t)-f_t(x_t^*)$. Based on the restarting conditions and the definition of nonstationary measures $S$,$V$, we can prove an upper bound on the total number of epochs. Combining this with the conclusions for each epoch, we can derive a bound for the final regret. Please refer to Section \ref{sec:extension} for the general theory and Section \ref{appendix:extension-proof} for its detailed proof.

\begin{remark}[Time-varying Inner-layer Functions]
Note that we assume that the inner layer function $F(x,d) = h(x-d)^+ +b(d-x)^+$ is invariant over time. While our analysis also works for the problem with time-varying inner-layer functions $F_t(\cdot,\cdot) ,t \in [T]$ under slight modification. For example, the per-unit underage cost $h$ and per-unit overage cost $b$ may vary over time.
\end{remark}

\textbf{Improved regret bounds.} Similar to the stationary newsvendor problem, the performance of Algorithm \ref{alg:saa-uncensored-nvp} can be improved for demand distributions satisfying the global minimal separation condition \citep{lin22datadriven,lyu2024closing}. We introduce the condition as follows.

\begin{definition}[Global minimal separation condition]\label{def:global-separation} A random variable $D$ satisfies the $\alpha$-global minimal separation condition, if its CDF $G_D(x)$ is differentiable and $G_D'(x)$ is bounded from below by $\alpha$ for all $x\in[0,\bar x]$, i.e., $G'_D(x) \geq \alpha$, for all $x \in [0,\bar x]$.
\end{definition}

For the newsvendor problem with cost function $f(x)$, we know that $f''(x) = (h+b)F'_D(x)$. Therefore, demand $D$ satisfies the $\alpha$-global minimal separation condition implies that $f(x)$ is strongly convex. Similar to the general demand case, we can extend the analysis to a general problem with strong convexity, and defer the proof in this section as a special case of the general theory.

\begin{theorem}\label{thm:uncensored-nvp-ub-pl} Suppose that for any $t \in [T]$, the random variable $D_t$ satisfies the $\alpha$-global minimal separation condition.
With probability at least $1-\delta$, the following conclusions hold.
\begin{enumerate}
    \item \textit{Switch number $S$:} the dynamic regret of Algorithm~\ref{alg:saa-uncensored-nvp} can be upper bounded by 
        \[
    \mathcal{R}^{\pi}(T)\leq {52(h+b)}{\alpha^{-1}}(S+1)\ln (2T^2/\delta)(1+\ln (T/S)) + 2(h+b)(1+\alpha^{-1})\bar{x}(S+1).
    \]
    \item \textit{Total variation $V$:} the dynamic regret of Algorithm~\ref{alg:saa-uncensored-nvp} can be upper bounded by 
\[ \mathcal{R}^{\pi}(T) \leq 4(h+b)(1+\alpha^{-1})\bar{x}(V^{2/3}T^{1/3}+1) + 192(h+b)\alpha^{-1} (V^{2/3}T^{1/3}+1)\ln T\ln^{1/2}(2T^2/\delta).
    \]
\end{enumerate}
Therefore, we have $\E[\mcR^{\pi}(T)] = \O(S\log(T/S)\log T )$ and $\E[\mcR^{\pi}(T)] = \O(V^{2/3}T^{1/3}\log^{3/2} T )$.
\end{theorem}

The proof of the above theorem is similar to the general distribution case. The main difference is that we use the condition of the nonstationarity test to provide an upper bound on $\Vert g_t(x_t)-g_t(x_t^*) \Vert$. Based on an inequality implied by the strong convexity ($f_t(x_t)-f_t(x^*_t) \leq \mathcal{O}(\Vert g_t(x_t)-g_t(x_t^*) \Vert^2)$, we can achieve an improved  regret bound compared with the conclusions for general functions.   Please refer to Section~\ref{sec:extension} for the general theory and Section \ref{subsec:Proof of Theorem thm:uncensored-nvp-ub-pl} for its detailed proof.

\subsection{Regret Analysis: Lower Bound}
\label{subsec:uncensored-lb}

In this section, we establish lower bound results. For simplicity, we consider the newsvendor problem with $h=b=1$, and the idea of hard instance design can be easily generalized to general $h$ and $b$.

We denote an instance of the nonstationary newsvendor problem by its $T$-period demand vector $\bm{D}:=\{D_1,~D_2,\dots,D_T\}$. Let $\Pi$ be the policy class that includes all admissible ordering policies. 
 Since the lower bound in the expectation setting implies the lower bound result in the high probability setting, in this paper we establish the lower bound result in the expectation setting.
For convenience, we  define the $T$-period total  dynamic regret under demand distributions $\bm{D}$ and policy $\pi$ as 
\[\mathcal R^{\pi}(T, \bm{D}) = \sum_{t=1}^T  f_t(x_t)  - \sum_{t=1}^T  f_t(x_t^*). \]

Compared with the definition of dynamic regret $\mathcal{R}^{\pi}(T)$ given by Eq.~\eqref{eq:regret}, the above notation $\mathcal{R}^{\pi}(T,\bm{D})$ emphsize the dependence on distribution vector $\bm{D}$.

\textbf{Lower bound for general demand distributions.} Let $\mathcal D_{S}^{(1)}$ be the instance class including all continuous demand distributions $\bD=(D_1, D_2, \dots, D_T)$ that satisfy
$\sum_{t=2}^{T}\mathbb{I}\bigl[D_t \overset{d}{\neq} D_{t-1} \bigr]\le S$ and $\mathcal D_{V}^{(1)}$ be the instance class including all continuous demand distributions $\bD=(D_1, D_2, \dots, D_T)$ that satisfy
$\sum_{t=2}^{T}\|D_{t}-D_{t-1}\|_{\mathrm{TV}}\le V$.  The following lower bound results for the instance set $\mathcal{D}_S^{(1)}$ and $\mathcal D_{V}^{(1)}$ matches the upper bound results in Theorem~\ref{thm:uncensored-nvp-ub-general}, respectively. These results demonstrate the optimality of our algorithm with respect to both $T$ and $S$ up to a logarithmic factor.

\begin{theorem}\label{thm:lb_cxv} The following conclusions hold.
\begin{enumerate}
    \item For instance class $\mathcal D_{S}^{(1)}$, the minimax regret can be bounded from below as
    \[\inf_{{\pi\in\Pi}}\sup_{\bm{D}\in\mathcal D_{S}^{(1)}}\E[\mathcal R^{\pi}(T, \bm{D})]\geq \frac{\exp(-16)}{4}\sqrt{ST}. \]

    \item For instance class $\mathcal D_{V}^{(1)}$, the minimax regret can be bounded from below as
    \[\inf_{{\pi\in\Pi}}\sup_{\bm{D}\in\mathcal D_{V}^{(1)}}\E[\mathcal R^{\pi}(T, \bm{D})]\geq \frac{\exp(-16)}{16} V^{1/3}T^{2/3}. \]
\end{enumerate}
 
\end{theorem}

To prove the above theorems, we divide the $T$-period horizon into batches of $\Delta_T$ periods. In each batch, the demand follows either distribution $D_a$ or $D_b$ with equal probability, and the demand distributions remain unchanged within each batch. The design of $D_a$ and $D_b$ must: 
i) ensure that the nonstationarity between batches satisfies the degree of nonstationarity ($S$ or $V$), and 
ii) balance the similarity between these two distributions with the cost of misspecification within batches. 
After designing the hard instances, we complete the proofs based on Kullback-Leibler divergence and the Bretagnolle-Huber inequality. The formal proof is presented in Sections \ref{sec:Proof of  Theoremthm:lb_cxv}.

\textbf{Lower bound under global minimal separation condition.} Let $\alpha =1/2$ be a constant. We define $\mathcal D_{S}^{(2)}$ be the instance class, including all continuous demand distributions $\bD=(D_1, D_2, \dots, D_T)$ that satisfy $G_{D_i}(x) \geq \alpha$ for all $x \in [0,\bar{x}]$, $ i \in [T]$ and  $\sum_{t=2}^{T}\mathbb{I}\bigl[D_t \overset{d}{\neq} D_{t-1} \bigr]\le S$, where $G_D(x)$ is the PDF of demand $D$. Similarly, we define $\mathcal D_{V}^{(2)}$ be the instance class, including all continuous demand distributions $\bD=(D_1, D_2, \dots, D_T)$ that satisfy $G_{D_i}(x) \geq\alpha$ for all $x \in [0,\bar{x}]$, $ i \in [T]$ and  $\sum_{t=2}^{T}\|D_{t}-D_{t-1}\|_{\mathrm{TV}}\le V$. We can prove the following theorem.

\begin{theorem}\label{thm:lb-str}
The following conclusions hold.
\begin{enumerate}
\item For instance class $\mathcal D_{S}^{(2)}$, the minimax regret can be bounded from below as
 \[\inf_{{\pi\in\Pi}}\sup_{\bm{D}\in\mathcal D_{S}^{(2)}}\E[\mathcal R^{\pi}(T, \bm{D})]\geq \frac{1}{400\pi^2} S \log (T/8S). \]

    \item For instance class $\mathcal D_{V}^{(2)}$, the minimax regret can be bounded from below as\[\inf_{{\pi\in\Pi}}\sup_{\bm{D}\in\mathcal D_{V}^{(2)}}\E[\mathcal R^{\pi}(T, \bm{D})]\geq \frac{\exp(-16)}{128}V^{2/3}T^{1/3}. \]

\end{enumerate}
\end{theorem}

The segmentation approach used in the proofs of the above theorem is similar to that in Theorems~\ref{thm:lb_cxv}. The main difference lies in the design of the hard instance. Additionally, unlike the argument based on Kullback-Leibler divergence and the Bretagnolle-Huber inequality, we design different instances and employ the van Trees inequality in the proof \citep{lyu2024closing}. For the detailed proofs, please refer to Section \ref{subsec:Proof of Theorem thm:lb-str}. Note that for the convenience of presentation, we fix $\alpha=1/2$ in the construction of the hard instances satisfying the global minimal separation condition. By some simple scaling, our results can be generalized to arbitrarily $\alpha>0$.
\section{Censored Demand: Algorithm Design and Analysis} \label{sec:censored}

Demand censoring is pervasive in real-world contexts, for instance, in retail operations with stockouts, in online marketplaces with flash sales, and in agricultural or perishable goods supply chains. 
In this section, we extend our framework to address the nonstationary newsvendor problem in the presence of {censored demand}, where we only observe sales data $\hat{D}_t = \min \{ x_t, D_t \}$.

\textbf{Unique challenges due to censored demand.}
When demand data are fully observable, as discussed in Section~\ref{sec:uncensored}, we can approximate the cost function using the sample average $\hat{f}_{i,j}(x)$, and obtain the approximated solution $x_t$ by optimizing $\hat{f}_{l_\tau,t}(x)$. Furthermore, we can construct the empirical distribution function of demand $\hat{G}_{i,j}(y)$ for any $y \in \mathbb{R}$, and perform distributional-detection-and-restart.

When demand is censored, we observe only $\{ \hat{D}_k = \min\{x_k, D_k\} \}_{k=1}^{t-1}$ rather than $\{ D_k \}_{k=1}^{t-1}$.
In the censored-demand setting, accurately estimating $f_t(x)$ becomes substantially more challenging because the conventional approach of using historical demand samples is not directly applicable. Specifically, there are two main challenges:
\begin{enumerate}
    \item The SAA estimator $\hat{f}_{i,j}(x)$ cannot be computed directly, since $F(x, D_t)$ depends on the unobserved portion of $D_t$ whenever $D_t > x_t$.
    \item The empirical distribution $\hat{G}_{i,j}(y)$ cannot be constructed exactly, because censored observations do not reveal the upper tail of the demand distribution beyond the ordering level $x_t$.
\end{enumerate}

Thus, the previous standard empirical estimation, decision-making based on SAA, and detection-restart procedures must be adapted to account for censoring. The subsequent section addresses these issues and proposes methodological extensions for censored observations.

\subsection{Algorithm Design for Censored Demand}
This section presents an extension of our framework to handle the challenges arising from censored demand by designing several new subroutines. These enhancements maintain the structure of our detection layer while integrating it into an elimination-based simulation optimization method.

Note that ${G}_{i,j}(x) =  1/(j-i+1) \sum_{k=i}^{j} G_{k}(x)$ and
    $\hat{G}_{i,j}(x) :=  1/(j-i+1)  \sum_{k=i}^{j} \mathbb{I}[D_k \le x]$ For convenience, for $i \le j$, we define
\begin{align*}
    g_{i,j}(x) &:= \frac{1}{j-i+1} \sum_{k=i}^{j} \nabla f_k(x) 
    = (h+b)\,{G}_{i,j}(x) - b, ~~~
    \hat{g}_{i,j}(x) := 
     (h+b) \cdot \hat{G}_{i,j}(x) - b.
\end{align*}
For the special case $i = j$, we abbreviate $g_{i,i}(x)$ as $g_i(x)$ and $\hat{g}_{i,i}(x)$ as $\hat{g}_i(x)$.

Our algorithm for the censored newsvendor problem is presented in Algorithm~\ref{alg:saacensored}. 
At a high level, this extension differs from the uncensored version (Algorithm~\ref{alg:saa-uncensored-nvp}) in three key aspects:
\begin{enumerate}
    \item \textit{Decision rule.} Instead of selecting $x_t$ as the empirical cost minimizer (which cannot be directly computed under censoring), we choose the \emph{largest} available decision in the current active set. 
    \item \textit{Simulation under censoring.} We propose a simulation step by leveraging the censored demand observation $\hat{D}_t = \min\{ D_t, x_t \}$, which enables the estimation of $\hat{g}_{l_{\tau},t}(x)$ for each $x$ in the active set $\mathcal{X}_t$, since the demand is observed under the upper bound $x_t$ of $\mathcal{X}_t$.
    \item \textit{Elimination with censored data.} We design an elimination rule that accounts for censoring, ensuring that suboptimal decisions are effectively removed from the active set. 
\end{enumerate}

\begin{algorithm}[h!]
\caption{Nonstationary SAA-based Elimination  with Adaptive Restarts
}
\label{alg:saacensored}
\begin{algorithmic}[1]
\State \textbf{Initialization:} Input confidence level $\delta\in(0,1)$. Let the epoch number $\tau=1$, and initial time step $t=1$, $l_1=1$. Initialize the active set as $\mathcal X_{l_1}=[0,\bar x]$.
\For{$\tau =1,2,\cdots$}\label{Line:for2}
\For{$t=l_\tau,l_\tau+1,\cdots$}
\State Select the largest level in $\mathcal X_t$ and make the decision $x_t:=\max{\mathcal X_t}$.

\State Observe the censored demand $\hat{D}_t:=\min\{x_t,D_t\}$.
\State \label{line:simulation} For each $x\in\mathcal X_t$, use the censored demand $\hat{D}_t$ to calculate 
\[
\hat{g}_{l_\tau,t}(x)
= \frac{1}{t-l_\tau+1}\sum_{k=l_\tau}^{t}\big[(h+b)\,\mathbb{I}[D_k\le x]-b\big].
\]
\State  \label{line:elimination} Update the active set to $\mathcal X_{t+1}$ 
\[\mathcal X_{t+1}:=\left\{x: x\in \mathcal X_t \text{~and~}  \hat{g}_{l_{\tau},t}(x)\leq 2(h+b) \sqrt{\frac{\ln (2T^2/\delta)}{t-l_{\tau}+1}}\right\}\]

\State  \label{Line:detect2} \textbf{If} either of the following two events occurs:
\begin{enumerate}[leftmargin=6em] 
    \item $\mathcal X_{t+1}=\emptyset$; 
        \item there exists $s\in [l_\tau,t]$, such that 
    \[
\max_{y\in [0,x_t]} \vert\hat{G}_{l_\tau,t-1}(y) - \hat G_{s,t}(y) \vert > 2\sqrt{\frac{\ln (2T^2/\delta)}{t-l_\tau}}+ 2\sqrt{\frac{\ln (2T^2/\delta)}{t-s+1}}.
\]

\end{enumerate}
\State \textbf{then} set $\tau=\tau+1$, $l_\tau = t+1$ and $\mathcal{X}_{l_\tau}=[0,\bar x]$, then start a new epoch at Line \ref{Line:for2}.
\EndFor
\State \textbf{end for}
\EndFor
\State \textbf{end for}
\end{algorithmic}
\end{algorithm}

Algorithm~\ref{alg:saacensored} also proceeds in epochs indexed by $\tau$.
At the start of epoch $\tau$, the active set is initialized to
    $\mathcal{X}_{l_\tau} = [0, \bar{x}].$\footnote{In implementation, we may discretize the initial active set $\mathcal{X}_{l_{\tau}} = [0, \bar{x}]$ into $T$ discrete decisions. This discretization does not affect the order of regret.}
In contrast to Algorithm~\ref{alg:saa-uncensored-nvp}, where $x_t$ is chosen by optimizing the empirical cost estimate, here we set $x_t$ to be the \emph{largest} available decision in $\mathcal{X}_t$. This choice ensures that, although only censored demand $\hat{D}_t$ is observed, we retain enough information to simulate and estimate the gradient approximation $\hat{g}_{l_\tau, t}(x)$ for all $x$ in the active set.
After each simulation step, we update the active set by eliminating any $x \in \mathcal{X}_t$ whose estimated gradient $\hat{g}_{l_\tau, t}(x)$ exceeds a prespecified threshold. This yields the refined active set $\mathcal{X}_{t+1}$.
We also specify two \emph{epoch-restarting conditions}: whenever one of these is triggered, a new epoch $\tau + 1$ is initiated, with the starting time $l_{\tau+1} = t + 1$ and the active set reinitialized to $[0, \bar{x}]$.
This elimination-based, epoch-driven strategy balances effective exploration with robust active set refinement in the nonstationary environment, ensuring convergence despite the loss of full demand observations.

\subsection{New Subroutines}
Below, we provide a detailed description of the new subroutines.

\textbf{Simulation computation.}
For an epoch starting at time $l_\tau$, define the gradient estimate as
\[
\hat{g}_{l_\tau,t}(x)
= \frac{1}{t-l_\tau+1}\sum_{k=l_\tau}^{t}\big[(h+b)\,\mathbb{I}[D_k\le x]-b\big].
\]
Note the above expression requires accurate computation of $\I[D_k \leq x]$. Next,we explain how to obtain $\mathbb{I}[D_t\le x]$ for all $x\in\mathcal{X}_t$ from the single observation $\hat{D}_t$.

Since  $x\le x_t$ for any $x\in\mathcal{X}_t$. The following cases cover all possibilities of $\hat{D}_t=\min\{D_t,x_t\}$:
\begin{itemize}
  \item If $\hat{D}_t < x_t$, then no censoring occurred and we have $\mathbb{I}[D_t\le x]=\mathbb{I}[\hat{D}_t\le x]$, $\forall x\in\mathcal{X}_t$.
  \item If $\hat{D}_t \ge x_t$, then censoring occurred and $D_t\ge x_t \ge x$   for all $x\in\mathcal{X}_t$. We discuss how to estimate $\I[D_t \leq x]$ for continuous and discrete demand separately. 
  \begin{itemize}
      \item For continuous distributions, since \(\P[D_t=x_t]=0\), when censoring occurs we almost surely have \(D_t>x_t\), and hence \(\mathbb{I}[D_t\le x]=0\) for all \(x\in\mathcal{X}_t\). Therefore, we can estimate $\I[D_t \leq x]$ by $\I[\hat D_t \leq x]$, which gives a unbiased estimate of $\P[D_t \leq x_t]$.
      \item For discrete distributions, observing \(\hat{D}_t=x_t\) implies \(D_t\ge x_t\); for every \(x\in\mathcal{X}_t\setminus\{x_t\}\) we then have \(D_t\ge x_t>x\), so \(\mathbb{I}[D_t\le x]=0\) for all \(x\in\mathcal{X}_t\setminus\{x_t\}\). For \(\hat{D}_t=x_t\), we assume an additional \textit{lost-sale indicator}\footnote{Availability of lost-sales indicator is a widely adopted assumption for the learning of discrete products\citep{huh2009nonparametric}, which informs the decision-maker whether a strict stockout happened.} is available and thus \(\mathbb{I}[D_t \leq x_t]\) can be computed accurately. 
  \end{itemize}

\end{itemize}

Note also that the active sets are shrinking over an epoch:
\[
\mathcal{X}_{l_\tau}\supseteq\mathcal{X}_{l_\tau+1}\supseteq\cdots\supseteq \mathcal{X}_{{l_\tau+1}-1}.
\]
Hence the chosen decisions satisfy $x_{l_\tau}\ge x_{l_\tau+1}\ge\cdots\ge x_{{l_\tau+1}-1}$, and historical data collected before period $t-1$ are available for the gradient computation of point $x\in\mathcal{X}_t$. 

\textbf{Elimination of the active set.}
Let $g_t(x)$ denote the true (population) gradient at time $t$; the epoch-specific optimizer $x_t^*$ satisfies $g_t(x_t^*)=0$. Using concentration arguments, one can show that, with high probability,
\[
\big|\hat{g}_{l_\tau,t}(x_t^*)\big| = \tilde{\mathcal{O}}\big((t-l_\tau+1)^{-1/2}\big).
\]
To avoid discarding the optimal candidates prematurely and guarantee the performance of all candidate decisions in the active set, we eliminate only those $x\in\mathcal{X}_t$ whose estimated gradient is significantly large. Specifically, with confidence parameter $\delta$, update the active set by
\[
\mathcal{X}_{t+1}
:= \Big\{x\in\mathcal{X}_t :
\hat{g}_{l_\tau,t}(x)
\le 2(h+b)\sqrt{\frac{\ln(2T^2/\delta)}{t-l_\tau+1}}
\Big\}.
\]

\textbf{Epoch-restarting conditions.}
We use two restart conditions. If either holds at time $t$ we start a new epoch $\tau+1$ with $l_{\tau+1}=t+1$ and reset the active set to be $[0,\bar{x}]$.

\textit{1. Empty active set:} $\mathcal{X}_{t+1}=\varnothing$, which guarantees the next decision $x_{t+1}$ is well-defined.

\textit{2. Nonstationarity detection:} there exists $s\in[l_\tau,t]$  such that
  \[
  \max_{y\in[0,x_t]} \big|\hat{G}_{l_{\tau},t-1}(y)-\hat{G}_{s,t}(y)\big|
  > 2\sqrt{\frac{\ln(2T^2/\delta)}{t-l_{\tau}}}
  + 2\sqrt{\frac{\ln(2T^2/\delta)}{t-s+1}}.
  \]
  This condition mirrors the nonstationarity detection used in Algorithm~\ref{alg:saa-uncensored-nvp}. Although in the censored-demand setting we can only compute both $\hat{G}_{l_{\tau},t-1}(y)$ and $\hat{G}_{s,t}(y)$ for $y\in[0,x_t]$ using the simulation procedure described above; restricting the comparison to $[0,x_t]$ is sufficient for the regret analysis.

 {Similar simulation techniques also have appeared in \cite{yuan2021marrying,chen2024learning}.}
 {Elimination-based method is a classic technique from the bandit literature for maintaining a high-performing active set \citep{even2006action}. Here, however, our elimination criterion is carefully designed for the gradient rather than the cost (or reward) commonly used in prior works. This tailored design enables us to establish the optimal regret using a gradient-based analysis.}
{The design of Algorithm \ref{alg:saacensored} represents a novel synthesis of our distributional-detection-and-restart framework with techniques from simulation optimization and bandit literature. While elimination-based methods and simulation techniques are known tools for exploration in stationary settings, our contribution lies in embedding this mechanism within an adaptive-restarting layer to specifically address the joint challenge of nonstationarity and demand censoring, a problem that, to our knowledge, has not been previously addressed with theoretical regret guarantees.
 Together with Algorithm \ref{alg:saa-uncensored-nvp}, we show the flexibility and effectiveness of our framework.
}

\subsection{Regret Analysis: Upper Bound}

In this section, we give upper bound results on the performance of Algorithm \ref{alg:saacensored} in the censored demand setting, which include the upper bound for general distribution and the improved upper bound for the distributions satisfying the global minimal separation condition. 
All bounds are of the same order as their counterparts in the uncensored setting.

\begin{theorem}\label{thm:censored-nvp-ub-general}
With probability at least $1-\delta$, the following conclusions hold.
\begin{enumerate}
    \item \textit{Switch number $S$:} the dynamic regret of Algorithm~\ref{alg:saacensored} can be upper bounded by 
    \[ 
    \mcR^{\pi}(T) \leq  32 \bar x(h+b) \sqrt{(S+1)T\ln (2 T^2/\delta)}.
    \]
    \item \textit{Total variation $V$:} the dynamic regret of Algorithm~\ref{alg:saacensored} can be upper bounded by 
\[ 
\mcR^{\pi}(T) \leq
   80\bar x(h+b) (V^{1/3}+1)T^{2/3}\ln^{1/2}(2T^2/\delta).
\]
\end{enumerate}
Therefore, we have $\E[\mcR^{\pi}(T)] = \O(V^{1/3}T^{2/3}\log^{1/2}T)$ and $\E[\mcR^{\pi}(T)] = \O(\sqrt{ST\log T})$.
\end{theorem}

The next theorem gives improved upper bounds for demand distributions that satisfy the $\alpha$-global minimal separation condition.

\begin{theorem}\label{thm:censored-nvp-ub-pl} Suppose that for any $t \in [T]$, the random variable $D_t$ satisfies the $\alpha$-global minimal separation condition.
With probability at least $1-\delta$, the following conclusions hold.
\begin{enumerate}
    \item \textit{Switch number $S$:} the dynamic regret of Algorithm~\ref{alg:saacensored} can be upper bounded by 
        \[
    \mathcal{R}^{\pi}(T)\leq \frac{64(h+b)}{\alpha}(S+1)\ln (2T^2/\delta)(1+\ln (T/S)).
    \]
    \item \textit{Total variation $V$:} the dynamic regret of Algorithm~\ref{alg:saacensored} can be upper bounded by 
\[ \mathcal{R}^{\pi}(T) \leq \frac{128(h+b)}{\alpha} (V^{2/3}T^{1/3}+1)\ln T\ln^{1/2}(2T^2/\delta).
    \]
\end{enumerate}
Therefore, we have $\E[\mcR^{\pi}(T)] = \O(S\log(T/S)\log T )$ and $\E[\mcR^{\pi}(T)] = \O(V^{2/3}T^{1/3}\log^{3/2} T )$.
\end{theorem}

The regret upper bounds above match the lower bounds established for the uncensored newsvendor problem in Section~\ref{subsec:uncensored-lb}, demonstrating the robustness and flexibility of our algorithmic design in handling the challenges posed by censored demand observations.

The proof of the above theorem is presented in Section~\ref{sec: Proofs Omitted in sec:censored} in the appendix.
The analysis differs from the uncensored case in two aspects. First, the optimization step in Algorithm~\ref{alg:saa-uncensored-nvp} guarantees the performance of $x_t$ for the empirical function. In contrast, the regret analysis for Algorithm~\ref{alg:saacensored} exploits properties and the design of the active set to ensure a good performance guarantee for $x_t$. Second, because the restarting condition is different in the censored setting, we prove a new technical lemma that bounds the number of epochs produced by Algorithm~\ref{alg:saacensored}, which is essential for the regret guarantees.

\section{Extension: General Framework and Further Applications}\label{sec:extension}
As mentioned earlier, the idea of algorithm design and analysis for the uncensored newsvendor problem can be extended to a general optimization problem.
We consider a nonstationary stochastic optimization problem over $T$ periods as follows
\begin{equation}
\label{eq:general-opt-problem}
    \min_{x\in{\mathcal X}}f_t(x) := \mathbb{E}_{D_t \sim G_t}[F(x,D_t)],
\end{equation}
where $F(x,d)$ is a known function (not necessarily the newsvendor cost) and for each $t \in [T]$, $D_t$ is a random variable 
with CDF $G_t(y)$. Let $\mathcal{H}_t = \{D_i \}_{i=1}^{t-1}$ 
be the history data at the beginning of period $t$. Similar to the uncensored newsvendor problem (Section \ref{sec: Problem Formulation}), we can define an admissible policy, dynamic regret Eq.~\eqref{eq:regret}, and measures of nonstationarity for random variables (number of distribution switches Eq.~\eqref{def:S} and total variation Eq.~\eqref{def:V}). 

In this section, we impose the following basic assumptions, which require $\cX$, $f_t(x)$ to be bounded and random variables $
\{D_t \mid t \in [T]\}$ to be independent over time.

\begin{assumption}[Boundedness and Independence]\label{ass:bounded}
The following conditions hold.
\begin{enumerate}
    \item The domain $\mathcal{X}\subseteq\mathbb R$ is convex and compact with $\mathrm{diam}(\mathcal X):=\max_{x,x'\in\mathcal X}\|x-x'\|_2\leq D_{\mathcal{X}}<\infty$.
    \item There exists a constant $C_f$ such that $|f_t(x)| \leq C_f$ for all $x \in \cX$, $t \in [T]$.
    \item Random variables $D_1,D_2,\dots,D_T$ are independent.
\end{enumerate}
\end{assumption}

Different from the newsvendor problem, in this section, the decision variable is allowed to be $n$-dimensional. The algorithm design and analysis can also be generalized to the case with multi-dimensional random variables $D_t,t \in [T]$. While, for the convenience of presentation, we only consider the one-dimensional case; see Remark \ref{remark:multi-dim-extension} for a more detailed discussion.

Recall that our algorithm incorporates a nonstationarity detection
layer into SAA. Therefore, the algorithm needs to solve a general optimization problem in each period, and we introduce the following optimization oracle. The general approach is presented in Algorithm \ref{alg:saageneral}.

\textbf{SAA and Optimization Oracle.}  
During each period $t$ within epoch $\tau$, the decision $x_t$ is determined by solving the deterministic optimization problem of the SAA method (see Line \ref{Line:SAA}). Specifically, given the data from the starting period $l_{\tau}$ to the current period $t$ within an epoch, the following optimization problem is solved with an optimization oracle:  
\[
\min_{x \in \mathcal{X}} \frac{1}{t-l_{\tau}+1} \sum_{k=l_{\tau-1}}^{t} F(x, D_k).
\]  
For the implementation of our nonstationary SAA algorithm, we assume the availability of an optimization oracle $\mathcal{O}$, capable of solving this optimization problem to a desired accuracy $\epsilon$. Additional details about the optimization oracle are provided in the subsequent sections and applications.
\begin{algorithm}[h!]
\caption{Nonstationary SAA with Adaptive Restarts}
\label{alg:saageneral}
\begin{algorithmic}[1]
\State \textbf{Input:} Confidence level $\delta\in(0,1)$, optimization oracle $\O$, and accuracy $\epsilon_t, t \in [T]$.
\State \textbf{Initialization:} Set the epoch number $\tau=1$, and initial time step $t=1$, $l_1=1$.
\For{$\tau =1,2,\cdots$}
\For{$t=l_\tau,l_\tau+1,\cdots$}\label{Line:for3}
\State \textbf{If} $t = l_{\tau}$ \textbf{then} set $x_t$ randomly from $\mathcal{X}$.  
\State \label{Line:SAA} \textbf{If} $t > l_{\tau}$ \textbf{then} call the optimization oracle $\mathcal{O}$ to solve the following problem with accuracy $\epsilon_t$ and obtain the decision $x_t$:
\[
\min_{x\in\mathcal X}\frac{1}{t-l_\tau} \sum_{k=l_\tau}^{t-1} F(x,D_k).
\]

\State  \label{Line:detect3} \textbf{If} there exists $s\in [l_\tau,t]$ such that
\[
\max_{y\in \R} \vert\hat{G}_{l_\tau,t-1}(y) - \hat G_{s,t}(y) \vert > 2\sqrt{\frac{\ln (2T^2/\delta)}{t-l_\tau}}+ 2\sqrt{\frac{\ln (2T^2/\delta)}{t-s+1}}
\]
\State \textbf{then} set $\tau=\tau+1$ and $l_\tau = t+1$, and start a new epoch at Line \ref{Line:for3}.
\EndFor
\State \textbf{end for}
\EndFor
\State \textbf{end for}
\end{algorithmic}
\end{algorithm}
\subsection{Regret Analysis: General Case}
\label{subsec:general-ub}
In this section, we establish upper bound results for the general (non-convex) inner-layer function $F(\cdot,\cdot)$. 
For convenience of discussion, we first introduce some notations and assumptions as follows.

For $i\leq j$, we define
\begin{align*}
    f_{i,j}(x) &:=\int_{\R} F(x,y) \intd G_{i,j}(y) =\frac{1}{j-i+1} \sum_{k=i}^{j} f_k(x),\\
    \hat{f}_{i,j}(x) &:= \int_{\R} F(x,y) \intd \hat G_{i,j}(y) =\frac{1}{j-i+1} \sum_{k=i}^{j} F(x,D_k).
\end{align*}
 For $i = j$, we simply denote $f_{i,i}(x)$ by $f_{i}(x)$, and $\hat{f}_{i,i}(x)$ by $\hat{f}_{i}(x)$, using $i$ to denote subscribe $i,j$.

To formalize the optimization step (Line \ref{Line:SAA}) in Algorithm \ref{alg:saageneral}, we assume there is an optimization oracle that can find an $\epsilon$-optimal solution. Generally, the design of oracles depends on the structure of optimization problems. If $F(\cdot,d)$ is convex, online gradient descent, as an optimization oracle, is sufficient to achieve the optimal regret; see Remark \ref{rem:opt} after Theorem \ref{thm:general-ub}.

\begin{assumption}\label{ass:opt-oracle}
The optimization oracle $\O$ satisfies that for any $\epsilon >0$, and $1 \leq i \leq j \leq T$, applying the oracle $\O$ to the problem
\[
\min_{x \in \cX} \hat{f}_{i,j}(x) ,
\]
it can output a solution $x_{\epsilon}$ satisfying $\hat{f}_{i,j}(x_{\epsilon}) - \inf_{x\in \X}\hat{f}_{i,j}(x) \leq \epsilon$.
\end{assumption}

The following assumption requires that \( F(x, \cdot) \) possesses a well-behaved integral continuity property, ensuring that a good approximation of CDF \( G_t(y) \) leads to an accurate approximation of \( f_t(x) \) for all \( x \). This, in turn, guarantees a reliable estimation of \( f_t(x) \) at each step given the stationary data from the beginning of each epoch. 

\begin{assumption}
\label{ass:DKW-to-error}
Let ${G}(y)$ and $\hat{G}(y)$ be two CDFs satisfying $\sup_{y \in \R}|{G}(y)-\hat G(y)| \leq \epsilon$.
There exists a universal constant $L_F$ (independent of distributions) such that for all $x \in \X$,
\[
\left|\int_{\R} F(x,y) \intd G(y)- \int_{\R} F(x,y) \intd \hat G(y)\right| \leq L_F \epsilon. 
\]
\end{assumption}

The above assumption is quite mild and satisfied by a wide range of problems (see Section \ref{subsec:app} for instances). We also give a helpful lemma to verify Assumption~\ref{ass:DKW-to-error}; see Lemma \ref{lem:DKW-to-error}. 

\begin{remark}[Multi-dimensional Random Variable Extension]\label{remark:multi-dim-extension}
Note that our algorithm framework monitors the stationarity of the one-dimensional distribution $D_t, t \in [T]$.   Our algorithm framework can be further extended to monitor the multivariate CDF of $D_t, t \in [T]$, and the corresponding assumption can be modified to a multivariate integral version. It is possible to generalize the algorithm and analysis to a multi-dimensional case. Due to the limited space, we omit more complex multivariate constructions and theories for brevity.
\end{remark}

Conditioned on the above good event, we can prove the following theorem on the performance of Algorithm \ref{alg:saageneral} in terms of Switch number $S$ and total variation $V$.

\begin{theorem}\label{thm:general-ub} Suppose that Assumptions~\ref{ass:bounded}, \ref{ass:opt-oracle}, and \ref{ass:DKW-to-error} hold for the general stochastic nonstationary problem Eq.~\eqref{eq:general-opt-problem}.
With probability at least $1-\delta$, the following conclusions hold.
\begin{enumerate}
  \item \textit{Switch number $S$:} the dynamic regret of Algorithm~\ref{alg:saageneral} can be upper bounded by 
    \[ 
    \mcR^{\pi}(T) \leq 40 L_F \sqrt{(S+1)T\ln (2 T^2/\delta)}  + \sum_{t=1}^{T}\epsilon_t +2C_f (S+1).
    \]
 
    \item \textit{Total variation $V$:} the dynamic regret of Algorithm~\ref{alg:saageneral} can be upper bounded by 
\[ 
\mcR^{\pi}(T) \leq
    4C_f(V^{2/3}T^{1/3}+1) + \sum_{t=1}^{T} \epsilon_t + 96L_F (V^{1/3}+1)T^{2/3}\ln^{1/2}(2T^2/\delta).
\]
\end{enumerate}
If $\epsilon_t = \O(1/\sqrt{t})$, we have $\E[\mcR^{\pi}(T)] = \O(\sqrt{ST\log T})$ and $\E[\mcR^{\pi}(T)] = \O(V^{1/3}T^{2/3}\log^{1/2}T)$.
\end{theorem}

The main idea behind the proof is given after the conclusions of the newsvendor problem (Theorem \ref{thm:uncensored-nvp-ub-general}). Please refer to Section \ref{sec:Proof of Theorem thm:uncensored-nvp-ub-general} for detailed proof of the above theorem.

\begin{remark}[Optimization Oracle]
\label{rem:opt}
It is remarkable that if $F(x,d)$ is convex, an optimization oracle is easy to find. Specifically, in each epoch, if we maintain a sequence of $\{ y_t\}_{t\geq l_{\tau}}$ and update it by Online Gradient Descent (OGD), by the regret analysis of OGD \citep{zinkevich2003online}, we know that for $t$ in epoch $\tau$, 
\[
\frac{1}{t-l_{\tau}}\sum_{k=l_{\tau}}^{t-1} F(y_k,D_k) - \min_{x \in \mathcal X} \frac{1}{t-l_{\tau}}\sum_{k=l_{\tau}}^{t-1} F(x,D_k) \leq \mathcal{O} \left(\frac{1}{\sqrt{t-l_{\tau}}} \right).
\]
If we take $x_t$ to be the average of $\{ y_t : t \in [l_{\tau},t-1] \}$, by the convexity of $F(\cdot,D)$, we know $\epsilon_t \leq \mathcal{O}(1/\sqrt{t-l_{\tau}})$. Then similar to the proof of Theorem~\ref{thm:general-ub}, we can show
\[\sum_{t=1}^{T}
 \epsilon_t  \leq \O\left(\max \{\sqrt{ST\log T},V^{1/3}T^{2/3}\log^{1/2}T\} \right).
\]
Thus, online gradient descent is a desired oracle when $F(\cdot,D)$ is convex.
\end{remark}

\begin{remark}[Lower Bound]
The above upper bound results hold for general problems satisfying Assumptions~\ref{ass:bounded}, \ref{ass:opt-oracle}, and \ref{ass:DKW-to-error}, and the newsvendor problem can be covered as a special case of the general problem. Therefore, the lower bound results for the newsvendor problem (Theorem \ref{thm:lb_cxv}) also provide a lower bound on the general problem, which matches the upper bounds of Theorem \ref{thm:general-ub}. 
\end{remark}

\subsection{Regret Analysis: Polyak-Lojasiewicz Condition}\label{subsec:PL-ub} 
In this section, we establish an improved regret bound for Algorithm \ref{alg:saageneral} under the Polyak-Lojasiewicz (PL) condition. The PL condition not only generalizes strong convexity (satisfied by the newsvendor problem) but also covers a wider range of problems; see Section \ref{subsec:app} for applications.

\begin{assumption}[Differentiability and PL Condition]\label{ass:pl} For any realization $D$, we can find a subgradient function $\nabla_x F(x,D)$ for $F(x,D)$. For any $t \in [T]$,
$f_t(x)$ is differentiable and $\nabla f_t(x) = \E[\nabla_x F(x,D_t)]$. There exists a $\mu >0$ such that for any $t \in [T]$ and $x \in \mathcal{X}$,
\begin{equation}
    \label{eq:pl}
f_t (x) - f_t(x_t^*) \leq \mu \Vert \nabla f_t(x) \Vert^2_2.
\end{equation}
\end{assumption}

\textbf{Motivation of PL condition.} We know that under the global minimal separation condition (Definition \ref{def:global-separation}), the cost function of the newsvendor problem is strongly convex, which implies that
\begin{equation}
        f_t\left({x}\right)-f_t(x^*_t) 
    \leq \frac{1}{2\gamma}\Vert \nabla f_t(x)\Vert^2_2.
\end{equation}
Therefore, we know the PL condition (Eq.~\eqref{eq:pl}) is milder than strong convexity. 

Moreover, many optimization problems in OM with supply or demand uncertainty \citep{feng2018supply,chen2023network} satisfy an assumption of transformed convexity (Definition \ref{def:hidden-cvx}). In Section \ref{subsec:app2}, we give an example with supply uncertainty, and we show that the PL condition is milder than the transformed strongly convexity. Moreover, a growing literature focuses on designing algorithms for problems with transformed convexity \citep{miao2021network,chen2022efficient,chen2023learning}. A recent paper \cite{chen2024landscape} showed that many policy optimization problems satisfy the PL condition. These studies demonstrate the wide application of the PL condition.

Similar to the notation for objective functions, for $i\leq j$, we further define
\begin{align*}
    g_{i,j}(x) &:=\int_{\R} \nabla F(x,y) \intd G_{i,j}(y)= \frac{1}{j-i+1} \sum_{k=i}^{j} \nabla f_k(x),\\
    \hat{g}_{i,j}(x) &:=\int_{\R} \nabla F(x,y) \intd \hat G_{i,j}(y)= \frac{1}{j-i+1} \sum_{k=i}^{j} \nabla F(x,D_k).
\end{align*}
If $i = j$, we denote $g_{i,i}(x)$ by $g_{i}(x)$ and $\hat{g}_{i,i}(x)$ by $\hat{g}_{i}(x)$ similarly.

We introduce the following assumptions for the regret analysis under the PL condition.

\begin{assumption}\label{ass:opt-oracle2}
The optimization oracle $\O$ satisfies that for any $\epsilon >0$ and $1 \leq i \leq j \leq T$, applying the oracle $\O$ to the problem,
\[
\min_{x \in \cX} \hat{f}_{i,j}(x) ,
\]
it can output a solution $x_{\epsilon}$ satisfying that $ \vert \hat{g}_{i,j}(x_{\epsilon})  \vert \leq \epsilon$.
\end{assumption}

The above assumption requires that we can find an $\epsilon$-gradient approximate solution for the optimization step (Line \ref{Line:SAA}) in Algorithm \ref{alg:saa-uncensored-nvp}.  Intuitively, at the optimal solution $x$, we cannot find a feasible direction such that $\hat{f}_{i,j}(x)$ is decreasing along this direction. Therefore, if the optimal solution $x^*_{i,j}$ can be achieved in the interior of $\mathcal{X}$, we have $\vert \hat{g}_{i,j}(x^*_{i,j})  \vert=0$.

The following assumption is the gradient version of Assumption \ref{ass:DKW-to-error}.

\begin{assumption}
\label{ass:DKW-to-error2}
Let ${G}(y)$ and $\hat{G}(y)$ be two CDFs satisfying $\sup_{y \in \R}|{G}(y)-\hat G(y)| \leq \epsilon$. There exists a universal constant $L_g$ (independent of distributions) such that for all $x \in \X$,
\[
\left\Vert\int_{\R} \nabla_x F(x,y) \intd G(y)- \int_{\R} \nabla_x F(x,y) \intd \hat G_2(y)\right\Vert_2 \leq L_g \epsilon. 
\]
\end{assumption}

Conditioned on the good event (Lemma \ref{lem:good-event}), we can prove the following theorem on the performance of Algorithm \ref{alg:saa-uncensored-nvp} for functions satisfying PL condition.

\begin{theorem}
   \label{thm:pl-ub}
Suppose that Assumptions~\ref{ass:bounded}, \ref{ass:pl}, \ref{ass:opt-oracle2}, and \ref{ass:DKW-to-error2} hold for the general stochastic nonstationary problem Eq.~\eqref{eq:general-opt-problem}.
With probability at least $1-\delta$, the following conclusions hold.
\begin{enumerate}
    \item \textit{Switch number $S$:} the dynamic regret of Algorithm~\ref{alg:saageneral} can be upper bounded by 
    \[
    \mathcal{R}(T)\leq 104\mu L_g^2(S+1)\ln (2T^2/\delta)(1+\ln (T/S))+ \sum_{t=1}^T 2\mu \epsilon^2_t + C_f(S+1).
    \]
 
    \item \textit{Total variation $V$:} the dynamic regret of Algorithm~\ref{alg:saageneral} can be upper bounded by 
\[ \mathcal{R}(T) \leq 4C_f(V^{2/3}T^{1/3}+1) + 2\mu\sum_{t =1}^{T} \epsilon^2_t + 384\mu L_g^2 (V^{2/3}T^{1/3}+1)\ln T\ln^{1/2}(2T^2/\delta).
    \]
\end{enumerate}
Moreover, if $\epsilon_t = \O(1/\sqrt{t})$, we have $\mcR(T) = \O(S\log(T/S)\log T )$ and $\mcR(T) = \O(V^{2/3}T^{1/3}\log^{3/2} T )$.
\end{theorem}

Please refer to Section \ref{subsec:Proof of Theorem thm:uncensored-nvp-ub-pl} for detailed proof of the above theorem. Similar to the general case, the lower bound results for the newsvendor problem (Theorem \ref{thm:lb-str}) also provide a lower bound for the general optimization problem, which matches the upper bounds of Theorem \ref{thm:pl-ub}.

\subsection{Applications of Framework}\label{subsec:app}
In this section, we present several application examples of our algorithm~\ref{alg:saageneral}. In Section \ref{subsec:app1}, we first consider newsvendor problems and demonstrate that our general framework not only derives the upper bound results in Section \ref{sec:uncensored} but also applies to general inventory cost functions.
Second, in Section \ref{subsec:app2}, we extend the newsvendor problem and consider supply uncertainty. We show that this problem is non-convex but satisfies the PL condition. Finally, in Section \ref{subsec:app3}, we introduce a first-price auction problem, showing further applications of our framework beyond the scope of inventory management.  Due to limited space, we provide only a basic introduction to the problems. Please refer to Section \ref{appendix:verification} for detailed verifications of assumptions.
\subsubsection{Newsvendor Problems}\label{subsec:app1}
Consider the nonstationary newsvendor problem
\[
\min_{x \in \mathcal{X}} f_t(x):= \E_{D \sim G_t}[F(x,D)],
\]
where $F(x,D) = \theta(x-D)$, and $\theta(\cdot)$ is a general convex function.

We consider the following two classical specific inventory cost forms:
\begin{itemize}
    \item \textit{Linear cost.} Unsatisfied demand causes a penalty cost of $b(D-x)^+$, and unconsumed inventory causes an overage cost of $h(x-D)^+$. Thus, $\theta(x-D) = b(D-x)^+ + h(x-D)^+$. 
    \item \textit{Quadratic cost.} $\theta(x-D) = a(x-D)^2$, where $a >0$.  This cost function is an example of a strongly convex and smooth cost function.
\end{itemize}

\subsubsection{Procurement Problem under Supply Uncertainty}\label{subsec:app2}
Consider the procurement problem of a firm under supply uncertainty. In each period $t \in [T]$, the firm needs to decide the order quantity $x$ from the supplier. Due to supply uncertainty, the firm will receive $S=s (x, Z) $, where $ s (x,z)$ is the supply function and $Z$ is a random variable. The per-unit purchase cost is $c$. The stochastic cost minimization problem can be formulated as
\[
\min_{x \in \mathcal{X}} f_t(x):= \E_{(D,Z) \sim G_t}[F(x,D,Z)],
\]
where $F(x,D,Z) = a(s(x,Z)-D)^2$. We assume $D$ and $Z$ are observable.

For this problem, two random variables $D$ and $Z$ are involved. Our framework can be generalized to multi-dimensional random variables and applied to solve this problem; see Remark \ref{remark:multi-dim-extension}. All the assumptions can be verified in the same manner as in the previous newsvendor problem; we omit the details. The only exception is that, in the newsvendor problem, strong convexity (and PL condition) could be directly established. However, in the current problem, due to supply uncertainty, the problem may be non-convex \citep{chen2022efficient}. To address this, we prove the problem satisfies the PL condition under random supply and random yield; see Section \ref{appendix:app2} for details.

We also discuss the potential application of our framework by pointing out the relation of the PL condition and a transformed convexity\citep{feng2018supply}; see Section \ref{appendix:app2}.

\subsubsection{Bidding in First-Price Auction}\label{subsec:app3}
Consider a bidder who wants to win a product in a repeated first-price auction. The bidder knows the value of the product $v \in [0,1]$ to him/her. In each period $t \in [T]$, he/she needs to decide the bid price $x$ to win the product, and we assume the highest bid of other bidders is a random variable\footnote{In each period, other bidders are randomly sampled from a very large population, therefore we can assume the highest bid to be a random variable. It is also a very common assumption in online first-price auctions.} $D_t \in [0,1]$, which is nonstationary over time. If the bidder wins the auction, i.e., $x \geq D_t$, he/she will obtain a utility of $(v-x)$; otherwise, the utility is $0$. Therefore, the utility obtained by the bidder is $\E[(v-x)\cdot \I[x \geq D]]$ and we formulate the maximization of utility \textit{equivalently} as
\[
\min_{x \in \mathcal{X}} f_t(x):= \E_{D \sim G_t}[F(x,D)],
\]
where $F(x,D) = (x-v)\cdot \I[x \geq D]$ and $G_t(\cdot)$ is the CDF of the highest bid given by other bidders.

Similar to the previous applications, we can verify Assumptions \ref{ass:bounded}, \ref{ass:opt-oracle}, and \ref{ass:DKW-to-error} for this problem. Therefore, we can apply Algorithm \ref{alg:saageneral} and Theorem \ref{thm:general-ub} to this problem; see Section \ref{appendix:app3}.

One may wonder whether the improved regret bounds of Theorem \ref{thm:pl-ub} are attainable for this problem. Unfortunately, because the objective function $F(x,D) = (x-v)\cdot \I[x \geq D]$ is discontinuous, finding an unbiased gradient estimator is challenging, which prevents a straightforward application of our analysis. Indeed, even for the simpler stationary first-price auction, it remains an open question whether a sample complexity of $\O(1/\epsilon)$ is achievable \citep{guo2019settling}. The nonstationary version is substantially more difficult, and we leave it as an important direction for future research.

\section{Case Studies on Real-world Data}\label{sec:Numerical Experiments}
In this section, we compare the NSAA algorithm with several other algorithms on the nonstationary newsvendor problem using two real-world datasets to demonstrate the effectiveness of our NSAA algorithm. We consider six different algorithms: NSAA, SAA, MSAA, RSAA, W2SE, and M2SE. Please refer to Section \ref{appendix:numerical-experiments} for a more detailed description of these algorithms.

\cite{keskin2023nonstationary} also conducted experiments for the newsvendor problem and tested many of the above algorithms. For the fairness of comparison, we follow the hyper-parameter tuning of \cite{keskin2023nonstationary}, which sets $\kappa=1$ for all algorithms throughout all experiments. Note that we use the same data source as \cite{keskin2023nonstationary} to collect our dataset, but we gather more recent data, resulting in a larger dataset. For specific details, please refer to the description below.

All experiments were conducted using MATLAB R2023a on an AMD Ryzen 3.00 GHz PC.
\subsection{Description of Real-world Datasets}
\begin{figure}[th]
\begin{center}
\includegraphics[width =0.48\textwidth]{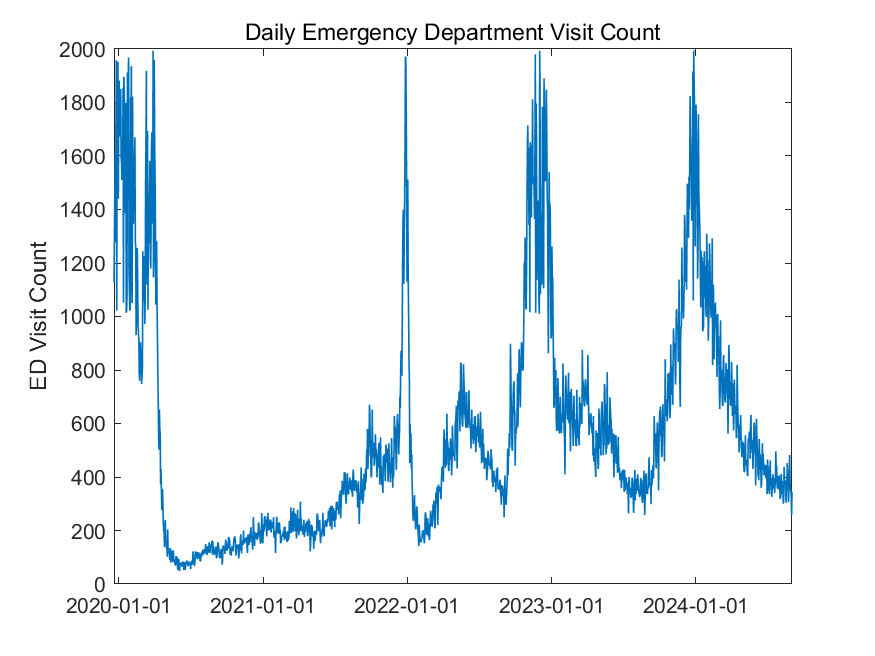}
\includegraphics[width =0.48\textwidth]{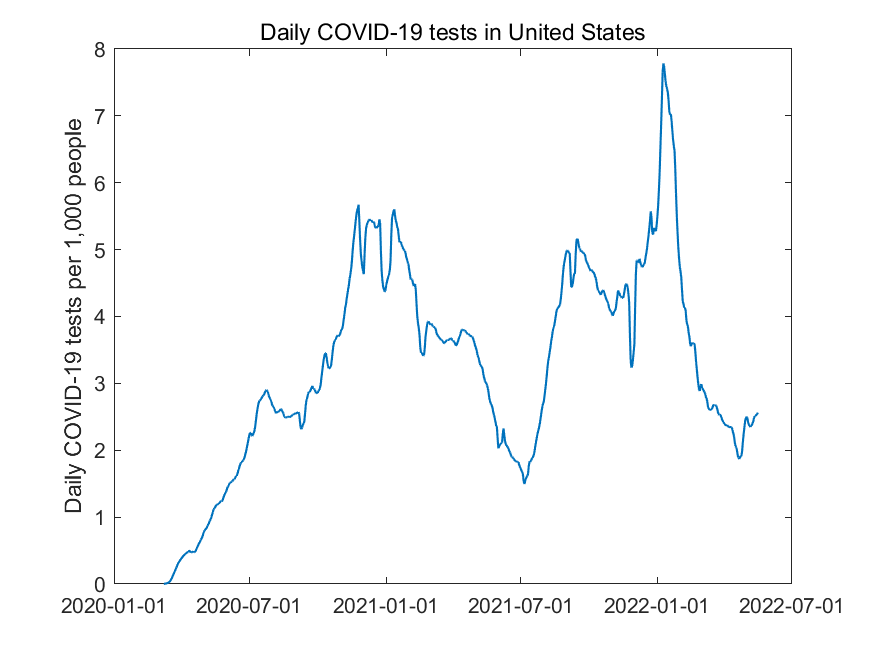}
\end{center}
\caption{{Visualization of Nursing and COVID-19 Datasets}}\label{fig:dataset}
\end{figure}

In this section, we examine two real-world datasets. The first dataset is \textit{Nursing}, and the second dataset is \textit{COVID-19}. The visualization of Nursing and COVID-19 is shown in Figure~\ref{fig:dataset}.

\textit{Nursing.} This dataset \citep{nyc_health_2020}  collects the number of emergency department patients with influenza-like illness in New York City from December 20, 2019 to September 1, 2024. It contains a total of 1,718 data points. From the figure, we can observe that the number of patients reported daily is highly nonstationary and exhibits certain seasonal characteristics. Specifically, we can see that peaks mostly occur at the end of the year, indicating that winter is a high season for flu. Hospitals must allocate enough nurses to accommodate emergency patients. On the other hand, over-scheduling will result in labor redundancy. This trade-off of cost can be well captured by the newsvendor problem formulation.

\textit{COVID-19.} This dataset includes the number of daily COVID-19 tests conducted per 1,000 people in the U.S. starting from March 8, 2020 to July 20, 2022 daily over a period of 800 days. It contains a total of 800 data points. From the graph, we can observe that as the COVID-19 pandemic spread, the number of tests increased, followed by some fluctuations, making the data highly nonstationary. COVID-19 testing requires both human and material resources, so hospitals need to allocate the necessary resources and manpower in advance to meet the daily testing demand. Over-preparing will lead to a waste of resources, while insufficient preparation results in inadequate testing. Therefore, the prediction for COVID-19 tests can be formulated as a newsvendor problem.

\subsection{Comparison under Newsvendor Framework}
As discussed earlier, for forecasting emergency department patient numbers and the demand for COVID-19 testing, the framework of the newsvendor problem can be applied. There are also many benchmark algorithms and empirical studies related to real-world datasets on setting cost
parameters.
Therefore, in the subsequent experiments, we consider that the demand $D_t$ comes from datasets mentioned earlier, and the function $F(x,d) = h(x-d)^+ + b(d-x)^+ $. Here, we fix $h=1$ and adjust the critical ratio $r = b/(h+b)$ to test the performance of different algorithms.

\textit{Evaluation metrics. } We measure performance of an algorithm $\pi$ by the \emph{cumulative cost}:
\begin{equation}
\nonumber
   \mathcal{C}^{\pi}(T) =  \sum_{t=1}^T h(x_t^{\pi}-D_t)^+ + b(D_t-x_t^{\pi})^+,
\end{equation}
where $T$ is the length of the dataset and $D_t$ is the $t$-th sample in the dataset.

For the convenience of comparison, we also measure the performance of an algorithm $\pi$  by the \emph{relative cost} compared with the performance of the NSAA algorithm defined as
\begin{equation}
\nonumber
   {R}^{\pi}(T) = \frac{\mathcal{C}^{\pi}(T)}{\mathcal{C}^{NSAA}(T)} ,
\end{equation}
where ${\mathcal{C}^{NSAA}(T)}$ is the cumulative cost of the NSAA algorithm.

\textit{Choice of critical ratio $r$.} A critical feature of the newsvendor problem is the value of critical ratio $r = b/(h+b)$. In this subsection, we set $r = 0.7$ the same as \cite{keskin2023nonstationary}. This choice of critical ratio is motivated by \cite{kolker2018optimal}, which conducted an empirical study to find the optimal staffing levels for healthcare systems under the newsvendor framework. The studies find that the critical ratios are about 0.76, 0.73, and 0.67 in different examples. Thus, for both the Nursing and COVID-19 datasets, we set the critical ratio to be 0.7.

\begin{figure}[th]
\begin{center}
\includegraphics[width =0.48\textwidth]{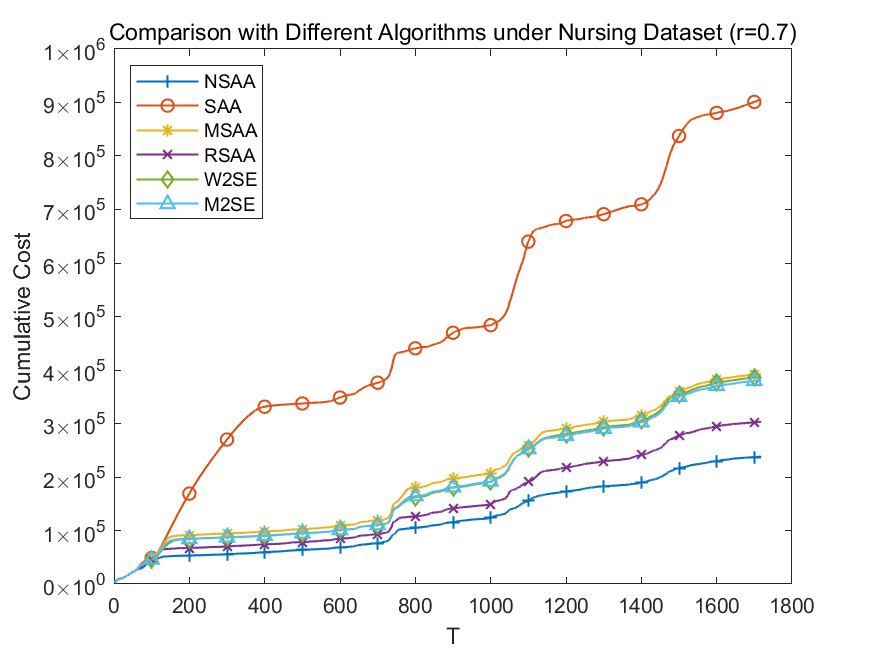}
\includegraphics[width =0.48\textwidth]{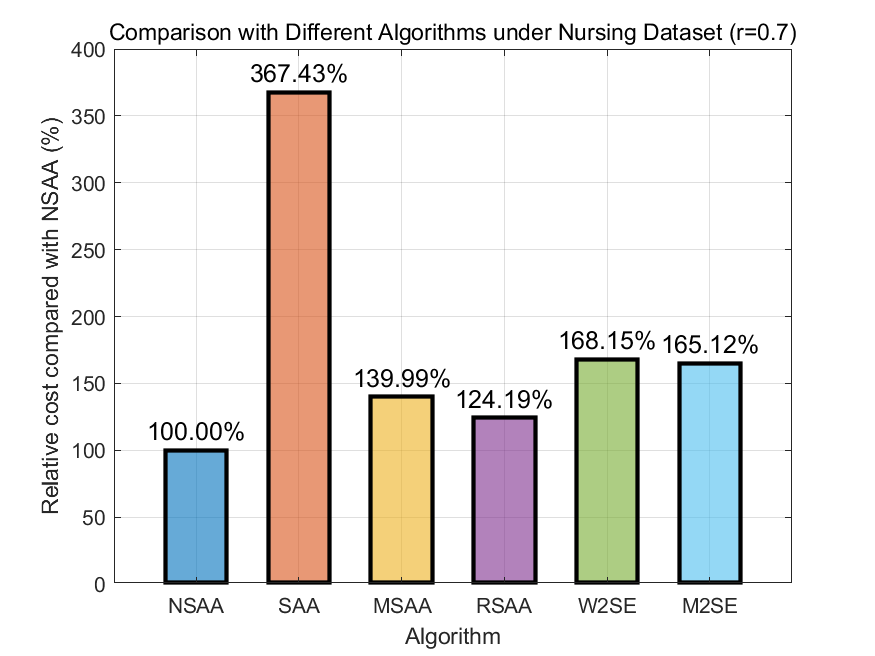}
\end{center}
\caption{{Comparison with Different Algorithms under Nursing Dataset (r=0.7)}}\label{fig:nursing-comparison}
\end{figure}

\begin{figure}[th]
\begin{center}
\includegraphics[width =0.48\textwidth]{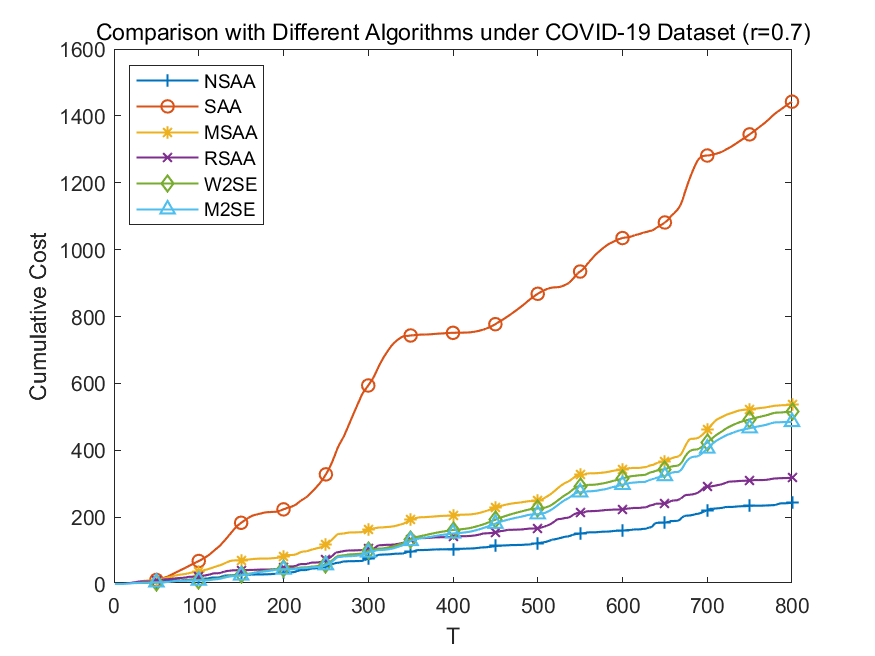}
\includegraphics[width =0.48\textwidth]{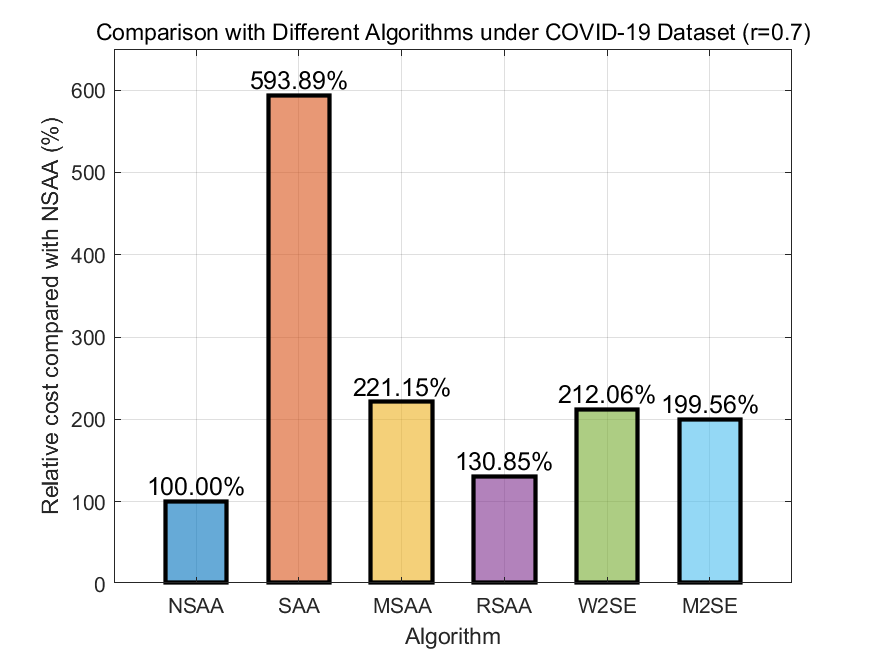}
\end{center}
\caption{{Comparison with Different Algorithms under COVID-19 Dataset (r=0.7)}}\label{fig:covid-comparison}
\end{figure}

In Figure~\ref{fig:nursing-comparison} and \ref{fig:covid-comparison}, we present the comparison results with different algorithms under the Nursing and COVID-19 datasets, respectively. We can observe from the figures that, compared to the naive SAA algorithm, all heuristic algorithms significantly reduce the cumulative cost, highlighting the importance of designing algorithms for nonstationary problems. Moreover, for both datasets, NSAA consistently performs as the best algorithm.

\subsection{Robustness to Critical Ratio}
In the previous subsection, we set the critical ratio $r= 0.7$ suggested by empirical studies. To test the robustness of the NSAA algorithm, we compare it with other algorithms under different critical ratios, $r$. Due to the limited space, please refer to Section \ref{appendix:numerical-experiments} for detailed results.

\section{Conclusion}\label{sec: Conslusion}
In this paper, we provide a comprehensive study of the nonstationary newsvendor problem across various demand distribution classes and under different observation settings and measures of nonstationarity.
We first introduce a direct, distribution-level characterization of nonstationarity within a fully nonparametric demand model, and we propose two distribution-level measures of nonstationarity that capture both abrupt and gradual changes commonly encountered in practice: the switch-number budget ($S$) and the total-variation budget ($V$). 
Then, we propose a new distributional-detection-and-restart framework for learning in nonstationary environments without prior knowledge of the degree of nonstationarity, and we instantiate it with two novel algorithms for the  uncensored and censored demand settings. Our theoretical results show that these algorithms attain minimax-optimal dynamic regret in $T$, $S$, and $V$ across all cases. This provides a comprehensive characterization that spans the full spectrum of nonstationary environments. The analytic framework is novel and broadly applicable, combining nontrivial proof techniques that should be of independent interest for other nonstationary decision-making problems. 

In addition to the newsvendor problem, the distributional-detection-and-restart framework applies broadly to a wide class of nonstationary stochastic optimization problems in OM.
Furthermore, we also extend the analysis to objective functions satisfying the PL condition, thereby covering a class of nonconvex problems that arise in OM contexts.
We believe this approach can be extended to a wider range of nonstationary decision-making and operations management problems beyond the scope of this paper.

\bibliographystyle{informs2014}
\bibliography{refs.bib}

\begin{thebibliography}{56}
\providecommand{\natexlab}[1]{#1}
\providecommand{\url}[1]{\texttt{#1}}
\providecommand{\urlprefix}{URL }

\bibitem[{An et~al.(2025)An, Li, Moseley, \protect\BIBand{} Ravi}]{an2025nonstationary}
An L, Li AA, Moseley B, Ravi R (2025) The nonstationary newsvendor with (and without) predictions. \emph{Manufacturing \& Service Operations Management} 27(3):881--896.

\bibitem[{Auer et~al.(2002)Auer, Cesa-Bianchi, Freund, \protect\BIBand{} Schapire}]{auer2002nonstochastic}
Auer P, Cesa-Bianchi N, Freund Y, Schapire RE (2002) The nonstochastic multiarmed bandit problem. \emph{SIAM journal on computing} 32(1):48--77.

\bibitem[{Auer et~al.(2019)Auer, Gajane, \protect\BIBand{} Ortner}]{auer2019adaptively}
Auer P, Gajane P, Ortner R (2019) Adaptively tracking the best bandit arm with an unknown number of distribution changes. \emph{Conference on Learning Theory}, 138--158.

\bibitem[{Besbes et~al.(2015)Besbes, Gur, \protect\BIBand{} Zeevi}]{besbes2015non}
Besbes O, Gur Y, Zeevi A (2015) Non-stationary stochastic optimization. \emph{Operations research} 63(5):1227--1244.

\bibitem[{Besbes et~al.(2019)Besbes, Gur, \protect\BIBand{} Zeevi}]{besbes2019optimal}
Besbes O, Gur Y, Zeevi A (2019) Optimal exploration--exploitation in a multi-armed bandit problem with non-stationary rewards. \emph{Stochastic Systems} 9(4):319--337.

\bibitem[{Besbes \protect\BIBand{} Mouchtaki(2023)}]{besbes2023big}
Besbes O, Mouchtaki O (2023) How big should your data really be? data-driven newsvendor: Learning one sample at a time. \emph{Management Science} 69(10):5848--5865.

\bibitem[{Besbes \protect\BIBand{} Muharremoglu(2013)}]{besbes2013implications}
Besbes O, Muharremoglu A (2013) On implications of demand censoring in the newsvendor problem. \emph{Management Science} 59(6):1407--1424.

\bibitem[{Cao et~al.(2020)Cao, Zhang, \protect\BIBand{} Poor}]{cao2020online}
Cao X, Zhang J, Poor HV (2020) Online stochastic optimization with time-varying distributions. \emph{IEEE Transactions on Automatic Control} 66(4):1840--1847.

\bibitem[{Chen et~al.(2024{\natexlab{a}})Chen, Jiang, Zhang, \protect\BIBand{} Zhou}]{chen2024learning}
Chen B, Jiang J, Zhang J, Zhou Z (2024{\natexlab{a}}) Learning to order for inventory systems with lost sales and uncertain supplies. \emph{Management Science} 70(12):8631--8646.

\bibitem[{Chen et~al.(2022)Chen, He, Hu, \protect\BIBand{} Ye}]{chen2022efficient}
Chen X, He N, Hu Y, Ye Z (2022) Efficient algorithms for minimizing compositions of convex functions and random functions and its applications in network revenue management. \emph{arXiv preprint arXiv:2205.01774} .

\bibitem[{Chen et~al.(2024{\natexlab{b}})Chen, Hu, \protect\BIBand{} Zhao}]{chen2024landscape}
Chen X, Hu Y, Zhao M (2024{\natexlab{b}}) Landscape of policy optimization for finite horizon mdps with general state and action. \emph{arXiv preprint arXiv:2409.17138} .

\bibitem[{Chen et~al.(2023)Chen, Lyu, Yuan, \protect\BIBand{} Zhou}]{chen2023learning}
Chen X, Lyu J, Yuan S, Zhou Y (2023) Learning in lost-sales inventory systems with stochastic lead times and random supplies. \emph{Available at SSRN 4671416} .

\bibitem[{Chen et~al.(2019)Chen, Wang, \protect\BIBand{} Wang}]{chen2019nonstationary}
Chen X, Wang Y, Wang YX (2019) Nonstationary stochastic optimization under l p, q-variation measures. \emph{Operations Research} 67(6):1752--1765.

\bibitem[{Chen \protect\BIBand{} Shi(2023)}]{chen2023network}
Chen Y, Shi C (2023) Network revenue management with online inverse batch gradient descent method. \emph{Production and Operations Management} 32(7):2123--2137.

\bibitem[{Chen \protect\BIBand{} Ma(2024)}]{chen2024survey}
Chen Z, Ma W (2024) Survey of data-driven newsvendor: Unified analysis and spectrum of achievable regrets. \emph{arXiv preprint arXiv:2409.03505} .

\bibitem[{Cheung \protect\BIBand{} Simchi-Levi(2019)}]{cheung2019sampling}
Cheung WC, Simchi-Levi D (2019) Sampling-based approximation schemes for capacitated stochastic inventory control models. \emph{Mathematics of Operations Research} 44(2):668--692.

\bibitem[{Cheung et~al.(2022)Cheung, Simchi-Levi, \protect\BIBand{} Zhu}]{cheung2022hedging}
Cheung WC, Simchi-Levi D, Zhu R (2022) Hedging the drift: Learning to optimize under nonstationarity. \emph{Management Science} 68(3):1696--1713.

\bibitem[{Cheung et~al.(2023)Cheung, Simchi-Levi, \protect\BIBand{} Zhu}]{cheung2023nonstationary}
Cheung WC, Simchi-Levi D, Zhu R (2023) Nonstationary reinforcement learning: The blessing of (more) optimism. \emph{Management Science} 69(10):5722--5739.

\bibitem[{Chu et~al.(2023)Chu, Feng, Shanthikumar, Shen, \protect\BIBand{} Wu}]{chu2023solving}
Chu LY, Feng Q, Shanthikumar JG, Shen ZJM, Wu J (2023) Solving the price-setting newsvendor problem with parametric operational data analytics ({{ODA}}). \emph{Available at SSRN 4400568} .

\bibitem[{Even-Dar et~al.(2006)Even-Dar, Mannor, Mansour, \protect\BIBand{} Mahadevan}]{even2006action}
Even-Dar E, Mannor S, Mansour Y, Mahadevan S (2006) Action elimination and stopping conditions for the multi-armed bandit and reinforcement learning problems. \emph{Journal of machine learning research} 7(6).

\bibitem[{Fatkhullin et~al.(2023)Fatkhullin, He, \protect\BIBand{} Hu}]{fatkhullin2023stochastic}
Fatkhullin I, He N, Hu Y (2023) Stochastic optimization under hidden convexity. \emph{arXiv preprint arXiv:2401.00108} .

\bibitem[{Feng \protect\BIBand{} Shanthikumar(2018)}]{feng2018supply}
Feng Q, Shanthikumar JG (2018) Supply and demand functions in inventory models. \emph{Operations Research} 66(1):77--91.

\bibitem[{Feng \protect\BIBand{} Shanthikumar(2023)}]{feng2023framework}
Feng Q, Shanthikumar JG (2023) The framework of parametric and nonparametric operational data analytics. \emph{Production and Operations Management} 32(9):2685--2703.

\bibitem[{Feng et~al.(2025)Feng, Shanthikumar, \protect\BIBand{} Wu}]{feng2025contextual}
Feng Q, Shanthikumar JG, Wu J (2025) Contextual data-integrated newsvendor solution with operational data analytics (oda). \emph{Management Science} .

\bibitem[{Flaxman et~al.(2004)Flaxman, Kalai, \protect\BIBand{} McMahan}]{flaxman2004online}
Flaxman AD, Kalai AT, McMahan HB (2004) Online convex optimization in the bandit setting: gradient descent without a gradient. \emph{arXiv preprint cs/0408007} .

\bibitem[{Gill \protect\BIBand{} Levit(1995)}]{Richard1995ApplyofvanTrees}
Gill RD, Levit BY (1995) Applications of the van trees inequality: a bayesian cram{\'e}r-rao bound. \emph{Bernoulli} 59--79.

\bibitem[{Guo et~al.(2019)Guo, Huang, \protect\BIBand{} Zhang}]{guo2019settling}
Guo C, Huang Z, Zhang X (2019) Settling the sample complexity of single-parameter revenue maximization. \emph{Proceedings of the 51st Annual ACM SIGACT Symposium on Theory of Computing}, 662--673.

\bibitem[{Huang \protect\BIBand{} Wang(2023)}]{huang2023stability}
Huang C, Wang K (2023) A stability principle for learning under non-stationarity. \emph{arXiv preprint arXiv:2310.18304} .

\bibitem[{Huh et~al.(2009)Huh, Janakiraman, Muckstadt, \protect\BIBand{} Rusmevichientong}]{huh2009adaptive}
Huh WT, Janakiraman G, Muckstadt JA, Rusmevichientong P (2009) An adaptive algorithm for finding the optimal base-stock policy in lost sales inventory systems with censored demand. \emph{Mathematics of Operations Research} 34(2):397--416.

\bibitem[{Huh \protect\BIBand{} Rusmevichientong(2009)}]{huh2009nonparametric}
Huh WT, Rusmevichientong P (2009) A nonparametric asymptotic analysis of inventory planning with censored demand. \emph{Mathematics of Operations Research} 34(1):103--123.

\bibitem[{Keskin et~al.(2023)Keskin, Min, \protect\BIBand{} Song}]{keskin2023nonstationary}
Keskin NB, Min X, Song JSJ (2023) The nonstationary newsvendor: Data-driven nonparametric learning. \emph{Available at SSRN 3866171} .

\bibitem[{Kim et~al.(2022)Kim, Madden, \protect\BIBand{} Dall’Anese}]{kim2022online}
Kim S, Madden L, Dall’Anese E (2022) Online stochastic gradient methods under sub-weibull noise and the polyak-{\l}ojasiewicz condition. \emph{2022 IEEE 61st Conference on Decision and Control (CDC)}, 3499--3506 (IEEE).

\bibitem[{Kolker(2018)}]{kolker2018optimal}
Kolker A (2018) The optimal workforce staffing solutions with random patient demand in healthcare settings. \emph{Encyclopedia of Information Science and Technology, Fourth Edition}, 3711--3724 (IGI Global Scientific Publishing).

\bibitem[{Levi et~al.(2015)Levi, Perakis, \protect\BIBand{} Uichanco}]{levi2015data}
Levi R, Perakis G, Uichanco J (2015) The data-driven newsvendor problem: new bounds and insights. \emph{Operations Research} 63(6):1294--1306.

\bibitem[{Levi et~al.(2007)Levi, Roundy, \protect\BIBand{} Shmoys}]{levi2007nearoptimal}
Levi R, Roundy RO, Shmoys DB (2007) Provably near-optimal sampling-based policies for stochastic inventory control models. \emph{Mathematics of Operations Research} 32(4):821--839.

\bibitem[{Lin et~al.(2022)Lin, Huh, Krishnan, \protect\BIBand{} Uichanco}]{lin22datadriven}
Lin M, Huh WT, Krishnan H, Uichanco J (2022) Technical note—data-driven newsvendor problem: Performance of the sample average approximation. \emph{Operations Research} 70(4):1996--2012.

\bibitem[{Lyu et~al.(2023)Lyu, Xie, Yuan, \protect\BIBand{} Zhou}]{lyu2023minibatch}
Lyu J, Xie J, Yuan S, Zhou Y (2023) A minibatch-sgd-based learning meta-policy for inventory systems with myopic optimal policy. \emph{Available at SSRN 4390778} .

\bibitem[{Lyu et~al.(2024)Lyu, Yuan, Zhou, \protect\BIBand{} Zhou}]{lyu2024closing}
Lyu J, Yuan S, Zhou B, Zhou Y (2024) Closing the gaps: Optimality of sample average approximation for data-driven newsvendor problems. \emph{arXiv preprint arXiv:2407.04900} .

\bibitem[{Mao et~al.(2024)Mao, Zhang, Zhu, Simchi-Levi, \protect\BIBand{} Ba{\c{s}}ar}]{mao2024model}
Mao W, Zhang K, Zhu R, Simchi-Levi D, Ba{\c{s}}ar T (2024) Model-free nonstationary reinforcement learning: Near-optimal regret and applications in multiagent reinforcement learning and inventory control. \emph{Management Science} .

\bibitem[{Miao \protect\BIBand{} Wang(2021)}]{miao2021network}
Miao S, Wang Y (2021) Network revenue management with nonparametric demand learning:$\backslash$sqrt $\{$T$\}$-regret and polynomial dimension dependency. \emph{Available at SSRN 3948140} .

\bibitem[{Mokhtari et~al.(2016)Mokhtari, Shahrampour, Jadbabaie, \protect\BIBand{} Ribeiro}]{mokhtari2016online}
Mokhtari A, Shahrampour S, Jadbabaie A, Ribeiro A (2016) Online optimization in dynamic environments: Improved regret rates for strongly convex problems. \emph{2016 IEEE 55th Conference on Decision and Control (CDC)}, 7195--7201 (IEEE).

\bibitem[{{NYC Health}(2024)}]{nyc_health_2020}
{NYC Health} (2024) Syndromic surveillance data. \url{https://a816-health.nyc.gov/hdi/epiquery/visualizations?PageType=ps&PopulationSource=Syndromic}, accessed on September 14, 2024.

\bibitem[{Pun \protect\BIBand{} Shames(2024)}]{pun2024online}
Pun YM, Shames I (2024) Online non-stationary stochastic quasar-convex optimization. \emph{arXiv preprint arXiv:2407.03601} .

\bibitem[{Qin et~al.(2022)Qin, Simchi-Levi, \protect\BIBand{} Wang}]{qin2022data}
Qin H, Simchi-Levi D, Wang L (2022) Data-driven approximation schemes for joint pricing and inventory control models. \emph{Management Science} 68(9):6591--6609.

\bibitem[{Qin et~al.(2023)Qin, Simchi-Levi, \protect\BIBand{} Zhu}]{qin2023sailing}
Qin H, Simchi-Levi D, Zhu R (2023) Sailing through the dark: Provably sample-efficient inventory control. \emph{Available at SSRN 4652347} .

\bibitem[{Russac et~al.(2019)Russac, Vernade, \protect\BIBand{} Capp{\'e}}]{russac2019weighted}
Russac Y, Vernade C, Capp{\'e} O (2019) Weighted linear bandits for non-stationary environments. \emph{Advances in Neural Information Processing Systems} 32.

\bibitem[{Shi et~al.(2016)Shi, Chen, \protect\BIBand{} Duenyas}]{shi2016nonparametric}
Shi C, Chen W, Duenyas I (2016) Nonparametric data-driven algorithms for multiproduct inventory systems with censored demand. \emph{Operations Research} 64(2):362--370.

\bibitem[{Wang(2023)}]{wang2023adaptivity}
Wang Y (2023) On adaptivity in nonstationary stochastic optimization with bandit feedback. \emph{Operations Research} .

\bibitem[{Wang et~al.(2023)Wang, Gao, Zhang, \protect\BIBand{} Zhou}]{wang2023hybrid}
Wang Z, Gao X, Zhang K, Zhou S (2023) A hybrid sampling based and gradient descent method with applications in inventory management. \emph{Available at SSRN} .

\bibitem[{Wei \protect\BIBand{} Luo(2021)}]{wei2021non}
Wei CY, Luo H (2021) Non-stationary reinforcement learning without prior knowledge: An optimal black-box approach. \emph{Conference on learning theory}, 4300--4354 (PMLR).

\bibitem[{Yuan et~al.(2021)Yuan, Luo, \protect\BIBand{} Shi}]{yuan2021marrying}
Yuan H, Luo Q, Shi C (2021) Marrying stochastic gradient descent with bandits: Learning algorithms for inventory systems with fixed costs. \emph{Management Science} 67(10):6089--6115.

\bibitem[{Zhang et~al.(2018)Zhang, Chao, \protect\BIBand{} Shi}]{zhang2018perishable}
Zhang H, Chao X, Shi C (2018) Perishable inventory systems: Convexity results for base-stock policies and learning algorithms under censored demand. \emph{Operations Research} 66(5):1276--1286.

\bibitem[{Zhang et~al.(2020)Zhang, Chao, \protect\BIBand{} Shi}]{zhang2020closing}
Zhang H, Chao X, Shi C (2020) Closing the gap: A learning algorithm for lost-sales inventory systems with lead times. \emph{Management Science} 66(5):1962--1980.

\bibitem[{Zhang et~al.(2021)Zhang, Gao, Wang, \protect\BIBand{} Zhou}]{zhang2021sampling}
Zhang K, Gao X, Wang Z, Zhou S (2021) Sampling-based approximation for serial multi-echelon inventory system. \emph{Management Science} Forthcoming.

\bibitem[{Zhao et~al.(2020)Zhao, Zhang, Jiang, \protect\BIBand{} Zhou}]{pmlr-v108-zhao20a}
Zhao P, Zhang L, Jiang Y, Zhou ZH (2020) A simple approach for non-stationary linear bandits. Chiappa S, Calandra R, eds., \emph{Proceedings of the Twenty Third International Conference on Artificial Intelligence and Statistics}, volume 108 of \emph{Proceedings of Machine Learning Research}, 746--755 (PMLR).

\bibitem[{Zinkevich(2003)}]{zinkevich2003online}
Zinkevich M (2003) Online convex programming and generalized infinitesimal gradient ascent. \emph{Proceedings of the 20th international conference on machine learning (icml-03)}, 928--936.

\end{thebibliography}
\ECSwitch


\begin{center}
 		\Large{\bf Supplementary Materials to {``Learning When to Restart: Nonstationary Newsvendor from Uncensored to Censored Demand''} }
 	\end{center}

 	\vspace{10pt}

\section{Technical Lemmas}

\begin{lemma}[Nonstationary Dvoretzky–Kiefer–Wolfowitz inequality]
\label{lem:DKW}
    For independent random variables $\{ D_{s},D_{s+1},\dots,D_t \}$, we have
    \[
    \P \left[\sup_{y \in \R} |\hat{G}_{s,t}(y) -\E[ G_{s,t}(y)]|>\lambda \right] \leq 2e^{-2(t-s+1)\lambda^2},
    \]
    for $\lambda >0$ and every $t\in \N^+$.
\end{lemma}

\begin{lemma}[Bretagnolle-Huber inequality]\label{le:Bretagnolle-Huber inequality} 
Let $P$ and $Q$ be probability measures on the same measurable space $(\Omega, \mathcal{F})$, and let $A \in \mathcal{F}$ be an arbitrary event. Then,
$$
P(A)+Q\left(A^c\right) \geq \frac{1}{2} \exp (-\mathrm{KL}(P, Q)),
$$
where $A^c=\Omega \backslash A$ is the complement of $A$.
\end{lemma}

\begin{lemma}[The Van Trees inequality \citep{Richard1995ApplyofvanTrees}]\label{le:vantree} 
Let $(\mathcal{X},\mathcal{F},P_\theta:\theta\in\Theta,\mu)$ be a dominated family of distributions, where $\mathcal{X}$ is the sample space, $\mathcal{F}$ is the $\sigma$-algebra, $P_\theta$ is the probability measure with parameter $\theta$,  the parameter space $\Theta$ is a closed interval on the real line, $\mu$ is the dominating measure.  Let $f(x|\theta)$ denote the density of $P_\theta$ with respect to $\mu$, and $\pi$ be some probability distribution on $\Theta$ with a density $\lambda(\theta)$ with respect to the Lebesgue measure.
Suppose the following assumptions hold:\begin{enumerate}
    \item $\lambda(\cdot)$ and $f(x|\cdot)$ are absolutely continuous ($\mu$-almost surely),
    \item for any fixed $\theta\in \Theta$, $\E_\theta\left[\frac{\partial}{\partial\theta}\left(\log f\left(X|\theta\right)\right)\right]=0,$ where $\E_\theta$ denotes the expectation over $X\sim P_\theta$.
    \item $\lambda(\cdot)$ converges to zero at the endpoints of the interval $\Theta$.
\end{enumerate}
Then,  for any estimator $\hat{\theta}(X)$ based on sample $X\sim P_\theta$, it holds that
$$
    \E\left[\left(\hat{\theta}\left(X\right)-\theta\right)^2\right]\geq \frac{1}{\E\left[I\left(\theta\right)\right]+I\left(\lambda\right)}, 
$$
where $\E$ denotes the expectation over the ensuing joint distribution of $X$ and $\theta$, $I(\theta) = \E_\theta\left[(\frac{\partial}{\partial\theta}\left(\log f\left(X|\theta\right)\right))^2\right]$ and $I(\lambda) = \E\left[(\frac{\partial}{\partial\theta}\left(\log \lambda\left(\theta\right)\right))^2\right]$ is the the Fisher information for $\theta$ and $\lambda$ respectively.
\end{lemma}

\section{Proofs Omitted in Section~\ref{sec:uncensored}}
\label{sec:Proofs Omitted in Sectionsec:uncensored}
In this section, we give the proof omitted in Section \ref{sec:uncensored}.

\subsection{Proof of  Theorem~\ref{thm:lb_cxv}}
\label{sec:Proof of  Theoremthm:lb_cxv}
\subsubsection{Proof of $\Omega(\sqrt{ST})$ }
\label{subsec:Proof of Omega(sqrtST)}
We first construct the problem instances used to establish the lower bound results.
We divide the $T$-period horizon into batches of $\Delta_T$ periods (the periods of the last batch might be less than $\Delta_T$). In each batch,  the demand randomly follows either $D_a$ or $D_b$ with equal probability, and the demand distributions of all periods are kept unchanged within a batch. Specifically, we define $D_a$ and $D_b$ as follows
\[
f_{D_a}(x)= 
\left\{
\begin{aligned}
&\frac{1}{2}+\epsilon,  && 0 \leq x<1 \\
&0, && 1 \leq x<3 \\
&\frac{1}{2}- \epsilon, && 3 \leq x<4
\end{aligned}
\right.
~~~\text{and}~~~
f_{D_b}(x)= 
\left\{
\begin{aligned}
&\frac{1}{2}-\epsilon,  && 0 \leq x<1 \\
&0, && 1 \leq x<3 \\
&\frac{1}{2}+\epsilon, && 3 \leq x<4,
\end{aligned}
\right.
\]
where we assume $\epsilon\in(0,1/4)$.

Let  $x_t^\pi(\mathcal H_{t-1})$ be the order-up-to level in period $t$ for policy $\pi$ with history information $\mathcal H_{t-1}:=\{D_1,\dots,D_{t-1}\}$. Recall that $F(x,d) = h\cdot(x-d)^++b\cdot(d-x)^+$ and we denote the optimal solution under demand $D$ as $x^*_D$.
For the first batch, we have
\begin{align}\nonumber
    &\sup_{D\in\{D_a,D_b\}}\left\{\sum_{t=1}^{\Delta_T}\E^{\pi}_D[F(x_t^\pi(\mathcal H_{t-1}),D) - F(x^*_{D},D)] \right\}\\&\geq \frac{1}{2}
    \sum_{t=1}^{\Delta_T}\E^{\pi}_{D_a}[F(x_t^\pi(\mathcal H_{t-1}),D_a) - F(x^*_{D_a},D_a)]+ \frac{1}{2}\sum_{t=1}^{\Delta_T}\E^{\pi}_{D_b}[F(x_t^\pi(\mathcal H_{t-1}),D_b) - F(x^*_{D_b},D_b)] \label{eq:lbthm1_eq_1} 
\end{align}
Recalling $h=1$ and $b=1$, we have $x^*_{D_a}\in (0,1)$ and $x^*_{D_b}\in (3,4)$. Thus, for $x\geq 2$, it holds that 
\begin{align*}
    F(x,D_a) - F(x^*_{D_a},D_a)&\geq   F(2,D_a) - F(1,D_a)
    \\&= \P[D_a\leq 1] -\P[D_a\geq 2]= 2\epsilon,  
\end{align*}
where the first inequality is because $F(x,D_a)$ is non-decreasing as $x\geq 2\geq x^*_{D_a}$, and the other inequalities are by the definition of $D_a$ and direct computation.

By the above inequality, we obtain
\begin{align} \label{eq:lbthm1_eq_2}
    \sum_{t=1}^{\Delta_T}\E^{\pi}_{D_a}[F(x_t^\pi(\mathcal H_{t-1}),D_a) - F(x^*_{D_a},D_a)]\geq 2 \epsilon \sum_{t=1}^{\Delta_T}\P^{\pi}_{D_a}[x_t^\pi(\mathcal H_{t-1}) > 2].
\end{align}
Similarly, it holds that
\begin{align} \label{eq:lbthm1_eq_3}
    \sum_{t=1}^{\Delta_T}\E^{\pi}_{D_b}[F(x_t^\pi(\mathcal H_{t-1}),D_b) - F(x^*_{D_b},D_b)]\geq 2 \epsilon \sum_{t=1}^{\Delta_T}\P^{\pi}_{D_b}[x_t^\pi(\mathcal H_{t-1}) \leq 2].
\end{align}
By Bretagnolle–Huber inequality (Lemma \ref{le:Bretagnolle-Huber inequality}),  it holds that 
\begin{align} \label{eq:lbthm1_eq_4}
    \P^{\pi}_{D_a}[x_t^\pi(\mathcal H_{t-1}) > 2] + \P^{\pi}_{D_b}[x_t^\pi(\mathcal H_{t-1}) \leq 2]&\geq \frac{1}{2}
   \exp\left(-\mathrm{KL}(\P_{D_a}^{(t-1)}, \P_{D_b}^{(t-1)})\right),
\end{align}
where $\mathrm{KL}(\P_{D_a}^{(t-1)}, \P_{D_b}^{(t-1)})$ is the Kullback-Leibler divergence divergence between the distribution of $\{D_1, \dots, D_t\}$ under $D_a$ and $D_b$.
By the chain rule of the Kullback-Leibler divergence, we have 
\begin{align}\nonumber
    \mathrm{KL}(\P_{D_a}^{(t-1)}, \P_{D_b}^{(t-1)}) &= \E_{X \sim D_a}\left[\ln\left(\frac{f_{D_a}(X)}{f_{D_b}(X)}\right)\right]
\\& =(t-1)\left[(1/2+\epsilon) \ln \left(\frac{1/2+\epsilon}{1/2-\epsilon}\right)+(1/2-\epsilon) \ln \left(\frac{1/2-\epsilon}{1/2+\epsilon}\right)\right]\nonumber
\\& \leq (t-1)\left[(1/2+\epsilon)\left(4 \epsilon+8 \epsilon^2\right)-4(1/2-\epsilon) \epsilon\right] \nonumber
\\& \leq {16 (t-1) \epsilon^2}, \label{eq:lbthm1_eq_5}
\end{align}
where the first inequality is due to $2 x \leq \ln ((1+x)/(1-x)) \leq 2 x+2 x^2$ as $x \in(0,1 / 2)$, and the second inequality is due to $\epsilon\leq 1/4$.

Combining Eqs.~(\ref{eq:lbthm1_eq_1}, \ref{eq:lbthm1_eq_2}, \ref{eq:lbthm1_eq_3}, \ref{eq:lbthm1_eq_4}, \ref{eq:lbthm1_eq_5}), we obtain the worst regret incurred in a batch
\begin{align}\label{eq:lbthm1_eq_6}
    \sup_{D\in\{D_a,D_b\}}\left\{\sum_{t=1}^{\Delta_T}\E^{\pi}_D[F(x_t^\pi(\mathcal H_{t-1}),D) - F(x^*_{D},D)] \right\}\geq \frac{1}{2}\epsilon \sum_{t=1}^{\Delta_T}\exp\left(-16 (t-1) \epsilon^2 \right)
\end{align}
Recalling the design of the hard instance, we obtain
\begin{align}\nonumber
    \sup_{D\in\{D_a,D_b\}}\left\{\sum_{t=1}^{T}\E^{\pi}_D[F(x_t^\pi(\mathcal H_{t-1}),D) - F(x^*_{D},D)] \right\}&\geq \frac{1}{2} \left\lceil\frac{T}{\Delta_T}\right\rceil \epsilon \sum_{t=1}^{\Delta_T}\exp\left(-16 (t-1) \epsilon^2 \right)\\ \nonumber
    &\geq \frac{1}{2} \left\lceil\frac{T}{\Delta_T}\right\rceil \epsilon \sum_{t=1}^{\min\{\epsilon^{-2},\Delta_T\}}\exp(-16)\\
     &\geq \frac{\exp(-16)}{2} \left\lceil\frac{T}{\Delta_T}\right\rceil \min\{\epsilon^{-1},\epsilon\Delta_T\},\label{eq:lb-inequality}
\end{align}
where the second inequality is due to $(t-1)\epsilon^2\leq 1$.

Setting $\Delta_T = \lceil T/S\rceil$ and $\epsilon = 1/\sqrt{\Delta_T}$, by the above inequality, we have  
\begin{align*}
&\sup_{\bm{D}\in\mathcal D_{S}^{(1)}}\left\{\sum_{t=1}^{T}\E^{\pi}_D[F(x_t^\pi(\mathcal H_{t-1}),D) - F(x^*_{D},D)] \right\}\\
&\geq \sup_{D\in\{D_a,D_b\}}\left\{\sum_{t=1}^{\Delta_T}\E^{\pi}_D[F(x_t^\pi(\mathcal H_{t-1}),D) - F(x^*_{D},D)] \right\}\\
&\geq \frac{\exp(-16)}{2} \left\lceil\frac{T}{{\Delta_T}}\right\rceil \sqrt{\Delta_T} \geq \frac{\exp(-16)}{4} \sqrt{TS},
\end{align*}
where the last inequality is by some computation, and the first inequality is due to the problem instance satisfying that
\[\sum_{t=2}^{T}\mathbb{I}\bigl[D_t \overset{d}{\neq} D_{t-1} \bigr]\le \lceil T/\Delta_T\rceil \leq \lceil T/ \lceil T/S \rceil \rceil \leq S,\]
where we complete the proof.
\Halmos
\endproof

\subsubsection{Proof of $\Omega(V^{1/3}T^{2/3})$}
\label{subsec:proof-of-V-general-lb}
We follow the main steps of the proof of $\Omega(\sqrt{ST})$. 
We define $D_a$ and $D_b$ as the same as the proof of $\Omega(\sqrt{ST})$ and divide the $T$-period horizon into batches of $\Delta_T$ periods. In each batch, the demand randomly follows either $D_a$ or $D_b$. Note that we will provide a different batch design to achieve the lower bound in terms of the variation budget $V$.

Since $\mathrm{TV}(D_a, D_b)=4\epsilon$,
setting $\Delta_T = \lceil (T/V)^{2/3}\rceil$  and $\epsilon = 1/(4\sqrt{\Delta_T})$,
 we have  
\[\sum_{t=2}^{T}\|D_{t}-D_{t-1}\|_{\mathrm{TV}}\leq 4\epsilon \cdot \lceil T/\Delta_T\rceil \leq V.\]

Similar to Eq.~\eqref{eq:lb-inequality}, we have
\begin{align*}
 \sup_{\bm{D}\in\mathcal D_{V}^{(1)}}\left\{\sum_{t=1}^{T}\E^{\pi}_D[F(x_t^\pi(\mathcal H_{t-1}),D) - F(x^*_{D},D)] \right\}&\geq \frac{\exp(-16)}{2} \left\lceil\frac{T}{\Delta_T}\right\rceil \min\{\epsilon^{-1},\epsilon\Delta_T\}
   \\ &\geq \frac{\exp(-16)}{8} \left\lceil\frac{T}{\Delta_T}\right\rceil \sqrt{\Delta_T}\\ &\geq \frac{\exp(-16)}{16} V^{1/3}T^{2/3},
\end{align*}
where we complete the proof.
\Halmos
\endproof

\subsection{Proof of Theorem \ref{thm:lb-str}.}
\label{subsec:Proof of Theorem thm:lb-str}

\subsubsection{Proof of $\Omega(S\log(T/S))$}

We first construct the problem instances used to establish the lower bound result. We divide the $T$-period horizon into batches of $\Delta_T$ periods (the periods of the last batch might be less than $\Delta_T$). In each batch,  the demand randomly follows a distribution family parameterized by $\theta$ with a certain density $q(\theta)$, and the demand distributions of all periods are kept unchanged within a batch. Note this construction is based on the hard instances in \cite{lyu2024closing}, and we will omit all verification steps that can be found in \cite{lyu2024closing} (with $\alpha$ and $\rho$ therein equal to $1/2$).

Specifically, we define density function $q(\theta)$ as follows
$$
        q(\theta)=40\cos^2\left({20\pi}\theta\right)\mathbb{I}\left[\theta\in\left[-\frac{1}{40},\frac{1}{40}\right]\right],
$$
where $[-1/40,1/40]$  is the parameter space and we denote it by $\Theta$. 

For $\theta\in\Theta$, we define density function $\rho(x\mid \theta)$ of demand as
$$
        \rho \left(x|\theta\right)=\left\{\begin{aligned}
		&3/2,&x\in\left[0,l_1(\theta)\right]\cup\left[r_1(\theta),1\right],\\
            &1/2,&x\in\left(l_2(\theta),r_2(\theta)\right],\\
                &1/2+\left(\cos \left(w_1\left(x-l_1(\theta)\right)\right)+1\right)/2,&x\in\left[l_1(\theta),l_2(\theta)\right],\\
                &1/2 +\left(\cos \left(w_2\left(r_1(\theta)-x\right)+1\right)\right)/2,&x\in\left[r_2(\theta),r_1(\theta)\right],
                \end{aligned}
            \right.
    $$
    where $l_1(\theta)=3/24$, $l_2(\theta)=l_1(\theta)+1/4$, $r_2(\theta)=l_2(\theta)+1/4$, $r_1(\theta)=r_2(\theta)+1/4$, and $w_1=w_2=4\pi$.

Let $\mathcal{F}$ be all distributions such that $f(x)$ is strongly convex. We have
\begin{align}\label{eq:integral_expression_lb}
        \underset{\pi\in\Pi}{\inf}\underset{F\in\mathcal{F}_{\alpha}}{\sup}\{R^{\pi}\left(F,\Delta_T\right)\}\geq \underset{\pi\in\Pi}{\inf}\int_{-\alpha/20}^{\alpha/20}\mathcal{R}^{\pi}(F_\theta,\Delta_T)q(\theta)\mathrm{d}\theta.
    \end{align}

Recalling the design of the hard instance, we obtain
\begin{align}\nonumber
\sum_{t=1}^{T}\E^{\pi}_D[F(x_t^\pi(\mathcal H_{t-1}),D) - F(x^*_{D},D)] &\geq \left\lceil\frac{T}{\Delta_T}\right\rceil \cdot  \underset{\pi\in\Pi} {\inf}\underset{F\in\mathcal{F}_{\alpha}}{\sup}\{R^{\pi}\left(F,\Delta_T\right)\}\\ \nonumber 
&\geq \left\lceil\frac{T}{\Delta_T}\right\rceil \cdot \frac{1}{200\pi^2}(\ln \Delta_T-3\ln 2),
\end{align}
where the last inequality is by Van Trees inequality and the construction of hard instance (Please refer to the proof of Theorem 2 in \cite{lyu2024closing} for detailed proof of this inequality).

Setting $\Delta_T = \lceil T/S\rceil$ and $\epsilon = 1/\sqrt{\Delta_T}$, by the above inequality, we have  
\begin{align*}
\sup_{\bm{D}\in\mathcal D_{S}^{(1)}}\left\{\sum_{t=1}^{T}\E^{\pi}_D[F(x_t^\pi(\mathcal H_{t-1}),D) - F(x^*_{D},D)] \right\} \geq \frac{1}{400\pi^2} S \log (T/8S),
\end{align*}
where we complete the proof.
\Halmos
\endproof
\subsubsection{Proof of $\Omega(V^{2/3}T^{1/3})$}

We first construct the problem instances that are utilized to establish the lower bound results. The idea is similar to the proof of Theorem \ref{thm:lb_cxv}, and we divide the $T$-period horizon into batches of $\Delta_T$ periods. In each batch, the demand randomly follows either $D_a$ or $D_b$. In each batch, the demand randomly follows either $D_a$ or $D_b$ and the demand distributions of all periods are kept unchanged within a batch. $D_a$ and $D_b$ are defined as follows,
\[
f_{D_a}(x)= 
\left\{
\begin{aligned}
&\frac{1}{2}+\epsilon,  && 0 \leq x<1 \\
&\frac{1}{2}- \epsilon, && 1 \leq x<2
\end{aligned}
\right.
~~~\text{and}~~~
f_{D_b}(x)= 
\left\{
\begin{aligned}
&\frac{1}{2}-\epsilon,  && 0 \leq x<1 \\
&\frac{1}{2}+\epsilon, && 1 \leq x<2  ,
\end{aligned}
\right.
\]
where we assume $\epsilon\in(0,1/4)$.

Let  $x_t^\pi(\mathcal H_{t-1})$ be the order-up-to level in period $t$ for policy $\pi$ with history information $\mathcal H_{t-1}:=\{D_1,\dots,D_{t-1}\}$. 
Recall that $F(x,d) = h\cdot(x-d)^++b\cdot(d-x)^+$ and we denote the optimal order-up-to level under demand $D$ as $x^*_D$.
For the first batch
\begin{align}\nonumber
    &\sup_{D\in\{D_a,D_b\}}\left\{\sum_{t=1}^{\Delta_T}\E^{\pi}_D[F(x_t^\pi(\mathcal H_{t-1}),D) - F(x^*_{D},D)] \right\}\\&\geq \frac{1}{2}
    \sum_{t=1}^{\Delta_T}\E^{\pi}_{D_a}[F(x_t^\pi(\mathcal H_{t-1}),D_a) - F(x^*_{D_a},D_a)]+ \frac{1}{2}\sum_{t=1}^{\Delta_T}\E^{\pi}_{D_b}[F(x_t^\pi(\mathcal H_{t-1}),D_b) - F(x^*_{D_b},D_b)] \label{eq:str-lbthm1_eq_1} 
\end{align}
Recalling $h=1$ and $b=1$, we have $$x^*_{D_a} = \frac{1}{1+2\epsilon}\in (0,1) \text{~~~~and~~~~}  x^*_{D_b} = 1+\frac{2\epsilon}{1+2\epsilon} \in (1,2).$$ 

Thus, for $x\geq 1$, it holds that 
\begin{align*}
    F(x,D_a) - F(x^*_{D_a},D_a)&\geq   F(1,D_a) - F(x^*_{D_a},D_a)\geq\frac{1}{4}(1-x^*_{D_a})^2 \\&\geq\frac{1}{4}\frac{4\epsilon^2}{(1+2\epsilon)^2} \geq\frac{\epsilon^2}{16},
\end{align*}
where the first inequality is because $F(x,D_a)$ is non-decreasing as $x\geq x^*_{D_a}$, the second inequality is derived by the $1/2$-strong convexity of $ F(x,D_a)$, and the last inequality is due to $\epsilon\in(0,1/4)$.

By the above inequality, we obtain 
\begin{align} \label{eq:str-lbthm1_eq_2}
    \sum_{t=1}^{\Delta_T}\E^{\pi}_{D_a}[F(x_t^\pi(\mathcal H_{t-1}),D_a) - F(x^*_{D_a},D_a)]\geq \frac{\epsilon^2}{16} \sum_{t=1}^{\Delta_T}\P^{\pi}_{D_a}[x_t^\pi(\mathcal H_{t-1}) > 1].
\end{align}
Similarly, we can prove that
\begin{align} \label{eq:str-lbthm1_eq_3}
    \sum_{t=1}^{\Delta_T}\E^{\pi}_{D_b}[F(x_t^\pi(\mathcal H_{t-1}),D_b) - F(x^*_{D_b},D_b)]\geq \frac{\epsilon^2}{16} \sum_{t=1}^{\Delta_T}\P^{\pi}_{D_b}[x_t^\pi(\mathcal H_{t-1}) \leq 1].
\end{align}
By Bretagnolle–Huber inequality,  it holds that 
\begin{align} \label{eq:str-lbthm1_eq_4}
    \P^{\pi}_{D_a}[x_t^\pi(\mathcal H_{t-1}) > 1] + \P^{\pi}_{D_b}[x_t^\pi(\mathcal H_{t-1}) \leq 1]&\geq \frac{1}{2}
   \exp\left(-\mathrm{KL}(\P_{D_a}^{(t-1)}, \P_{D_b}^{(t-1)})\right),
\end{align}
where $\mathrm{KL}(\P_{D_a}^{(t-1)}, \P_{D_b}^{(t-1)})$ is the Kullback-Leibler divergence divergence between the distribution of $\{D_1, \dots, D_t\}$ under $D_a$ and $D_b$.
By the chain rule of the Kullback-Leibler divergence, we have 
\begin{align}\nonumber
    \mathrm{KL}(\P_{D_a}^{(t-1)}, \P_{D_b}^{(t-1)}) &= \E_{D_a}\left[\log\left(\frac{f_{D_a}(X)}{f_{D_b}(X)}\right)\right]
\\& =(t-1)\left[(1/2+\epsilon) \log \left(\frac{1/2+\epsilon}{1/2-\epsilon}\right)+(1/2-\epsilon) \log \left(\frac{1/2-\epsilon}{1/2+\epsilon}\right)\right]\nonumber
\\& \leq (t-1)\left[(1/2+\epsilon)\left(4 \epsilon+8 \epsilon^2\right)-4(1/2-\epsilon) \epsilon\right] \nonumber
\\& \leq {16 (t-1) \epsilon^2}, \label{eq:str-lbthm1_eq_5}
\end{align}
where the first inequality is due to $2 x \leq \log ((1+x)/(1-x)) \leq 2 x+2 x^2$ as $x \in(0,1 / 2)$, and the second inequality is due to $\epsilon\in(0,1/4)$.

Combining Eqs.~(\ref{eq:str-lbthm1_eq_1}, \ref{eq:str-lbthm1_eq_2}, \ref{eq:str-lbthm1_eq_3}, \ref{eq:str-lbthm1_eq_4}, \ref{eq:str-lbthm1_eq_5}), we obtain
\begin{align}\label{eq:str-lbthm1_eq_6}
    \sup_{D\in\{D_a,D_b\}}\left\{\sum_{t=1}^{\Delta_T}\E^{\pi}_D[F(x_t^\pi(\mathcal H_{t-1}),D) - F(x^*_{D},D)] \right\}\geq \frac{\epsilon}{64} \sum_{t=1}^{\Delta_T}\exp\left(-16 (t-1) \epsilon^2 \right).
\end{align}

Since $\mathrm{TV}(D_a, D_b)=4\epsilon$,
setting $\Delta_T = \lceil (T/V)^{2/3}\rceil$  and $\epsilon = 1/(4\sqrt{\Delta_T})$,
 we can verify that 
\[\sum_{t=2}^{T}\|D_{t}-D_{t-1}\|_{\mathrm{TV}}\leq V,\]
and
\begin{align*}
 \sup_{\bm{D}\in\mathcal D_{V}^{(1)}}\left\{\sum_{t=1}^{T}\E^{\pi}_D[F(x_t^\pi(\mathcal H_{t-1}),D) - F(x^*_{D},D)] \right\}&\geq \frac{\epsilon^2}{64} \left\lceil\frac{T}{\Delta_T}\right\rceil  \sum_{t=1}^{\Delta_T}\exp\left(-16 (t-1) \epsilon^2 \right)\\
    &\geq \frac{\epsilon^2}{64} \left\lceil\frac{T}{\Delta_T}\right\rceil \sum_{t=1}^{\min\{\epsilon^{-2},\Delta_T\}}\exp(-16)\\
     &= \frac{\epsilon^2}{64} \left\lceil\frac{T}{\Delta_T}\right\rceil \cdot \epsilon^{-2}\exp(-16)\\
    &\geq \frac{\exp(-16)}{128} V^{2/3}T^{1/3},
\end{align*}
where we complete the proof.
\Halmos
\endproof

\section{Proofs Omitted in Section \ref{sec:censored}}
\label{sec: Proofs Omitted in sec:censored}
Recall that the Dvoretzky–Kiefer–Wolfowitz (DKW) inequality states that the concentration phenomenon of empirical distribution functions around the true distribution function for one-dimensional independent random variables. Applying the DKW inequality, we can prove the good event that all marginal CDFs we construct are concentrated around the true CDFs.
\begin{lemma}
\label{lem:good-event}
For any input $\delta \in (0,1)$, we have the following event
\[
\G = \left\{ \sup_{y \in \R} \vert \hat{G}_{s,t}(y) - G_{s,t}(y)\vert \leq \sqrt{\frac{\ln (2T^2/\delta)}{t-s+1}},~\forall 1\leq s\leq t \leq T, \text{ and } i \in [m] \right\}
\]
happens with probability at least $1-\delta$.
\end{lemma}
\subsection{Proof of Lemma~\ref{lem:good-event}}

\proof{Proof of Lemma~\ref{lem:good-event}.}
From the DKW inequality (Lemma \ref{lem:DKW}), we know for any $1 \leq s \leq t \leq T$ and $i \in [m]$,
\begin{align*}
\P\left[ \sup_{y \in \R} \vert \hat{G}_{s,t}(y) - G_{s,t}(y)\vert > \sqrt{\frac{\ln (2T^2/\delta)}{t-s+1}}\right] &\leq 2\exp\left(-2(t-s+1)\cdot \frac{\ln (2T^2/\delta)}{t-s+1} \right) \\
&\leq 2\exp(-\ln (2T^2/\delta)) \leq \frac{\delta}{mT^2}.
\end{align*}
Note that the set $\{(s,t,i)|1 \leq s \leq t \leq T, i \in [m]\}$ has at most  $mT^2$ pairs. Therefore, we know
\begin{align*}
    \P[\G] = 1 - \P[\G^\complement] &\geq 1 - \sum_{1 \leq s \leq t \leq T} \sum_{i=1}^m \P\left[ \sup_{y \in \R} \vert \hat{G}_{s,t}(y) - G_{s,t}(y)\vert > \sqrt{\frac{\ln (2T^2/\delta)}{t-s+1}}\right]\\
    &\geq 1 - mT^2 \cdot \frac{\delta}{mT^2} = 1-\delta,
\end{align*}
where we complete the proof.
\Halmos\endproof

\subsection{Proof of Theorem \ref{thm:censored-nvp-ub-general}.}
\label{subsec:Proof of Theorem thm:censored-nvp-ub-general}

\begin{lemma}\label{le:censoredlemmaS}
    Under the implementation of Algorithm~\ref{alg:saacensored}, the number of epochs $\tau_{max}\leq \min\{S+1, \tau_{max}\leq 2V^{2/3}T^{1/3}+1\}$. 
\end{lemma}

\subsubsection{Proof of  $ {\O}(\sqrt{ST\log (T/\delta)} )$. }
By the definition of $S$, there exists a disjoint interval partition of time horizon, i.e.,~ $[1,T]=\mathcal{I'}_1 \cup \mathcal{I'}_2 \cup \cdots \cup \mathcal{I'}_{S+1}$, such that in each interval $\mathcal{I}_k', k \in [S]$, random variables $\{D_i, i \in \mathcal{I'}_k\}$ have the same distribution.
For each  $\mathcal{E}_{\tau}=[l_{\tau}, l_{\tau+1})$, we project the partition of the whole time horizon  $[1,T]=\mathcal{I'}_1 \cup \mathcal{I'}_2 \cup \cdots \cup \mathcal{I'}_{S+1}$ on  $\mathcal{E}_{\tau}$ and obtain the following partition of 
\[
\mathcal{E}_{\tau} = \mathcal{I}_1 \cup \mathcal{I}_2 \cup \dots \cup \mathcal{I}_{K_{\tau}}.
\]
By the definition of the partition and projection, it is easy to note that for each $k\in[K_\tau]$,~$ S_{\mathcal I_k}=0$.

Consider a fixed $t \in [l_{\tau}+1,l_{\tau+1}-1)$, let $t \in \cI_k$ for some $k$, and let $s_k$ be the starting period of $\cI_k$.
By the definition of the epoch, we know that in period $t$, the distribution is detected to be stationary. It follows that,
\[
\sup_{y\in[0, x_t]} |\hat G_{l_{\tau},t-1}(y) - \hat G_{s_{k},t}(y)| \leq 2\sqrt{\frac{\ln (2T^2/\delta)}{t-l_{\tau}}} + 2\sqrt{\frac{\ln (2T^2/\delta)}{t-s_k+1}},
\]
which further implies
\begin{equation}
\label{eq:f-detect-error-censored}
|\hat{g}_{l_{\tau},t-1}(x) - \hat{g}_{s_k,t}(x)| \leq 2(h+b)\left(\sqrt{\frac{\ln (2T^2/\delta)}{t-l_{\tau}}} + \sqrt{\frac{\ln (2T^2/\delta)}{t-s_k+1}}\right),
\end{equation}
for all $x\in\X_t$ by    $ \hat{g}_{i,j}(x) := 
     (h+b) \cdot \hat{G}_{i,j}(x) - b$. 

Under the event $\G$, we have
\[
\sup_{y\in\mathbb{R}} |\hat G_{s_{k},t}(y) -G_{s_{k},t}(y)| \leq \sqrt{\frac{\ln (2T^2/\delta)}{t-s_k+1}}. 
\]
Thus, we know for any $x \in \X_t$,
\begin{equation}
\label{eq:f-con-error2-censored}
|\hat g_{s_{k},t}(x) -g_{s_{k},t}(x)| \leq (h+b)\sqrt{\frac{\ln (2T^2/\delta)}{t-s_k+1}}.
\end{equation}
From the definition of $s_k$, we know $g_t(x) = g_{s_k,t}(x)$.
Combining with Eq.~\eqref{eq:f-detect-error-censored} and Eq.~\eqref{eq:f-con-error2-censored}, we know
\begin{align}\nonumber
    &|\hat{g}_{l_{\tau},t-1}(x) - g_t(x) | \\\nonumber
&\leq |\hat{g}_{l_{\tau},t-1}(x) - \hat{g}_{s_k,t}(x) | + |\hat{g}_{s_k,t}(x) - {g}_{s_k,t}(x) | + |{g}_{s_k,t}(x) - g_t(x) | \\
&\leq (h+b)\left(2\sqrt{\frac{\ln (2T^2/\delta)}{t-l_{\tau}}} + 3\sqrt{\frac{\ln (2T^2/\delta)}{t-s_k+1}}\right).\label{eq:thmcensoredeq1}
\end{align}
Since $t \in [l_\tau+1, l_{\tau+1} - 1)$,  by the epoch-restarting conditions in Algorithm~\ref{alg:saacensored}, $\mathcal{X}_t  \neq \emptyset$. Thus, at the beginning of each period $t$, the set $\mathcal{X}_t$ is non-empty. The non-emptiness of $\mathcal{X}_t$ guarantees the well-definedness of $x_t$.
Note that $x_t\in \mathcal X_t$, which means that 
\[
\hat{g}_{l_{\tau},t-1}(x_t)\leq 2(h+b)\sqrt{\frac{\ln (2T^2/\delta)}{t-l_{\tau}}}.
\]

Combining the above two inequalities, we have 
\begin{align}\nonumber
| g_t(x_t)|  &\leq | \hat{g}_{l_{\tau},t-1}(x_t)| +| g_t(x_t) - \hat{g}_{l_{\tau},t-1}(x_t) |
\\\label{eq:thmgtxtbound-censor}
&\leq  (h+b)\left(4\sqrt{\frac{\ln (2T^2/\delta)}{t-l_{\tau}}} + 3\sqrt{\frac{\ln (2T^2/\delta)}{t-s_k+1}}\right),   
\end{align}
By the convexity of $f_t(x)$ and boundedness of the domain, we have 
\begin{align}
    f_t (x_t) - f_t(x_t^*) & \leq (x_t-x_t^*)^\top\nabla f_t(x_t)
    \leq \bar x \Vert g_t(x_t) \Vert,\label{eq:thm_convex_eq_1-censored}
\end{align}
where $x_t^*=\argmin_{x\in[0,x_t]}f_t(x)$.
Combining Eqs.~(\ref{eq:thmgtxtbound-censor},\ref{eq:thm_convex_eq_1-censored}), we obtain
\begin{align}
\nonumber    f_t(x_t) - f_t(x^*_t)  \leq  \bar x(h+b)\left(4\sqrt{\frac{\ln 
 (2T^2/\delta)}{t-l_{\tau}}} +3\sqrt{\frac{\ln (2T^2/\delta)}{t-s_k+1}}\right).\label{eq:ft-reg-censored}
\end{align}
By that $\sum_{n=1}^N 1/\sqrt{n} \leq 2\sqrt{N} $, we know
\begin{align*}
    \sum_{t \in \cI_i} \sqrt{\frac{\ln (2T^2/\delta)}{t-s_k+1}} \leq 2\sqrt{|\cI_i|\ln (2T^2/\delta)}\text{~~~and~~~}
    \sum_{t \in \cE_{\tau}} \sqrt{\frac{\ln 
 (2T^2/\delta)}{t-l_{\tau}}} \leq 2\sqrt{|\cE_{\tau}|\ln (2T^2/\delta)}.
\end{align*}
Combining the above inequality with the boundedness of $f_t(\cdot)\leq \bar x(h+b)$, we have
\begin{align*}
       & \sum_{t \in \mathcal{E}_{\tau}} f_t (x_t) - f_t(x_t^*)\\&\leq f_{l_{\tau}}(x_{l_{\tau}})-f_{l_{\tau}}(x_{l_{\tau}}^*) + f_{l_{\tau+1}-1}(x_{l_{\tau+1}-1})-f_{l_{\tau+1}-1}(x_{l_{\tau+1}-1}^*) + \sum_{t= l_{\tau}+1}^{l_{\tau+1}-2}   f_t (x_t) - f_t(x_t^*) \\
    &\leq 2\bar x(h+b)( 1+ 4 \sqrt{|\cE_{\tau}|\ln (2 T^2/\delta)} + 3\sum_{i=1}^{K_{\tau}} \sqrt{|\cI_{i}|\ln (2 T^2/\delta)})\\
    &\leq 8\bar x(h+b)( \sqrt{|\cE_{\tau}|\ln (2 T^2/\delta)} + \sum_{i=1}^{K_{\tau}} \sqrt{|\cI_{i}|\ln (2 T^2/\delta)}).
\end{align*}
By Lemma~\ref{le:censoredlemmaS}, we have $\tau_{max}\leq S+1$. 
Moreover, due to restarts, we may divide a sequence of stationary variables into at most two different epochs. We know that
\[
\sum_{\tau=1}^{\tau_{max}} \sum_{i=1}^{K_{\tau}} 1 \leq 2S +1.
\]
By Cauchy–Schwarz inequality and the bound of regret incurred in each epoch, we know that with probability at least $1-\delta$, 
\begin{align*}
    \sum_{t=1}^T f_t (x_t) - f_t(x_t^*)=  \sum_{\tau=1}^{\tau_{max}} \sum_{t \in \mathcal{E}_{\tau}} f_t (x_t) - f_t(x_t^*) &\leq  32 \bar x(h+b) \sqrt{(S+1)T\ln (2 T^2/\delta)}.
\end{align*}
where we complete the proof.
\Halmos
\endproof

\subsubsection{Proof of  ${\O}( V^{1/3}T^{2/3}\log^{1/2}T  )$ }
For each  $\mathcal{E}_{\tau}=[l_{\tau}, l_{\tau+1})$, we project the partition of the whole time horizon constructed in Lemma~\ref{le:partition_V} on  $\mathcal{E}_{\tau}$ and obtain the following partition of 
\[
\mathcal{E}_{\tau} = \mathcal{I}_1 \cup \mathcal{I}_2 \cup \dots \cup \mathcal{I}_{K_{\tau}}.
\]
By Lemma~\ref{le:partition_V}, it is easy to note that
\begin{align}\label{eq:thm_general_V_eq_1-censoredv}
    V_{\mathcal I_k}\leq \sqrt{\frac{1}{|\mathcal I_k|}},~~\forall k\in[K_{\tau}].
\end{align}

Consider a fixed $t \in [l_{\tau}+1,l_{\tau+1}-1)$, let $t \in \cI_k$ for some $k$, and let $s_k$ be the starting period of $\cI_k$.
By the definition of the epoch, we know that in period $t$, the distribution is detected to be stationary. It follows that,
\[
\sup_{y\in[0, x_t]} |\hat G_{l_{\tau},t-1}(y) - \hat G_{s_{k},t}(y)| \leq 2\sqrt{\frac{\ln (2T^2/\delta)}{t-l_{\tau}}} + 2\sqrt{\frac{\ln (2T^2/\delta)}{t-s_k+1}},
\]
which further implies
\begin{equation}
\label{eq:f-detect-error-censoredv}
|\hat{g}_{l_{\tau},t-1}(x) - \hat{g}_{s_k,t}(x)| \leq 2(h+b)\left(\sqrt{\frac{\ln (2T^2/\delta)}{t-l_{\tau}}} + \sqrt{\frac{\ln (2T^2/\delta)}{t-s_k+1}}\right),
\end{equation}
for all $x\in\X_t$ by    $ \hat{g}_{i,j}(x) := 
     (h+b) \cdot \hat{G}_{i,j}(x) - b$. 

Under the event $\G$, we have
\[
\sup_{y\in\mathbb{R}} |\hat G_{s_{k},t}(y) -G_{s_{k},t}(y)| \leq \sqrt{\frac{\ln (2T^2/\delta)}{t-s_k+1}}. 
\]
Thus, we know for any $x \in \X_t$,
\begin{equation}
\label{eq:f-con-error2-censoredv}
|\hat g_{s_{k},t}(x) -g_{s_{k},t}(x)| \leq (h+b)\sqrt{\frac{\ln (2T^2/\delta)}{t-s_k+1}}.
\end{equation}
From the definition of $s_k$ and Eq.~\eqref{eq:thm_general_V_eq_1}, we know
\[|g_t(x)-g_{s_k,t}(x)|\leq (h+b) \sqrt{\frac{1}{|\mathcal I_k|}} \leq (h+b) \sqrt{\frac{\ln (2T^2/\delta)}{t-s_k+1}}. \]
Combining with Eq.~\eqref{eq:f-detect-error-censoredv} and Eq.~\eqref{eq:f-con-error2-censoredv}, we know
\begin{align}\nonumber
    &|\hat{g}_{l_{\tau},t-1}(x) - g_t(x) | \\\nonumber
&\leq |\hat{g}_{l_{\tau},t-1}(x) - \hat{g}_{s_k,t}(x) | + |\hat{g}_{s_k,t}(x) - {g}_{s_k,t}(x) | + |{g}_{s_k,t}(x) - g_t(x) | \\
&\leq (h+b)\left(2\sqrt{\frac{\ln (2T^2/\delta)}{t-l_{\tau}}} + 4\sqrt{\frac{\ln (2T^2/\delta)}{t-s_k+1}}\right).\label{eq:thmcensoredeq1-censoredv}
\end{align}
Since $t \in [l_\tau+1, l_{\tau+1} - 1)$,  by the epoch-restarting conditions in Algorithm~\ref{alg:saacensored}, $\mathcal{X}_t  \neq \emptyset$. Thus, at the beginning of each period $t$, the set $\mathcal{X}_t$ is non-empty. The non-emptiness of $\mathcal{X}_t$ guarantees the well-definedness of $x_t$.
Note that $x_t\in \mathcal X_t$, which means that \[\hat{g}_{l_{\tau},t-1}(x_t)\leq 2(h+b)\sqrt{\frac{\ln (2T^2/\delta)}{t-l_{\tau}}}\]
Combining the above two inequalities, we have 
\begin{align}\nonumber
| g_t(x_t)|  &\leq | \hat{g}_{l_{\tau},t-1}(x_t)| +| g_t(x_t) - \hat{g}_{l_{\tau},t-1}(x_t) |
\\\label{eq:thmgtxtbound-censoredv}
&\leq  (h+b)\left(4\sqrt{\frac{\ln (2T^2/\delta)}{t-l_{\tau}}} + 4\sqrt{\frac{\ln (2T^2/\delta)}{t-s_k+1}}\right),   
\end{align}
By the convexity of $f_t(x)=h\cdot\E\left[(x-D_t)^+\right]+b\cdot\E\left[(D_t-x)^+\right]$ and boundedness of the domain, we have 
\begin{align}
    f_t (x_t) - f_t(x_t^*) & \leq (x_t-x_t^*)^\top\nabla f_t(x_t)
    \leq \bar x \Vert g_t(x_t) \Vert,\label{eq:thm_convex_eq_1-censoredv}
\end{align}
where $x_t^*=\argmin_{x\in[0,x_t]}f_t(x)$.
Combining Eqs.~(\ref{eq:thmgtxtbound-censoredv},\ref{eq:thm_convex_eq_1-censoredv}), we obtain
\begin{align}
\nonumber    f_t(x_t) - f_t(x^*_t)  \leq  \bar x(h+b)\left(4\sqrt{\frac{\ln 
 (2T^2/\delta)}{t-l_{\tau}}} +4\sqrt{\frac{\ln (2T^2/\delta)}{t-s_k+1}}\right).\label{eq:ft-reg-censoredv}
\end{align}
By that $\sum_{n=1}^N 1/\sqrt{n} \leq 2\sqrt{N} $, we know
\begin{align*}
    \sum_{t \in \cI_i} \sqrt{\frac{\ln (2T^2/\delta)}{t-s_k+1}} \leq 2\sqrt{|\cI_i|\ln (2T^2/\delta)}\text{~~~and~~~}
    \sum_{t \in \cE_{\tau}} \sqrt{\frac{\ln 
 (2T^2/\delta)}{t-l_{\tau}}} \leq 2\sqrt{|\cE_{\tau}|\ln (2T^2/\delta)}.
\end{align*}
Combining the above inequality with the boundedness of $f_t(\cdot)\leq \bar x(h+b)$, we have
\begin{align*}
       & \sum_{t \in \mathcal{E}_{\tau}} f_t (x_t) - f_t(x_t^*)\\&\leq f_{l_{\tau}}(x_{l_{\tau}})-f_{l_{\tau}}(x_{l_{\tau}}^*) + f_{l_{\tau+1}-1}(x_{l_{\tau+1}-1})-f_{l_{\tau+1}-1}(x_{l_{\tau+1}-1}^*) + \sum_{t= l_{\tau}+1}^{l_{\tau+1}-2}   f_t (x_t) - f_t(x_t^*) \\
    &\leq 2\bar x(h+b)( 1+ 4 \sqrt{|\cE_{\tau}|\ln (2 T^2/\delta)} + 4\sum_{i=1}^{K_{\tau}} \sqrt{|\cI_{i}|\ln (2 T^2/\delta)})\\
    &\leq 10\bar x(h+b)( \sqrt{|\cE_{\tau}|\ln (2 T^2/\delta)} + \sum_{i=1}^{K_{\tau}} \sqrt{|\cI_{i}|\ln (2 T^2/\delta)}).
\end{align*}
By Lemma \ref{le:censoredlemmaS}, we know that  $\tau_{max}\leq 2 V^{2/3}T^{1/3}+1:=K$.
Moreover, due to restarts, we may divide an interval $\mathcal I'_{k'}$ into at most two different intervals $\mathcal I_{j}$ and $\mathcal I_{j'}$. It follows that 
\[
\sum_{\tau=1}^{\tau_{max}} \sum_{i=1}^{K_{\tau}} 1 \leq 2K +1.
\]

By Cauchy–Schwarz inequality, we know that
\begin{align*}
    \sum_{\tau =1}^{\tau_{max}} \sqrt{|\mathcal E_{\tau}|} \leq 2\sqrt{KT} \leq 4(V^{1/3}+1)T^{2/3},\\
    \sum_{\tau =1}^{\tau_{max}} \sum_{i =1}^{K_{\tau}}\sqrt{|\mathcal I_{i}|} \leq 2\sqrt{(2K+1)T} \leq 4(V^{1/3}+1)T^{2/3},
\end{align*}
Combining the above two equations with the bound of regret incurred in each epoch, we know that with probability at least $1-\delta$, 
\begin{align*}
     \sum_{t=1}^T f_t (x_t) - f_t(x_t^*)=  \sum_{\tau=1}^{\tau_{max}} \sum_{t \in \mathcal{E}_{\tau}} f_t (x_t) - f_t(x_t^*) &\leq 80\bar x(h+b) (V^{1/3}+1)T^{2/3}\ln^{1/2}(2T^2/\delta),
\end{align*}
where we complete the proof.
\Halmos
\endproof

\subsubsection{Proof of Lemma~\ref{le:censoredlemmaS}.}

In the following, we prove $\tau_{max}\leq S+1$ and $\tau_{max}\leq 2V^{2/3}T^{1/3}+1$ individually.

\textit{Proof of $\tau_{max}\leq S+1$.}
By the definition of $S$, there exists a disjoint interval partition of time horizon, i.e.,~ $[1,T]=\mathcal{I'}_1 \cup \mathcal{I'}_2 \cup \cdots \cup \mathcal{I'}_{S+1}$, such that in each interval $\mathcal{I'}_{k},~ k \in [S+1]$, random variables $\{D_i, i \in \mathcal{I'}_k\}$ have the same distribution.
Recall that $\mathcal{E}_{\tau}= [l_{\tau},l_{\tau+1})$ be the set of periods of epoch $\tau$.
\begin{claim}
    We claim that for any $k'\in [S+1]$, the interval $\mathcal I'_{k'}$ could include at most one end point of epochs, i.e., if $l_{\tau}-1\in[\mathcal I'_{k'}]$, then $l_{\tau-1}-1\notin[\mathcal I'_{k'}]$ and $l_{\tau+1}-1\notin[\mathcal I'_{k'}]$.
\end{claim}
 With this claim, we could simply conclude that $\tau_{max} \leq S+1$. Now we prove this claim by contradiction. Without loss of generality, we assume both $l_{\tau}-1$ and $l_{\tau+1}-1$ are contained in the interval $\mathcal I'_{k'}$.

We first show that $\mathcal X_{l_{\tau+1}}$ is non-empty. 
Under the event $\G$, we have
\[
\sup_{y\in\mathbb{R}} |\hat G_{l_{\tau},t}(y) -G_{l_{\tau},t}(y)| \leq \sqrt{\frac{\ln (2T^2/\delta)}{t-l_{\tau}+1}}. 
\]
By     $ \hat{g}_{i,j}(x) := 
     (h+b) \cdot \hat{G}_{i,j}(x) - b$ and the assumption that $l_{\tau}-1\in \mathcal I'_{k'}$ and $l_{\tau+1}-1\in \mathcal I'_{k'}$, we have  for any $t\in[l_{\tau},l_{\tau+1}-1]$ we have $g_{l_{\tau}}(x)={g}_{l_{\tau},t}(x)$  and thus it holds that
\begin{align*}
    |\hat{g}_{l_{\tau},t}(x) - g_{l_{\tau}}(x) | 
&= |\hat{g}_{l_{\tau},t}(x) - {g}_{l_{\tau},t}(x) |\\
&\leq (h+b)\sqrt{\frac{\ln (2T^2/\delta)}{t-l_{\tau}+1}}.
\end{align*}
Let $x_{l_\tau}^* = \argmin_{x\in[0,\bar x]}f_{l_\tau}(x)$, we know that $g_{l_\tau}(x_{l_\tau}^*)=\cdots=g_t(x_{l_\tau}^*)=0$ for any $t\leq l_{\tau+1}-1$.
Combining the above inequality, for any $t\in[l_{\tau},l_{\tau+1}-1]$ we have
\begin{align*}
    |\hat{g}_{l_{\tau},t}(x_{l_\tau}^*)  | \leq  (h+b)\sqrt{\frac{\ln (2T^2/\delta)}{t-l_{\tau}+1}}.
\end{align*}
Thus, $x_{l_\tau}^*$ won't be eliminated for any $t\in[l_{\tau},l_{\tau+1}-1]$.  It holds that $x_{l_\tau}^*\in \mathcal X_{l_{\tau+1}}$ and consequently $\mathcal X_{l_{\tau+1}}$ is non-empty.

Then we show that the nonstationary detection won't be triggered.
Under $\mathcal{G}$,  we know for any $s \in [l_{\tau},t]$, 
\begin{align*}
    \sup_{y \in \R} \vert \hat{G}_{l_{\tau},t-1}(y) - {G}_{l_{\tau},t-1}(y) \vert \leq \sqrt{\frac{\ln (2T^2/\delta)}{t-l_{\tau}}},\text{~~and~~}
    \sup_{y \in \R} \vert \hat{G}_{s,t}(y) - {G}_{s,t}(y) \vert \leq \sqrt{\frac{\ln (2T^2/\delta)}{t-s+1}}.
\end{align*}
Since $S_{[l_{\tau}-1,t]} = 0$, we know that $\sup_{y \in \R}\Vert G_{l_{\tau},t}(y) - G_{s,t}(y) \Vert=0$,
\[
\sup_{y \in \R}\vert \hat G_{l_{\tau},t-1}(y)- \hat G_{s,t}(y) \vert \leq  \sqrt{\frac{\ln (2T^2/\delta)}{t-l_{\tau}}}+ \sqrt{\frac{\ln (2T^2/\delta)}{t-s+1}}.
\]
Thus, we know that Algorithm \ref{alg:saacensored} will not restart a new epoch $\tau+1$ in period $t=l_{\tau+1}-1$, where we get a contradiction and prove the claim.

\textit{Proof of $\tau_{max}\leq 2V^{2/3}T^{1/3}+1$.}
By Lemma~\ref{le:partition_V},
There exists a disjoint interval partition of time horizon, i.e.,~ $[1,T]=\mathcal{I'}_1 \cup \mathcal{I'}_2 \cup \cdots \cup \mathcal{I'}_{K}$, such that $V_{\mathcal I'_k}\leq \sqrt{1/{|\mathcal I'_k|}}$, $\forall k\in[K]$, and $K\leq 2V^{2/3}T^{1/3}+1$.

Recall that $\mathcal{E}_{\tau}= [l_{\tau},l_{\tau+1})$ be the set of periods of epoch $\tau$.
\begin{claim}
    We claim that for any $k'\in [K]$, the interval $\mathcal I'_{k'}$ could include at most one end point of epochs, i.e., if $l_{\tau}-1\in[\mathcal I'_{k'}]$, then $l_{\tau-1}-1\notin[\mathcal I'_{k'}]$ and $l_{\tau+1}-1\notin[\mathcal I'_{k'}]$.
\end{claim}
 With this claim, we could simply conclude that $\tau_{max} \leq K$. Now we prove this claim by contradiction. Without loss of generality, we assume both $l_{\tau}-1$ and $l_{\tau+1}-1$ are contained in the interval $\mathcal I'_{k'}$.

We first show that $\mathcal X_{l_{\tau+1}}$ is non-empty. 
By $V_{\mathcal I'_{k'}}\leq \sqrt{1/{|\mathcal I'_{k'}|}}$ and the fact that $t\in \mathcal I'_{k'}$, we know that $\sup_{y \in R}\Vert G_{l_{\tau},t}(y)-G_{l_{\tau}}(y)\Vert \leq \sqrt{1/|\mathcal I'_{k'}|}$, which further implies that
\begin{equation}
\label{eq:g-con-nonstationoary_V-censored}
|{g}_{l_{\tau},t}(x) - g_{l_{\tau}}(x) |\leq (h+b)\sqrt{\frac{1}{|\mathcal I_k|}} \leq (h+b)\sqrt{\frac{\ln (2T^2/\delta)}{t-l_{\tau}+1}},
\end{equation}
Under the event $\G$, we have
\[
\sup_{y\in\mathbb{R}} |\hat G_{l_{\tau},t}(y) -G_{l_{\tau},t}(y)| \leq \sqrt{\frac{\ln (2T^2/\delta)}{t-l_{\tau}+1}}. 
\]
By     $ \hat{g}_{i,j}(x) := 
     (h+b) \cdot \hat{G}_{i,j}(x) - b$ and the assumption that $l_{\tau}-1\in \mathcal I'_{k'}$ and $l_{\tau+1}-1\in \mathcal I'_{k'}$, we have  for any $t\in[l_{\tau},l_{\tau+1}-1]$ it holds that
\begin{align*}
    &|\hat{g}_{l_{\tau},t}(x) - g_{l_{\tau}}(x) | \\
&\leq |\hat{g}_{l_{\tau},t}(x) - \hat{g}_{l_{\tau},t}(x) | +  |{g}_{l_{\tau},t}(x) - g_{l_{\tau}}(x) | \\
&\leq 2(h+b)\sqrt{\frac{\ln (2T^2/\delta)}{t-l_{\tau}+1}}.
\end{align*}

Let $x_{l_\tau}^* = \argmin_{x\in[0,\bar x]}f_{l_\tau}(x)$, we know that $g_{l_\tau}(x_{l_\tau}^*)=0$ for any $t\leq l_\tau-1$.
Combining the above inequality, for any $t\in[l_{\tau},l_{\tau+1}-1]$ we have
\begin{align*}
    |\hat{g}_{l_{\tau},t}(x_{l_\tau}^*)  | \leq  2(h+b)\sqrt{\frac{\ln (2T^2/\delta)}{t-l_{\tau}+1}}.
\end{align*}
Thus, $x_{l_\tau}^*$ won't be eliminated for any $t\in[l_{\tau},l_{\tau+1}-1]$.  It holds that $x_{l_\tau}^*\in \mathcal X_{l_{\tau+1}}$ and consequently $\mathcal X_{l_{\tau+1}}$ is non-empty.

Then we show that the nonstationary detection won't be triggered.
Under $\mathcal{G}$,  we know for any $s \in [l_{\tau},t]$, 
\begin{align*}
    \sup_{y \in \R} \vert \hat{G}_{l_{\tau},t-1}(y) - {G}_{l_{\tau},t-1}(y) \vert \leq \sqrt{\frac{\ln (2T^2/\delta)}{t-l_{\tau}}},\text{~~and~~}
    \sup_{y \in \R} \vert \hat{G}_{s,t}(y) - {G}_{s,t}(y) \vert \leq \sqrt{\frac{\ln (2T^2/\delta)}{t-s+1}}.
\end{align*}
Due to $t \in \mathcal I'_{k'}$, we know that $V_{[l_{\tau},t]} \leq \sqrt{1/|\mathcal I'_{k'}|}$. Therefore, we know that $\sup_{y \in \R}\vert G_{l_{\tau},t-1}(y) - G_{s,t}(y) \vert\leq \sqrt{1/|\mathcal I'_{k'}|} $ and 
\begin{align*}
\sup_{y \in \R}\vert\hat{G}_{l_{\tau},t-1}(y) - \hat G_{s,t}(y) \vert &\leq 
\sup_{y \in \R}\vert \hat{G}_{l_{\tau},t-1}(y) - {G}_{l_{\tau},t-1}(y) \vert + \sup_{y \in \R} \vert \hat{G}_{s,t}(y) - {G}_{s,t}(y) \vert + \sup_{y \in \R}\vert G_{l_{\tau},t-1}(y) - G_{s,t}(y) \vert\\
&\leq  \sqrt{\frac{\ln (2T^2/\delta)}{t-l_{\tau}}}+ \sqrt{\frac{\ln (2T^2/\delta)}{t-s+1}} +\sqrt{\frac{1}{|\mathcal I'_{k'}|}}\\
&\leq 2\sqrt{\frac{\ln (2T^2/\delta)}{t-l_{\tau}}}+ 2\sqrt{\frac{\ln (2T^2/\delta)}{t-s+1}}.
\end{align*}
Thus, we know that Algorithm \ref{alg:saacensored} will not restart a new epoch $\tau+1$ in period $t=l_{\tau+1}-1$, where we get a contradiction and prove the claim.
\Halmos

\subsection{Proof of Theorem \ref{thm:censored-nvp-ub-pl}.}
\label{subsec:Proof of Theorem thm:censored-nvp-ub-pl}
\subsubsection{Proof of $ \O(S\log(T/S)\log T)$.}

Similar to the proof of $\O(\sqrt{ST\log T})$ in Theorem \ref{thm:censored-nvp-ub-general}, we can prove that
\begin{align*}
    |\hat{g}_{l_{\tau},t-1}(x) - g_t(x) | 
\leq (h+b)\left(2\sqrt{\frac{\ln (2T^2/\delta)}{t-l_{\tau}}} + 3\sqrt{\frac{\ln (2T^2/\delta)}{t-s_k+1}}\right).
\end{align*}
Note that $x_t\in \mathcal X_t$, which means that \[\hat{g}_{l_{\tau},t-1}(x_t)\leq 2(h+b)\sqrt{\frac{\ln (2T^2/\delta)}{t-l_{\tau}}}\]
Combining the above two inequalities, we have 
\begin{align}\nonumber
| g_t(x_t)|  &\leq | \hat{g}_{l_{\tau},t-1}(x_t)| +| g_t(x_t) - \hat{g}_{l_{\tau},t-1}(x_t) |
\\\label{eq:thmgtxtbound-censor-svx}
&\leq  (h+b)\left(4\sqrt{\frac{\ln (2T^2/\delta)}{t-l_{\tau}}} + 3\sqrt{\frac{\ln (2T^2/\delta)}{t-s_k+1}}\right),   
\end{align}
By the $\alpha$-global minimal separation condition, we know that 
\[f_t''(x) = (h+b)G_t'(x)\geq \alpha(h+b).\] Thus we know that $f_t(x)$ is $\alpha(h+b)$-strongly convex, and the following inequality holds
\begin{align}
    f_t (x_t) - f_t(x_t^*)
    \leq \frac{1}{2\alpha(h+b)} | g_t(x_t) |^2,\label{eq:thm_convex_eq_1-censored-svx}
\end{align}
where $x_t^*=\argmin_{x\in[0,x_t]}f_t(x)$.
Combining Eqs.~(\ref{eq:thmgtxtbound-censor-svx},\ref{eq:thm_convex_eq_1-censored-svx}), we obtain
\begin{align*}
    f_t\left({x}_t\right)-f_t(x^*_t) \leq\frac{h+b}{2\alpha}  \left(\frac{16\ln (2T^2/\delta)}{t-l_{\tau}} + \frac{9\ln (2T^2/\delta)}{t - s_k+1} \right).
\end{align*}
Similar to the proof of $\O(\sqrt{ST\log T})$ in Theorem \ref{thm:censored-nvp-ub-general}, by some simple calculation, we know that with probability at least $1-\delta$,
\begin{align*}
    \sum_{t=1}^T f_t (x_t) - f_t(x_t^*) \leq \frac{64(h+b)}{\alpha}(S+1)\ln (2T^2/\delta)(1+\ln (T/S)),
\end{align*}
where we complete the proof.
\Halmos
\endproof
\subsubsection{Proof of $ \O(V^{2/3}T^{1/3}\log^{3/2}T)$.}
\label{subsec:Proof of  O(V23T13log32T)}
Similar to the proof of $ \O(S\log(T/S)\log T)$ above and proof of $ \O(V^{1/3}T^{2/3}\log^{1/2}T)$ in Theorem~\ref{thm:censored-nvp-ub-general}, 
\begin{align*}
    f_t\left({x}_t\right)-f_t(x^*_t) \leq\frac{h+b}{2\alpha}  \left(\frac{16\ln (2T^2/\delta)}{t-l_{\tau}} + \frac{16\ln (2T^2/\delta)}{t - s_k+1} \right).
\end{align*}
Similar to the proof of Theorem~\ref{thm:censored-nvp-ub-general}, we have
\begin{align*}
    \sum_{t=1}^T f_t (x_t) - f_t(x_t^*) &\leq  \frac{128(h+b)}{\alpha} (V^{2/3}T^{1/3}+1)\ln T\ln^{1/2}(2T^2/\delta),
\end{align*}

\section{Proof Omitted in Section \ref{sec:extension}}\label{appendix:extension-proof}

In this section, all the arguments are conditioned on the good event $\mathcal{G}$ defined in Lemma \ref{lem:good-event}.
\subsection{Proof of Theorem \ref{thm:general-ub}}
\label{sec:Proof of Theorem thm:uncensored-nvp-ub-general}
We first present the following critical technical lemma that bounds the number of epochs, and its proof is presented in Section~\ref{sec:proof of lemma uncensored lemma}.
\begin{lemma}\label{le:uncensoredlemmaS}
    Under the implementation of Algorithm~\ref{alg:saa-uncensored-nvp}, the number of epochs $\tau_{max}\leq \min\{S+1, \tau_{max}\leq 2V^{2/3}T^{1/3}+1\}$. 
\end{lemma}

\subsubsection{Proof of  ${\O}(\sqrt{ST\log T} )$ }

Let $\mathcal{E}_{\tau}= [l_{\tau},l_{\tau+1})$ be the set of periods of epoch $\tau$.
For the convenience of presentation, we define the projection of a time interval $[t_a, t_b]$ on epoch $\mathcal{E}_{\tau}=[l_{\tau}, l_{\tau+1})$ as 
\begin{equation}
    \left\{
\begin{aligned}
&\emptyset, ~~&& t_b< l_\tau,\\
   &[l_\tau, t_b], ~~&& t_a< l_\tau,~ l_\tau\leq t_b\leq l_{\tau+1}-1,\\
       &[t_a, t_b], ~~&& t_a\ge l_\tau,~l_\tau\leq t_b\leq l_{\tau+1}-1,\\
          &[t_a, l_{\tau+1}-1], ~~&&  l_\tau\leq t_a\leq l_{\tau+1}-1,~ t_b> l_{\tau+1}-1,\\
          &\emptyset, ~~&& t_a> l_{\tau+1}-1.\\
\end{aligned}
    \right.
\end{equation}
By the definition of $S$, there exists a disjoint interval partition of time horizon, i.e.,~ $[1,T]=\mathcal{I'}_1 \cup \mathcal{I'}_2 \cup \cdots \cup \mathcal{I'}_{S+1}$, such that in each interval $\mathcal{I}_k', k \in [S]$, random variables $\{D_i, i \in \mathcal{I'}_k\}$ have the same distribution.
For each  $\mathcal{E}_{\tau}=[l_{\tau}, l_{\tau+1})$, we project the partition of the whole time horizon  $[1,T]=\mathcal{I'}_1 \cup \mathcal{I'}_2 \cup \cdots \cup \mathcal{I'}_{S+1}$ on  $\mathcal{E}_{\tau}$ and obtain the following partition of 
\[
\mathcal{E}_{\tau} = \mathcal{I}_1 \cup \mathcal{I}_2 \cup \dots \cup \mathcal{I}_{K_{\tau}}.
\]
By the definition of the partition and projection, it is easy to note that for each $k\in[K_\tau]$,$ S_{\mathcal I_k}=0$.

Consider a fixed $t \in [l_{\tau}+1,l_{\tau+1}-1)$, let $t \in \cI_k$ for some $k$, and let $s_k$ be the starting period of $\cI_k$.
From the definition of optimization oracle and Algorithm \ref{alg:saa-uncensored-nvp}, we have
\begin{equation}
\label{eq:opt-f-error}
\hat{f}_{l_{\tau},t-1}(x_t) - \hat{f}_{l_{\tau},t-1}(\hat{x}^*_{l_{\tau},t-1}) \leq \epsilon_t,
\end{equation}
where $\hat{x}^*_{l_{\tau},t-1}$ is the optimal solution of $\min_{x \in \cX} \hat{f}_{l_{\tau},t-1}(x)$.

By the definition of epoch, we know that in period $t$, the distribution is detected to be stationary. It follows that
{\[
\sup_{y\in\mathbb R}| \hat G_{l_{\tau},t-1}(y) - \hat G_{s_{k},t}(y)|\leq 2\sqrt{\frac{\ln (2T^2/\delta)}{t-l_{\tau}}} + 2\sqrt{\frac{\ln (2T^2/\delta)}{t-s_k+1}},
\]}
which further implies
\begin{equation}
\label{eq:f-detect-error}
|\hat{f}_{l_{\tau},t-1}(x) - \hat{f}_{s_k,t}(x)| \leq 2L_F\left(\sqrt{\frac{\ln (2T^2/\delta)}{t-l_{\tau}}} + \sqrt{\frac{\ln (2T^2/\delta)}{t-s_k+1}}\right),
\end{equation}
for all $x$ by Assumption \ref{ass:DKW-to-error}.

{Under the event $\G$, we have 
\[
\sup_{y\in\mathbb R}|\hat G_{s_{k},t}(y) -G_{s_{k},t}(y)| \leq \sqrt{\frac{\ln (2T^2/\delta)}{t-s_k+1}}. 
\]}
By Assumption \ref{ass:DKW-to-error} and the above inequality, we know for any $x \in \X$,
\begin{equation}
\label{eq:f-con-error}
|\hat f_{s_{k},t}(x) -f_{s_{k},t}(x)| \leq L_F\sqrt{\frac{\ln (2T^2/\delta)}{t-s_k+1}}.
\end{equation}

From the definition of $s_k$, we know $f_t(x) = f_{s_k,t}(x)$.
Combining with Eq.~\eqref{eq:f-detect-error} and Eq.~\eqref{eq:f-con-error}, we know
\begin{align*}
    &|\hat{f}_{l_{\tau},t-1}(x_t) - f_t(x_t) | \\
&\leq |\hat{f}_{l_{\tau},t-1}(x_t) - \hat{f}_{s_k,t}(x_t) | + |\hat{f}_{s_k,t}(x_t) - {f}_{s_k,t}(x_t) | + |{f}_{s_k,t}(x_t) - f_t(x_t) | \\
&\leq L_F\left(2\sqrt{\frac{\ln (2T^2/\delta)}{t-l_{\tau}}} + 3\sqrt{\frac{\ln (2T^2/\delta)}{t-s_k+1}}\right).
\end{align*}

Similarly, we can prove that
\[
|\hat{f}_{l_{\tau},t-1}(x^*_t) - f_t(x^*_t) | \leq  L_F\left(2\sqrt{\frac{\ln (2T^2/\delta)}{t-l_{\tau}}} + 3\sqrt{\frac{\ln (2T^2/\delta)}{t-s_k+1}}\right).
\]

Therefore, we obtain
\begin{align}
\nonumber    f_t(x_t) - f_t(x^*_t) &\leq \hat{f}_{l_{\tau},t-1}(x_t) - \hat{f}_{l_{\tau},t-1}(x^*_t) \\ \nonumber &\phantom{\leq}+ L_F\left(4\sqrt{\frac{\ln 
 (2T^2/\delta)}{t-l_{\tau}}} +6\sqrt{\frac{\ln (2T^2/\delta)}{t-s_k+1}}\right) \\
 &\leq \epsilon_t + L_F\left(4\sqrt{\frac{\ln 
 (2T^2/\delta)}{t-l_{\tau}}} +6\sqrt{\frac{\ln (2T^2/\delta)}{t-s_k+1}}\right), \label{eq:ft-reg}
\end{align}
where the last inequality is by Eq.~\eqref{eq:opt-f-error}.

By that $\sum_{n=1}^N 1/\sqrt{n} \leq 2\sqrt{N} $, we know
\begin{align*}
    \sum_{t \in \cI_i} \sqrt{\frac{\ln (2T^2/\delta)}{t-s_k+1}} &\leq 2\sqrt{|\cI_i|\ln (2T^2/\delta)},\\
    \sum_{t= l_{\tau}+1}^{l_{\tau+1}-2} \sqrt{\frac{\ln 
 (2T^2/\delta)}{t-l_{\tau}}} &\leq 2\sqrt{|\cE_{\tau}|\ln (2T^2/\delta)}.
\end{align*}

Therefore, we have
\begin{align*}
  & \sum_{t \in \mathcal{E}_{\tau}} f_t (x_t) - f_t(x_t^*)\\&\leq f_{l_{\tau}}(x_{l_{\tau}})-f_{l_{\tau}}(x_{l_{\tau}}^*) + f_{l_{\tau+1}-1}(x_{l_{\tau+1}-1})-f_{l_{\tau+1}-1}(x_{l_{\tau+1}-1}^*) + \sum_{t= l_{\tau}+1}^{l_{\tau+1}-2}   f_t (x_t) - f_t(x_t^*) \\
    &\leq 2 C_f + \sum_{t= l_{\tau}+1}^{l_{\tau+1}-2} \epsilon_t+ 8 L_F\sqrt{|\cE_{\tau}|\ln (2 T^2/\delta)}  + \sum_{i=1}^{K_{\tau}} 12 L_F\sqrt{|\cI_{i}|\ln (2 T^2/\delta)}.
\end{align*}
By Lemma~\ref{le:uncensoredlemmaS}, we know that $\tau_{max}\leq S+1$.
Moreover, due to restarts, we may divide a sequence of stationary variables into at most two different epochs. We know that
\[
\sum_{\tau=1}^{\tau_{max}} \sum_{i=1}^{K_{\tau}} 1 \leq 2S +1.
\]

By Cauchy–Schwarz inequality and the bound of regret incurred in each epoch, we know that with probability at least $1-\delta$, 
\begin{align*}
    \sum_{t=1}^T f_t (x_t) - f_t(x_t^*) &\leq 40 L_F \sqrt{(S+1)T\ln (2 T^2/\delta)}  + \sum_{\tau=1}^{\tau_{max}} \sum_{t \in \mathcal{E}_{\tau}}\epsilon_t +2C_f (S+1),
\end{align*}
where we complete the proof.
\Halmos\endproof

\subsubsection{Proof of  ${\O}( V^{1/3}T^{2/3}\log^{1/2}(T)  )$}\label{subsec:proof-of-V-general-ub}

We first present the following technical lemma about the partition of time horizon.
\begin{lemma}\label{le:partition_V}
There exists a disjoint interval partition of time horizon, i.e.,~ $[1,T]=\mathcal{I'}_1 \cup \mathcal{I'}_2 \cup \cdots \cup \mathcal{I'}_{K}$, such that $V_{\mathcal I'_k}\leq \sqrt{1/{|\mathcal I'_k|}}$, $\forall k\in[K]$, and $K\leq 2V^{2/3}T^{1/3}+1$.
\end{lemma}

For each  $\mathcal{E}_{\tau}=[l_{\tau}, l_{\tau+1})$, we project the partition of the whole time horizon constructed in Lemma~\ref{le:partition_V} on  $\mathcal{E}_{\tau}$ and obtain the following partition of 
\[
\mathcal{E}_{\tau} = \mathcal{I}_1 \cup \mathcal{I}_2 \cup \dots \cup \mathcal{I}_{K_{\tau}}.
\]
By Lemma~\ref{le:partition_V}, it is easy to note that
\begin{align}\label{eq:thm_general_V_eq_1}
    V_{\mathcal I_k}\leq \sqrt{\frac{1}{|\mathcal I_k|}},~~\forall k\in[K_{\tau}].
\end{align}

Since in period $t\in [l_{\tau}+1,l_{\tau+1}-1)$ the nonstationarity detection fails, denoting $s_k$ as the starting period of $\cI_k, k \in [K_{\tau}]$, we have 
\[
\sup_{y\in\mathbb R}| \hat G_{l_{\tau},t-1}(y) - \hat G_{s_{k},t}(y)| \leq 2\sqrt{\frac{\ln (2T^2/\delta)}{t-l_{\tau}}} + 2\sqrt{\frac{\ln (2T^2/\delta)}{t-s_k+1}},
\]
which further implies
\begin{equation}
\label{eq:f-detect-error_V}
|\hat{f}_{l_{\tau},t-1}(x) - \hat{f}_{s_k,t}(x)| \leq  2L_F\left(\sqrt{\frac{\ln (2T^2/\delta)}{t-l_{\tau}}} + \sqrt{\frac{\ln (2T^2/\delta)}{t-s_k+1}}\right),
\end{equation}
for all $x$ by Assumption \ref{ass:DKW-to-error}.

{Under the event $\G$, we have
\[
\sup_{y\in\mathbb R}|\hat G_{s_{k},t}(y) -G_{s_{k},t}(y)| \leq \sqrt{\frac{\ln (2T^2/\delta)}{t-s_k+1}}. 
\]}
By Assumption \ref{ass:DKW-to-error} and the above inequality, we know for any $x \in \X$,
\begin{equation}
\label{eq:f-con-error_V}
|\hat f_{s_{k},t}(x) -f_{s_{k},t}(x)| \leq L_F\sqrt{\frac{\ln (2T^2/\delta)}{t-s_k+1}}.
\end{equation}
By Eq.~\eqref{eq:thm_general_V_eq_1} and the fact that $t\in \mathcal I_k$, we know that $\sup_{y \in \R}\vert G_{s_k,t}(y)-G_{t}(y)\vert \leq \sqrt{1/|\mathcal I_k|}$, which further implies that
\begin{equation}
\label{eq:f-con-nonstationoary_V}
|{f}_{s_k,t}(x_t) - f_t(x_t) |\leq L_F\sqrt{\frac{1}{|\mathcal I_k|}} \leq L_F\sqrt{\frac{\ln (2T^2/\delta)}{t-s_k+1}},
\end{equation}
where the last inequality is due to $s_k, t \in \mathcal I_k$ and $T \geq 2$.

Combining Eqs.~(\ref{eq:f-detect-error_V},\ref{eq:f-con-error_V},\ref{eq:f-con-nonstationoary_V}), we know
\begin{align*}
    &|\hat{f}_{l_{\tau},t-1}(x_t) - f_t(x_t) | \\
&\leq |\hat{f}_{l_{\tau},t-1}(x_t) - \hat{f}_{s_k,t}(x_t) | + |\hat{f}_{s_k,t}(x_t) - {f}_{s_k,t}(x_t) | + |{f}_{s_k,t}(x_t) - f_t(x_t) | \\
&\leq L_F\left(2\sqrt{\frac{\ln (2T^2/\delta)}{t-l_{\tau}}} + 4\sqrt{\frac{\ln (2T^2/\delta)}{t-s_k+1}}\right).
\end{align*}

Similarly, we can prove that
\[
|\hat{f}_{l_{\tau},t-1}(x^*_t) - f_t(x^*_t) | \leq L_F\left(2\sqrt{\frac{\ln (2T^2/\delta)}{t-l_{\tau}}} + 4\sqrt{\frac{\ln (2T^2/\delta)}{t-s_k+1}}\right).
\]

From the definition of optimization oracle and Algorithm \ref{alg:saa-uncensored-nvp}, we have
\begin{equation}
\label{eq:opt-f-error_V}
\hat{f}_{l_{\tau},t-1}(x_t) - \hat{f}_{l_{\tau},t-1}(\hat{x}^*_{l_{\tau},t-1}) \leq \epsilon_t,
\end{equation}
where $\hat{x}^*_{l_{\tau},t-1}$ is the optimal solution of $\min_{x \in \cX} \hat{f}_{l_{\tau},t-1}(x)$.

Therefore, we obtain
\begin{align}
\nonumber    f_t(x_t) - f_t(x^*_t) &\leq \hat{f}_{l_{\tau}-1,t-1}(x_t) - \hat{f}_{l_{\tau}-1,t-1}(x^*_t) \\ \nonumber &\phantom{\leq}+ L_F\left(4\sqrt{\frac{\ln 
 (2T^2/\delta)}{t-l_{\tau}}} +8\sqrt{\frac{\ln (2T^2/\delta)}{t-s_k+1}}\right) \\
 &\leq \epsilon_t + L_F\left(4\sqrt{\frac{\ln 
 (2T^2/\delta)}{t-l_{\tau}}} +8\sqrt{\frac{\ln (2T^2/\delta)}{t-s_k+1}}\right), \label{eq:ft-reg_V}
\end{align}
where the last inequality is by Eq.~\eqref{eq:opt-f-error_V}.

In the interval $\cI_i$, we have
\[
\sum_{t \in \cI_i} \sqrt{\frac{\ln (2T^2/\delta)}{t-s_k+1}} \leq 2\sqrt{|\cI_i|\ln (2T^2/\delta)}.
\]
In the epoch $\cE_{\tau}$, we have
\[
\sum_{t= l_{\tau}+1}^{l_{\tau+1}-2} \sqrt{\frac{\ln 
 (2T^2/\delta)}{t-l_{\tau}}} \leq 2\sqrt{|\cE_{\tau}|\ln (2T^2/\delta)}.
\]
Therefore, we have
\begin{align*}
  & \sum_{t \in \mathcal{E}_{\tau}} f_t (x_t) - f_t(x_t^*)\\&\leq f_{l_{\tau}}(x_{l_{\tau}})-f_{l_{\tau}}(x_{l_{\tau}}^*) + f_{l_{\tau+1}-1}(x_{l_{\tau+1}-1})-f_{l_{\tau+1}-1}(x_{l_{\tau+1}-1}^*) + \sum_{t= l_{\tau}+1}^{l_{\tau+1}-2}   f_t (x_t) - f_t(x_t^*) \\
    &\leq 2 C_f + \sum_{t= l_{\tau}+1}^{l_{\tau+1}-2} \epsilon_t+ 8 L_F\sqrt{|\cE_{\tau}|\ln (2 T^2/\delta)}  + \sum_{i=1}^{K_{\tau}} 16 L_F\sqrt{|\cI_{i}|\ln (2 T^2/\delta)}.
\end{align*}
By Lemma~\ref{le:uncensoredlemmaS}, we know that $\tau_{max}\leq 2V^{2/3}T^{1/3}+1:=K$.
Moreover, due to restarts, we may divide an interval $\mathcal I'_{k'}$ into at most two different intervals $\mathcal I_{j}$ and $\mathcal I_{j'}$. It follows that 
\[
\sum_{\tau=1}^{\tau_{max}} \sum_{i=1}^{K_{\tau}} 1 \leq 2K +1.
\]
By Cauchy–Schwarz inequality, we know that
\begin{align*}
    \sum_{\tau =1}^{\tau_{max}} \sqrt{|\mathcal E_{\tau}|} \leq 2\sqrt{KT} \leq 4(V^{1/3}+1)T^{2/3},\\
    \sum_{\tau =1}^{\tau_{max}} \sum_{i =1}^{K_{\tau}}\sqrt{|\mathcal I_{i}|} \leq 2\sqrt{(2K+1)T} \leq 4(V^{1/3}+1)T^{2/3},
\end{align*}
Combining the above two equations with the bound of regret incurred in each epoch, we know that with probability at least $1-\delta$, 
\begin{align*}
    \sum_{t=1}^T f_t (x_t) - f_t(x_t^*) &\leq 4C_f(V^{2/3}T^{1/3}+1) + \sum_{\tau =1}^{\tau_{max}} \sum_{t \in \mathcal E_{\tau}} \epsilon_t + 96L_F (V^{1/3}+1)T^{2/3}\ln^{1/2}(2T^2/\delta),
\end{align*}
where we complete the proof.
\Halmos\endproof
\subsubsection{Proof of Lemma~\ref{le:uncensoredlemmaS}.}
\label{sec:proof of lemma uncensored lemma}

In the following, we prove $\tau_{max}\leq S+1$ and $\tau_{max}\leq 2V^{2/3}T^{1/3}+1$ individually.

\textit{Proof of $\tau_{max}\leq S+1$.}
By the definition of $S$, there exists a disjoint interval partition of time horizon, i.e.,~ $[1,T]=\mathcal{I'}_1 \cup \mathcal{I'}_2 \cup \cdots \cup \mathcal{I'}_{S+1}$, such that in each interval $\mathcal{I'}_{k},~ k \in [S+1]$, random variables $\{D_i, i \in \mathcal{I'}_k\}$ have the same distribution.
Recall that $\mathcal{E}_{\tau}= [l_{\tau},l_{\tau+1})$ be the set of periods of epoch $\tau$.
\begin{claim}
    We claim that for any $k'\in [S+1]$, the interval $\mathcal I'_{k'}$ could include at most one end point of epochs, i.e., if $l_{\tau}-1\in[\mathcal I'_{k'}]$, then $l_{\tau-1}-1\notin[\mathcal I'_{k'}]$ and $l_{\tau+1}-1\notin[\mathcal I'_{k'}]$.
\end{claim}
 With this claim, we could simply conclude that $\tau_{max} \leq S+1$. Now we prove this claim by contradiction. Without loss of generality, we assume both $l_{\tau}-1$ and $l_{\tau+1}-1$ are contained in the interval $\mathcal I'_{k'}$.
To prove this claim, we only need to show that the nonstationary detection won't be triggered.
Under $\mathcal{G}$,  we know for any $s \in [l_{\tau},t]$, 
\begin{align*}
    \sup_{y \in \R} \vert \hat{G}_{l_{\tau},t-1}(y) - {G}_{l_{\tau},t-1}(y) \vert \leq \sqrt{\frac{\ln (2T^2/\delta)}{t-l_{\tau}}},\text{~~and~~}
    \sup_{y \in \R} \vert \hat{G}_{s,t}(y) - {G}_{s,t}(y) \vert \leq \sqrt{\frac{\ln (2T^2/\delta)}{t-s+1}}.
\end{align*}
Since $S_{[l_{\tau}-1,t]} = 0$, we know that $\sup_{y \in \R}\vert G_{l_{\tau},t-1}(y) - G_{s,t}(y) \vert=0$,
\[
\sup_{y \in \R}\vert \hat G_{l_{\tau},t-1}(y)- \hat G_{s,t}(y) \vert \leq  \sqrt{\frac{\ln (2T^2/\delta)}{t-l_{\tau}}}+ \sqrt{\frac{\ln (2T^2/\delta)}{t-s+1}}.
\]
Thus, we know that Algorithm \ref{alg:saa-uncensored-nvp} will not restart a new epoch $\tau+1$ in period $t=l_{\tau+1}-1$, where we get a contradiction and prove the claim.

\textit{Proof of $\tau_{max}\leq 2V^{2/3}T^{1/3}+1$.}
By Lemma~\ref{le:partition_V},
There exists a disjoint interval partition of time horizon, i.e.,~ $[1,T]=\mathcal{I'}_1 \cup \mathcal{I'}_2 \cup \cdots \cup \mathcal{I'}_{K}$, such that $V_{\mathcal I'_k}\leq \sqrt{1/{|\mathcal I'_k|}}$, $\forall k\in[K]$, and $K\leq 2V^{2/3}T^{1/3}+1$.

Recall that $\mathcal{E}_{\tau}= [l_{\tau},l_{\tau+1})$ be the set of periods of epoch $\tau$.
\begin{claim}
    We claim that for any $k'\in [K]$, the interval $\mathcal I'_{k'}$ could include at most one end point of epochs, i.e., if $l_{\tau}-1\in[\mathcal I'_{k'}]$, then $l_{\tau-1}-1\notin[\mathcal I'_{k'}]$ and $l_{\tau+1}-1\notin[\mathcal I'_{k'}]$.
\end{claim}
 With this claim, we could simply conclude that $\tau_{max} \leq K$. Now we prove this claim by contradiction. Without loss of generality, we assume both $l_{\tau}-1$ and $l_{\tau+1}-1$ are contained in the interval $\mathcal I'_{k'}$.

To prove this lemma, we only need to show that the nonstationary detection won't be triggered.
Under $\mathcal{G}$,  we know for any $s \in [l_{\tau}-1,t]$, 
\begin{align*}
    \sup_{y \in \R} \vert \hat{G}_{l_{\tau},t-1}(y) - {G}_{l_{\tau},t-1}(y) \vert \leq \sqrt{\frac{\ln (2T^2/\delta)}{t-l_{\tau}}},\text{~~and~~}
    \sup_{y \in \R} \vert \hat{G}_{s,t}(y) - {G}_{s,t}(y) \vert \leq \sqrt{\frac{\ln (2T^2/\delta)}{t-s+1}}.
\end{align*}
Due to $t \in \mathcal I'_{k'}$, we know that $V_{[l_{\tau}-1,t]} \leq \sqrt{1/|\mathcal I'_{k'}|}$. Therefore, we know that $\sup_{y \in \R}\vert G_{l_{\tau},t-1}(y) - G_{s,t}(y) \vert\leq \sqrt{1/|\mathcal I'_{k'}|} $ and 
\begin{align*}
\sup_{y \in \R}\vert\hat{G}_{l_{\tau},t-1}(y) - \hat G_{s,t}(y) \vert &\leq 
\sup_{y \in \R}\vert \hat{G}_{l_{\tau},t-1}(y) - {G}_{l_{\tau},t-1}(y) \vert \\
&\phantom{+}+ \sup_{y \in \R} \vert \hat{G}_{s,t}(y) - {G}_{s,t}(y) \vert + \sup_{y \in \R}\vert G_{l_{\tau},t-1}(y) - G_{s,t}(y) \vert\\
&\leq  \sqrt{\frac{\ln (2T^2/\delta)}{t-l_{\tau}}}+ \sqrt{\frac{\ln (2T^2/\delta)}{t-s+1}} +\sqrt{\frac{1}{|\mathcal I'_{k'}|}}\\
&\leq 2\sqrt{\frac{\ln (2T^2/\delta)}{t-l_{\tau}}}+ 2\sqrt{\frac{\ln (2T^2/\delta)}{t-s+1}}.
\end{align*}
Thus, we know that Algorithm \ref{alg:saa-uncensored-nvp} will not restart a new epoch $\tau+1$ in period $t=l_{\tau+1}-1$, where we get a contradiction and prove the claim.
\Halmos

\subsubsection{Proof of Lemma~\ref{le:partition_V}.} 
The proof idea of this lemma is from \cite{wei2021non}. 

Consider the constructive partitioning procedure presented in Algorithm~\ref{alg:partition1}. Invoking Algorithm~\ref{alg:partition1} with input $[1,T]$, we can get a partition $\mathcal{I}'_{1} \cup \mathcal{I}'_{2}\cup \cdots \cup  \mathcal{I}'_{K}$, where $K$ is the number of intervals in the partition. By construction of intervals (Line~\ref{line:ifcondition1}) of Algorithm~\ref{alg:partition1}, it is easy to note that $\forall k\in[K]$,
\[V_{\mathcal I'_k} \leq \sqrt{\frac{1}{|\mathcal I'_k|}}.\]  

\begin{algorithm}[h!]
\caption{Existence of Partition in  Lemma~\ref{le:partition_V}}
\label{alg:partition1}
\begin{algorithmic}[1]
\State \textbf{Initialization:} Input an interval $\mathcal{E}=[s, e]$. Let $k=1$, $s_1=s$, $t=s$.
\While{$t\leq e$}
    \If{$V_{[s_k,t]} \leq \sqrt{\frac{1}{t-s_k+1}}$ and $V_{[s_k,t+1]} > \sqrt{\frac{1}{t-s_k+2}}$}\label{line:ifcondition1}
        \State Let $e_k\leftarrow t$, $\mathcal{I}'_k\leftarrow [s_k,e_k]$
        \State $k\leftarrow k+1$, $s_k\leftarrow t+1$.
    \EndIf
    \State \textbf{end if}
    \State $t\leftarrow t+1$.
\EndWhile
\State \textbf{end while}
\If{$s_k\leq e$}
    \State $e_k\leftarrow e$, $\mathcal{I}'_{k}\leftarrow [s_k,e_k]$.
\EndIf
 \State \textbf{end if}
\State \textbf{Output:} partition $\mathcal{I}'_{1} \cup \mathcal{I}'_{2}\cup \cdots \cup  \mathcal{I}'_{k}$, $K=k$. 
\end{algorithmic}
\end{algorithm}

Next, it suffices to provide an upper bound on $K$. For $K=1$, the conclusion holds obviously. For $K>1$, by the condition in Line~\ref{line:ifcondition1} and the definition of $V_\mathcal I$, we obtain
\begin{align}
V_{[s, e]} \nonumber
&\geq V_{[s_1, e_1+1]}+V_{[s_2, e_2+1]}+\dots+V_{[s_{K-1}, e_{K-1}+1]}  \\&\geq \sum_{k=1}^{K-1} \sqrt{\frac{1}{e_k-s_k+2}}\nonumber \\&\geq \sum_{k=1}^{K-1} \sqrt{\frac{1}{|\mathcal I'_k|+1}}.  \label{eq:le_partition_eq_1}
\end{align}
By H\"{o}lder's inequality, it holds that
\begin{align}\nonumber
    (K-1)&=  \sum_{k=1}^{K-1}\frac{1}{(|\mathcal I'_k|+1)^{1/3}} (|\mathcal I'_k|+1)^{1/3}
    \\&\leq \left( \sum_{k=1}^{K-1} \sqrt{\frac{1}{|\mathcal I'_k|+1}} \right)^{\frac{2}{3}} \left(\sum_{k=1}^{K-1} (|\mathcal I'_k|+1) \right)^{\frac{1}{3}}\label{eq:le_partition_eq_2}
\end{align}
Combining Eq.\eqref{eq:le_partition_eq_1} and  Eq.\eqref{eq:le_partition_eq_2}, we get 
\begin{align*}
    K-1&\leq \left(\sum_{k=1}^{K-1} (|\mathcal I'_k|+1) \right)^{\frac{1}{3}}V_{[1, T]}^{2/3}\\& \leq (T+K-1)^{1/3}V^{2/3} \leq 2V^{2/3}T^{1/3},
\end{align*}
where we complete the proof of Lemma~\ref{le:partition_V}. 
\Halmos
\endproof
\subsection{Proof of Theorem \ref{thm:pl-ub}.}
\label{subsec:Proof of Theorem thm:uncensored-nvp-ub-pl}
\subsubsection{Proof of $ \O(S\log(T/S)\log T)$.}
\label{subsec:Proof of  O(Slog(TS)log T)}

Similar to the proof of Theorem \ref{thm:uncensored-nvp-ub-general}, we can prove that
\[
\Vert\hat{g}_{l_{\tau},t-1}(x_t) - g_t(x_t) \Vert_2 \leq L_g\left(2\sqrt{\frac{\ln (2T^2/\delta)}{t-l_{\tau}}} + 3\sqrt{\frac{\ln (2T^2/\delta)}{t-s_k+1}}\right).
\]
The above inequality implies that
\begin{align*}
 \Vert g_t(x_t)  \Vert_2  &\leq  \Vert \hat{g}_{l_{\tau},t-1}(x_t)  \Vert_2 +\Vert g_t(x_t) - \hat{g}_{l_{\tau},t-1}(x_t) \Vert_2 
\\
&\leq \epsilon_t + L_g\left(2\sqrt{\frac{\ln (2T^2/\delta)}{t-l_{\tau}}} + 3\sqrt{\frac{\ln (2T^2/\delta)}{t-s_k+1}}\right),    
\end{align*}
where $ \Vert \hat{g}_{l_{\tau},t-1}(x_t)  \Vert_2$ is bounded by the definition of $x_t$ and the optimization oracle.

Combining the above inequality with Assumption \ref{ass:pl}, we have
\begin{align*}
    f_t\left({x}_t\right)-f_t(x^*_t) \leq 2\mu \epsilon^2_t + 4\mu L_g^2 \left(\frac{4\ln (2T^2/\delta)}{t-l_{\tau}} + \frac{9\ln (2T^2/\delta)}{t - s_k+1} \right).
\end{align*}
Similar to the proof of Theorem \ref{thm:censored-nvp-ub-general}, by some simple calculation, we obtain
\begin{align*}
    \sum_{t=1}^T f_t (x_t) - f_t(x_t^*) \leq 104\mu L_g^2(S+1)\ln (2T^2/\delta)(1+\ln (T/S))+ \sum_{t=1}^T 2\mu \epsilon^2_t + 2C_f(S+1),
\end{align*}
Thus, we know that with probability at least $1-\delta$, 
\begin{align*}
    \sum_{t=1}^T f_t (x_t) - f_t(x_t^*) &\leq 104\mu L_g^2(S+1)\ln (2T^2/\delta)(1+\ln (T/S))+ \sum_{t=1}^T 2\mu \epsilon^2_t + 2C_f(S+1),
\end{align*}
where we complete the proof.
\Halmos
\endproof
\subsubsection{Proof of $ \O(V^{2/3}T^{1/3}\log^{3/2}T)$.}
\label{subsec:proof-of-V-pl-ub}
Similar to the proof of Theorem~\ref{thm:uncensored-nvp-ub-general}, for each  $\mathcal{E}_{\tau}=[l_{\tau}, l_{\tau+1})$, we project the partition of the whole time horizon constructed in Lemma~\ref{le:partition_V} on  $\mathcal{E}_{\tau}$ and obtain the following partition of 
\[
\mathcal{E}_{\tau} = \mathcal{I}_1 \cup \mathcal{I}_2 \cup \dots \cup \mathcal{I}_{K_{\tau}}.
\]
By Lemma~\ref{le:partition_V}, it is easy to note that
\begin{align}\label{eq:thm_general_V_eq_1(1)}
    V_{\mathcal I_k}\leq \sqrt{\frac{1}{|\mathcal I_k|}},~~\forall k\in[K_{\tau}].
\end{align}

Since in period $t\in [l_{\tau}+1,l_{\tau+1}-1)$ the nonstationarity detection fails, denoting $s_k$ as the starting period of $\cI_k, k \in [K_{\tau}]$, we have
\[
\sup_{y \in \R}\vert \hat G_{l_{\tau},t-1}(y) - \hat G_{s_{k},t}(y)\vert \leq 2\sqrt{\frac{\ln (2T^2/\delta)}{t-l_{\tau}}} + 2\sqrt{\frac{\ln (2T^2/\delta)}{t-s_k+1}},
\]
which further implies
\begin{equation}
\label{eq:g-detect-error_V}
\Vert \hat{g}_{l_{\tau},t-1}(x) - \hat{g}_{s_k,t}(x)\Vert_2 \leq  L_g\left(2\sqrt{\frac{\ln (2T^2/\delta)}{t-l_{\tau}}} + 2\sqrt{\frac{\ln (2T^2/\delta)}{t-s_k+1}}\right),
\end{equation}
for all $x$ by Assumption \ref{ass:DKW-to-error2}.

Under the event $\G$, we have
\[
\sup_{y \in \R}|\hat G_{s_{k},t}(y) -G_{s_{k},t}(y)| \leq \sqrt{\frac{\ln (2T^2/\delta)}{t-s_k+1}}. 
\]
By Assumption \ref{ass:DKW-to-error2} and the above inequality, we know for any $x \in \X$,
\begin{equation}
\label{eq:g-con-error_V}
\Vert \hat g_{s_{k},t}(x) -g_{s_{k},t}(x) \Vert_2 \leq L_g\sqrt{\frac{\ln (2T^2/\delta)}{t-s_k+1}}.
\end{equation}
By Eq.~\eqref{eq:thm_general_V_eq_1(1)} and the fact that $t\in \mathcal I_k$, we know that $\sup_{y \in \R}\vert G_{s_k,t}(y)-G_{t}(y)\vert \leq \sqrt{1/|\mathcal I_k|}$, which further implies that
\begin{equation}
\label{eq:g-con-nonstationoary_V}
\Vert {g}_{s_k,t}(x_t) - g_t(x_t) \Vert_2\leq L_g\sqrt{\frac{1}{|\mathcal I_k|}} \leq L_g\sqrt{\frac{\ln (2T^2/\delta)}{t-s_k+1}},
\end{equation}
where the last inequality is due to $s_k, t \in \mathcal I_k$ and $T \geq 2$.

It follows that
\begin{equation}
\label{eq:thm_str_eq_4}
\begin{split}
\nonumber
    \Vert g_t(x_t)- \hat{g}_{l_{\tau},t-1}(x_t)\Vert_2 &\leq \underbrace{\Vert g_t(x_t)-{g}_{s_k,t}(x_t) \Vert_2}_{\mathrm{term~I}}+\underbrace{\Vert {g}_{s_k,t}(x_t)-\hat{g}_{s_k,t}(x_t) \Vert_2}_{\mathrm{term~II}} \\&~~~~~+ \underbrace{\Vert \hat{g}_{s_k,t}(x_t) -\hat{g}_{l_{\tau},t-1}(x_t)  \Vert_2}_{\mathrm{term~III}} \\
    &\leq L_g\left(2\sqrt{\frac{\ln (2T^2/\delta)}{t-l_{\tau}}} + 4\sqrt{\frac{\ln (2T^2/\delta)}{t-s_k+1}}\right),
\end{split}
\end{equation}
where term I is upper bounded by Eq.~\eqref{eq:g-con-nonstationoary_V}, term II is bounded by Eq.~\eqref{eq:g-con-error_V}, term III is bounded by Eq.~\eqref{eq:g-detect-error_V}.

The above inequality implies that
\begin{align*}
 \Vert g_t(x_t)  \Vert_2  &\leq  \Vert \hat{g}_{l_{\tau},t-1}(x_t)  \Vert_2 +\Vert g_t(x_t) - \hat{g}_{l_{\tau},t-1}(x_t) \Vert_2 
\\
&\leq \epsilon_t + L_g\left(2\sqrt{\frac{\ln (2T^2/\delta)}{t-l_{\tau}}} + 4\sqrt{\frac{\ln (2T^2/\delta)}{t-s_k+1}}\right),    
\end{align*}
where $ \Vert \hat{g}_{l_{\tau},t-1}(x_t)  \Vert_2$ is bounded by the definition of $x_t$ and the optimization oracle.

Combining the above inequality with Assumption~\ref{ass:pl}, we have
\begin{align*}
    f_t\left({x}_t\right)-f_t(x^*_t) \leq 2\mu \epsilon^2_t + 16\mu L_g^2 \left(\frac{\ln (2T^2/\delta)}{t-l_{\tau}} + \frac{4\ln (2T^2/\delta)}{t - s_k+1} \right).
\end{align*}

Note that $\sum_{t=1}^n 1/t \leq 1+ \ln n$. In the interval $\cI_i$, we have
\[
\sum_{t \in \cI_i} {\frac{\ln (2T^2/\delta)}{t-s_k+1}} \leq \ln (2T^2/\delta)(1+\ln |\cI_i|) \leq \ln (2T^2/\delta)(1+\ln T)
\]
In the epoch $\cE_{\tau}$, we have
\[
\sum_{t= l_{\tau}+1}^{l_{\tau+1}-2}  {\frac{\ln 
 (2T^2/\delta)}{t-l_{\tau}}} \leq \ln (2T^2/\delta)(1+\ln |\cE_i|) \leq \ln (2T^2/\delta)(1+\ln T)
\]

By Lemma~\ref{le:uncensoredlemmaS}, we know that $\tau_{max}\leq 2V^{2/3}T^{1/3}+1:=K$.
Moreover, due to restarts, we may divide an interval $\mathcal I'_{k'}$ into at most two different intervals $\mathcal I_{j}$ and $\mathcal I_{j'}$. It follows that 
\[
\sum_{\tau=1}^{\tau_{max}} \sum_{i=1}^{K_{\tau}} 1 \leq 2K +1.
\]

By the upper bound of regret incurred in each epoch, some calculation (similar to the proof of Theorem~\ref{thm:censored-nvp-ub-general}) gives
\begin{align*}
    \sum_{t=1}^T f_t (x_t) - f_t(x_t^*) &\leq  4C_f(V^{2/3}T^{1/3}+1) + 2\mu\sum_{t =1}^{T} \epsilon^2_t + 384\mu L_g^2 (V^{2/3}T^{1/3}+1)\ln T\ln^{1/2}(2T^2/\delta),
\end{align*}

where we complete the proof.
\Halmos
\endproof

\subsection{Assumption Verification for Applications in Section \ref{subsec:app}}\label{appendix:verification}

\subsubsection{Proof Omitted in Section \ref{subsec:app1}}
\label{appendix:app1}

This section discusses the verification of assumptions required by Theorems \ref{thm:general-ub} and \ref{thm:pl-ub} for newsvendor problems. We first prove a useful lemma for the verification of Assumptions \ref{ass:DKW-to-error} and \ref{ass:DKW-to-error2}.

\begin{lemma}
\label{lem:DKW-to-error}

Suppose $h(y)$ is $L_{h}$-Lipschitz.
Let $G_1(y)$ and $G_2(y)$ be two CDFs satisfying

\begin{enumerate}
\item $\sup_{y \in \R} \vert G_1(y) - G_2(y) \vert \leq \epsilon$  for some $\epsilon >0$.
\item There exist $a$ and $b$ such that $G(a) = 0$, and $G(b) = 1$.
\end{enumerate}
Then we have
\[
\left|\int_a^b h(y) \intd G_1(y)- \int_a^b h(y) \intd G_2(y)\right| \leq (b-a) L_{h}\epsilon. 
\]
\end{lemma}

\proof{Proof of Lemma \ref{lem:DKW-to-error}.}
By the integration by parts formula for Lebesgue-Stieltjes integrals, we know
\[
\int_a^b h(y) \intd G(y) = h(b)G(b)-h(a-)G(a-) - \int_a^b G(y) \intd h(y) .
\]

For functions $G_1(y)$ and $G_2(y)$, from the second condition of the lemma, we know $h(b)G_1(b)=h(b)G_2(b)$ and $h(a-)G_1(a-)=h(a-)G_2(a-)$. By $h(y)$ is $L_h$-Lipschitz, we know it is differentiable a.e.~and $|h'(y) \leq L_h$. It follows that
\begin{align*}
    \left|\int_a^b h(y) \intd G_1(y)- \int_a^b h(y) \intd G_2(y)\right| &\leq \left|\int_a^b G_1(y) \intd h(y)- \int_a^b G_2(y) \intd h(y)\right| \\
    &\leq \int_a^b |G_1(y)-G_2(y)|\cdot |h'(y)| \intd y\\
    &\leq (b-a)L_h\epsilon,
\end{align*}
where the last inequality is by $|h'(y)| \leq L_h$ and $\sup_{y \in \R} \vert G_1(y) - G_2(y) \vert \leq \epsilon$.
\Halmos\endproof

Next, we prove that, under mild conditions, all assumptions required by the regret analysis of Algorithm \ref{alg:saageneral} can be verified.  Recall, we assume that $x \in [0, \bar{x}]$ and $D \in [0,\bar{x}]$.

\begin{lemma}
\label{lem:verify1} We have the following conclusions.
\begin{enumerate}
    \item $\cX$ is convex and compact with diameter $D_{\cX} = \bar{x}$. For the linear cost problem, we know $C_f \leq (h+b)\bar{x}$ and for the quadratic cost problem, we know $C_f \leq a \bar{x}^2$.
     Thus, Assumption \ref{ass:bounded} holds.
    \item The empirical optimization problem is a convex optimization problem, and Assumption \ref{ass:opt-oracle} holds.
    \item For the linear cost problem, we can set $L_F = (h+b)\bar{x}$; for the quadratic cost problem, we can set $L_F = 2a\bar{x}^2$ such that Assumption \ref{ass:DKW-to-error} holds.
    \item For the quadratic cost problem, we can set $\mu = 1/a$ such that Assumption \ref{ass:pl} holds. For the linear cost problem, suppose $D_t,t \in [T]$ have density function $G'_t(y)$, and $G'_t(y) \geq \alpha >0,~\forall y \in [0,\bar{x}]$, and we can set $\mu = 1/2(h+b)\alpha$ such that Assumption \ref{ass:pl} holds.
    \item For $\epsilon \geq (h+b)/(j-i+1)$, we can find an optimization oracle $\O$ such that Assumption \ref{ass:opt-oracle2} holds.
    \item For the quadratic cost problem, we can set $L_g = 2a\bar{x}^2$ such that Assumption \ref{ass:pl} holds. For the linear cost problem, we can set $L_g = (h+b)$ such that Assumption \ref{ass:DKW-to-error2} holds.
\end{enumerate}
\end{lemma}

Note that for Assumption \ref{ass:opt-oracle2}, we restrict the accuracy $\epsilon$ to be greater than $(h+b)/(j-i+1)$ for the discontinuous derivative of the linear cost problem, and it is sufficient for the conclusion. By Theorem \ref{thm:pl-ub},  this optimization error will contribute to the regret a term of
\[
\mu \sum_{t=1}^T \epsilon^2_t \leq \frac{h+b}{2\alpha}\sum_{\tau=1}^{\tau_{max}} \sum_{i=1}^{K_\tau} 1/i^2 \leq \frac{5(h+b)}{6\alpha}\cdot \tau_{max} \leq \frac{(h+b)}{\alpha}\cdot \min\{2V^{2/3}T^{1/3}+1, S\},
\]
where the upper bound on the number of epochs is from the proof of Theorem \ref{thm:pl-ub}.

Plugging in the constant into Theorem \ref{thm:general-ub} and Theorem \ref{thm:pl-ub} gives the proof of Theorems \ref{thm:uncensored-nvp-ub-general} and \ref{thm:uncensored-nvp-ub-pl}

\proof{Proof of Lemma \ref{lem:verify1}.}

We know $F(x,d)$ is $(h+b)$-Lipschitz for the linear cost problem, and $2a\bar{x}$-Lipschitz for the quadratic cost problem. By Lemma \ref{lem:DKW-to-error}, we know that
\[
\left|\int_{0}^{\bar{x}} \theta (x-y ) \intd G(y )- \int_{0}^{\bar{x}} \theta (x-y ) \intd \hat G(y )\right| \leq L_{\theta}\bar{x} \sup_{y \in [0,\bar{x}]} |G(y )-\hat G(y )|. 
\]

We know Assumption \ref{ass:DKW-to-error} is satisfied with $L_F = L_{\theta}\bar{x}$, with $L_\theta$ depending on the problem.

Now, we consider Assumption \ref{ass:pl}. Note $\theta(\cdot)$ is linear cost or continuously differentiable. There exists subgradient function $\nabla F_x(x,D)$ for $F(x,D)$, $f_t(x)$ is differentiable, and $\E[\nabla F_x(x,D)] = \nabla f_t(x)$. For the quadratic cost problem, we know $f_t(x)$ is $2a$-strongly convex. Therefore, we can set $\mu = 1/a$. 
For the linear cost problem, we know $f_t(x)$ is $2(h+b)\alpha$-strongly convex and $\mu = 1/2(h+b)\alpha$.

Next, we consider Assumption~\ref{ass:opt-oracle2}. For the quadratic cost problem, we know the empirical optimization problem is also quadratic with an optimal solution in $[0,\bar{x}]$. The oracle is easy to find. For the linear cost problem, we know the empirical problem has a discontinuous gradient of
\[
h\hat{G}_{i,j}(x) - b(1-\hat{G}_{i,j}(x)).
\]
This is a stair function and the optimal solution has gradient $\hat{g}_{i,j}(x^*_{i,j}) \leq (h+b)/(j-i+1)$.

Finally, we consider Assumption \ref{ass:DKW-to-error2}. For the quadratic cost problem, we know $\theta'(\cdot) \leq 2a\bar{x}$. By Lemma \ref{lem:DKW-to-error}, we know that
\[
\left|\int_0^{\bar{x}} \theta'(x-y) \intd \hat G(y)- \int_0^{\bar{x}} \theta'(x-y) \intd G(y)\right| \leq 2a\bar{x}^2\epsilon. 
\]

Next, we consider the linear cost.
Some computation gives $$\E[ \theta'(x-D)] = hG(x) - b(1-G(x)).$$ Thus, if $\sup_{y \in \R}\vert \hat G(y)-G(y)\vert\leq \epsilon$, we have
\[
\left|\int_0^{\bar{x}} \theta'(x-y) \intd \hat G(y)- \int_0^{\bar{x}} \theta'(x-y) \intd G(y)\right| \leq (b+h)\cdot \sup_{y \in \R}\vert \hat G(y)-G(y)\vert \leq (h+b) \epsilon,
\]
where we complete the proof.
\Halmos
\endproof

\subsubsection{Proof Omitted in Section \ref{subsec:app2}}\label{appendix:app2}

For the problem with random yield, we know that $F(x,D,Z) = a(Zx-D)^2$, where $Z \in [0,1]$. Therefore, we know $f_t(x)$ would be strongly convex provided $\E[Z]>0$.

Therefore, we consider the random capacity problem, where $F(x,D,z) = a(\min\{x,Z\}-D)^2$. We will prove that under mild conditions, $f_t(x)$ can be transformed into a strongly convex function, and we can further prove that the PL condition is satisfied. The following transformed convexity definition is proposed by  \citep{fatkhullin2023stochastic}  to capture the convexity after transformation.

\begin{definition}
\label{def:hidden-cvx}
For each $t \in [T]$, $f_t(x)$ is $(\gamma,c)$-transformed strongly convex with $\gamma,c >0$ that is there exists a function $\varphi_t(x)$ such that
\begin{enumerate}
    \item $\varphi_t: \X \mapsto \R^n$ is an invertible function and let $\varphi_t^{-1}$ denote its inverse.
    \item The Jacobian $J(u) = \nabla_u \varphi_t^{-1}(u)$ exists and is bounded  $\Vert J(u) \Vert_2 \leq c$,
    for any $u \in \cU_t := \varphi_t(\X)$.
    \item Let $H_t(u) = f_t(\varphi_t^{-1}(u))$, then $H_t(u)$ is $\gamma$-strongly convex. 
\end{enumerate}    
\end{definition}

We would like to highlight the relationship between the transformed strong convexity and the PL condition. Define $u = \varphi_t(x)$. From Definition \ref{def:hidden-cvx} and the chain rule
, it follows that
\begin{align*}
    \Vert \nabla_u H_t(u) \Vert_2 =\Vert \nabla_u f_t(\varphi^{-1}(u)) \Vert_2 &\leq \Vert g_t(x) \Vert_2 \cdot \Vert J(u) \Vert_2 \leq c \Vert g_t(x) \Vert_2.
\end{align*}
By that $H_t(u)$ is strongly convex, recall $\nabla f_t(x) = g_t(x)$ and we have
\begin{align*}
    f_t\left({x}\right)-f_t(x^*_t) &= H_t(u) - H_t(u_t^*)\\
    &\leq \frac{1}{2\gamma} \Vert \nabla_uH_t\left({u}_t\right)\Vert^2_2 \\
    &\leq \frac{c^2}{2\gamma}\Vert \nabla f_t(x)\Vert^2_2 .
\end{align*}

Therefore, the PL condition is milder than the transformed strongly convex.

\begin{lemma}
    If $\P[Z \geq \bar{x}]>0$, we define expectation transformation $y = \varphi(x) = \E[\min\{x,Z\}]$. The tranformation $\varphi(\cdot)$ is bi-Lipschitz and $f_t(\varphi^{-1}(y))$ is strongly convex.
\end{lemma}
\proof{Proof.}
It is easy to see that
\begin{align*}
    f_t(x) &= \E[a(s(x,Z)-D)^2] \\
    &= a\E[s^2(x,Z)] -2\E[D]\E[s(q,Z)] + \E[D^2]\\
    &=a\E[s^2(x,Z)] -2\E[D]y + \E[D^2]
\end{align*}

Therefore, we only need to prove that $h_t(y):= \E[s^2(\varphi^{-1}(y),Z)]$ is strongly convex in $y$.

By the chain rule and some computation, we know that
\begin{align*}
    \frac{\partial h_t(y)}{\partial y} &= \frac{\partial h_t(y)}{\partial x} \cdot\frac{\partial x}{\partial y} \\
    &= 2x(1-G_Z(x)) \cdot \frac{1}{1-G_Z(x)} = 2x.
\end{align*}
Therefore, it follows that
\begin{align*}
    \frac{\partial^2 h_t(y)}{\partial y^2} &= 2\frac{\partial x}{\partial y}=\frac{2}{1-G_Z(x)} = 2/\P[Z \geq \bar{x}] > 0.
\end{align*}

Therefore, we know that $h_t(y)$ is strongly convex. Moreover, by $\P[Z \geq \bar{x}]>0$, it is easy to verify that $\varphi(x)$ is bi-Lipschitz and we complete the proof.
\Halmos\endproof

\subsubsection{Proof Omitted in Section \ref{subsec:app3}}\label{appendix:app3}
Next, we prove that all assumptions required in Sections \ref{subsec:general-ub} hold under mild conditions. 

\begin{lemma}
\label{lem:verify3}
Suppose that $D_t, t \in [T]$ are independent over time. There exists an oracle $\O_3$ that can find the optimal solution of $\min_{x \in \cX} \hat{f}_{i,j}(x)$ in polynomial time. We can set $\mathcal{X} = [0,1]$, $\mathcal{D}_{\mathcal{X}} = 1$, and $C_f = L_F = 1$ such that Assumptions \ref{ass:bounded}, \ref{ass:opt-oracle}, and \ref{ass:DKW-to-error} hold. 
\end{lemma}

\proof{Proof of Lemma \ref{lem:verify3}.} Note for each realization $D$, the function $F(x,D)= (x-v)\I[x \geq D]$ is a two-piece, linear function. Therefore, we know that $\hat{f}_{i,j}(x)$ is a $(j-i+2)$-piece linear function. Therefore, we only need to check the endpoints of all the pieces to find the optimal solution, which gives a polynomial-time oracle $\O_3$. Assumption \ref{ass:opt-oracle} is satisfied.

We set $\mathcal{X} = [0,1]$, which is convex and with diameter $\mathcal{D}_{\mathcal{X}} = 1$. We have $f_t(x) = \E[F(x,D_t)] \leq 1$
and we can set $C_f = 1$.  Assumption \ref{ass:bounded} is satisfied.  Noting that $(x-v)$ is deterministic, we know $\E[F(x,D_t)] = (x-v)G_t(x)$. It is easy to see that Assumption \ref{ass:DKW-to-error} is satisfied with $L_F = 1$. 
\Halmos
\endproof

\section{Details of Numerical Experiments}
\label{appendix:numerical-experiments}
\subsection{Description of Algorithms}
We consider the following algorithms in the literature.
\begin{enumerate}
    \item \textbf{Nonstationary Sample Average Approximation (NSAA):} This algorithm is our general NSAA algorithm applied to the newsvendor problem.
    \item \textbf{Sample Average Approximation (SAA):} This algorithm is the naive sample average approximation method for stationary distributions applied to the newsvendor problem.
    \item \textbf{Moving-window Sample Average Approximation (MSAA):} In period $t$, the algorithm uses the latest $n$ data to estimate the empirical quantile of the newsvendor problem, where $n = \lceil \kappa T^{1/2}\rceil$ and $\kappa$ is a hyper-parameter.
    \item \textbf{Restarting Sample Average Approximation (RSAA):}  This algorithm divides the time horizon into several epochs of $n$ period, where $n = \lceil \kappa T^{1/2}\rceil$ and $\kappa$ is a hyper-parameter. The algorithm will restart at the beginning of each epoch and use only the historical data from this epoch to estimate the empirical data.
    \item \textbf{Decaying-Weight Two-Stage Estimation (W2SE):} This algorithm assigns a decaying weight to historical data and uses a two-stage estimation to estimate the quantile. 
    \item \textbf{Moving-Window Two-Stage Estimation (M2SE):} This algorithm assigns equal weight to historical data in a moving window and uses a two-stage estimation to estimate the quantile. 
\end{enumerate}

MSAA and RSAA are heuristics used by \cite{keskin2023nonstationary} as a benchmark, and there is no theoretical guarantee for these heuristics. \cite{keskin2023nonstationary} also designed and provided theoretical guarantees for W2SE and M2SE algorithm. 

\subsection{Robustness to Critical Ratio}
\begin{figure}[h]
\begin{center}
\includegraphics[width =0.48\textwidth]{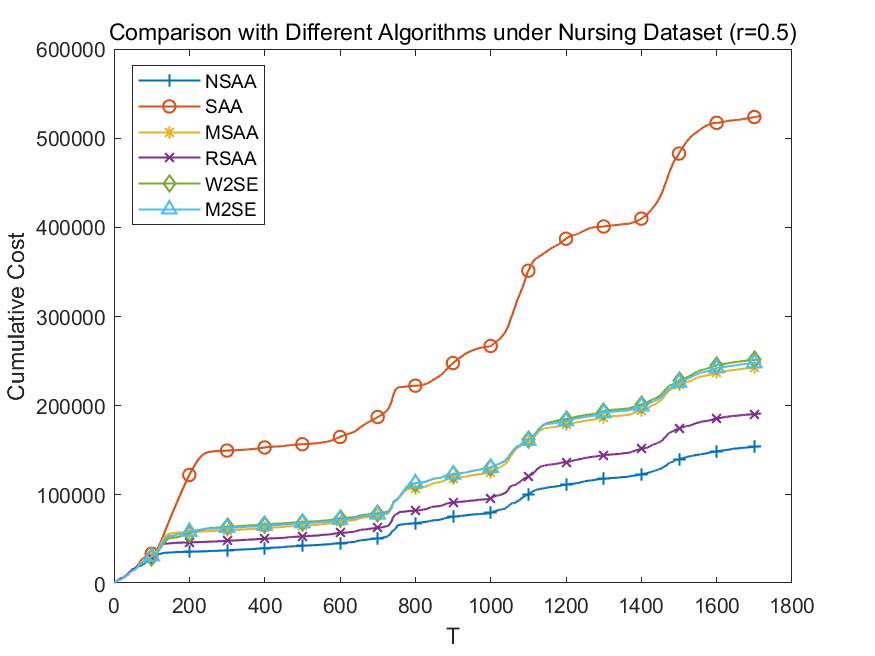}
\includegraphics[width =0.48\textwidth]{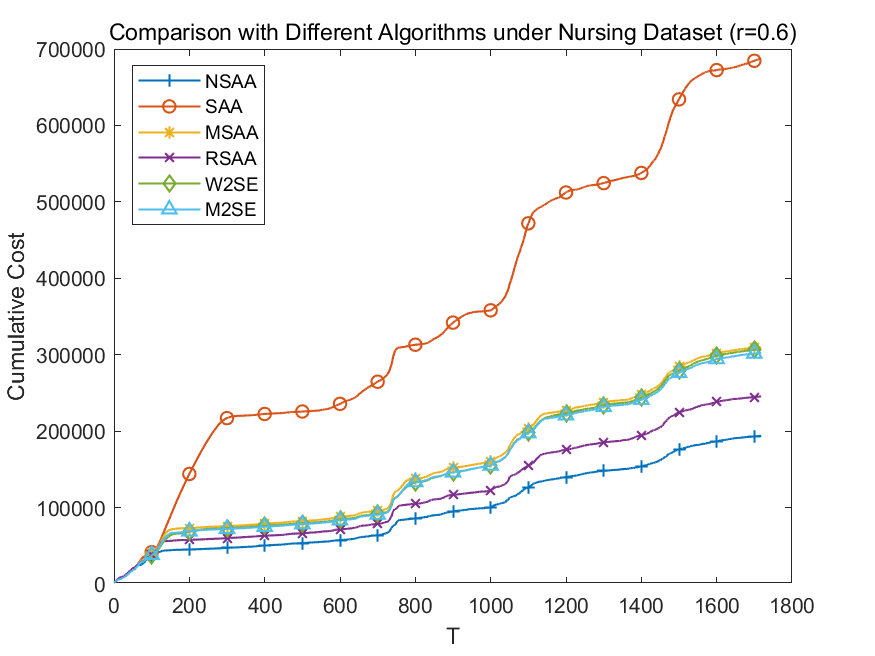}
\includegraphics[width =0.48\textwidth]{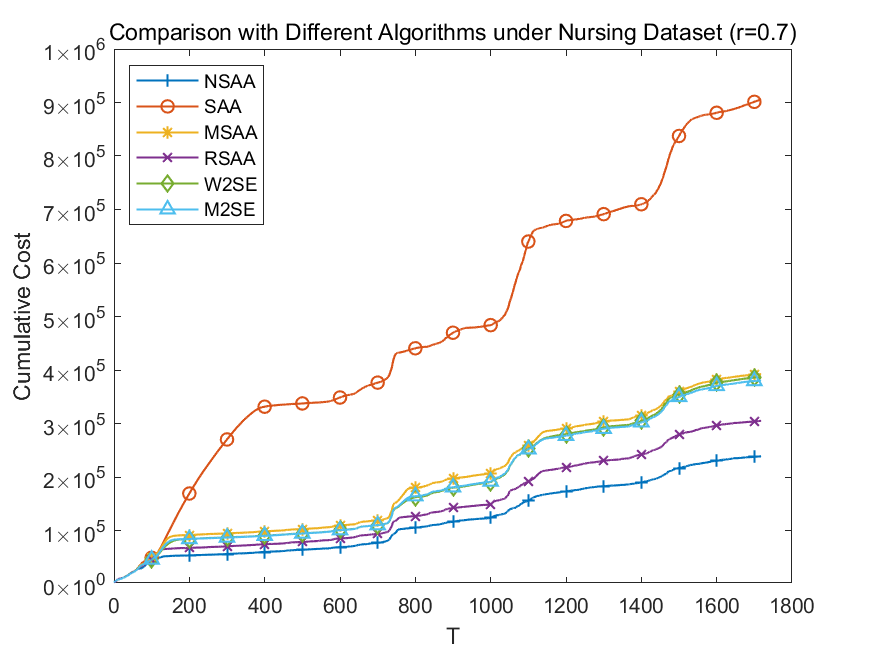}
\includegraphics[width =0.48\textwidth]{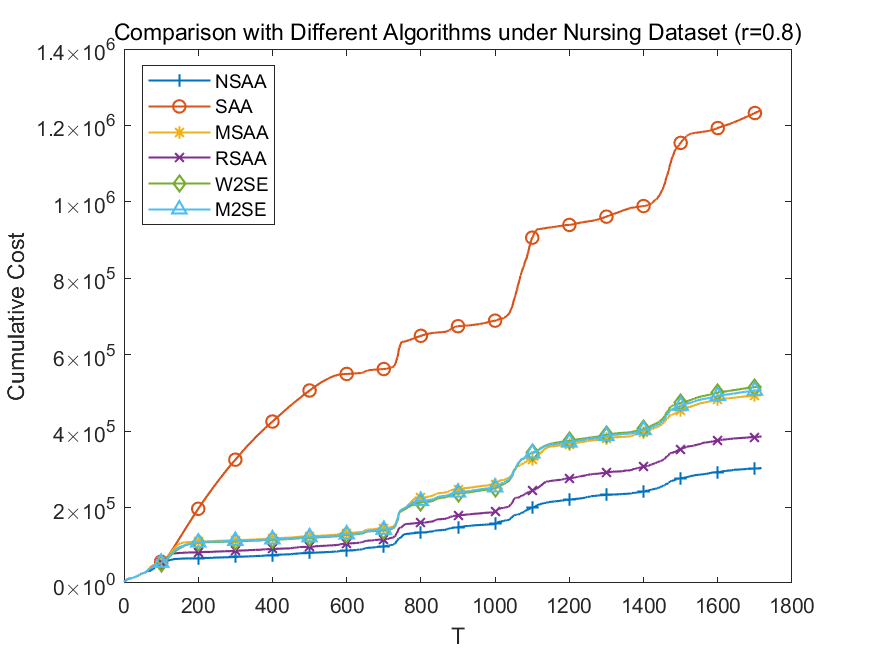}
\end{center}
\caption{{Performance of NSAA under Nursing Datasets and Different Ratio $r$}}\label{fig:robust1}
\end{figure}

\begin{figure}[h]
\begin{center}
\includegraphics[width =0.48\textwidth]{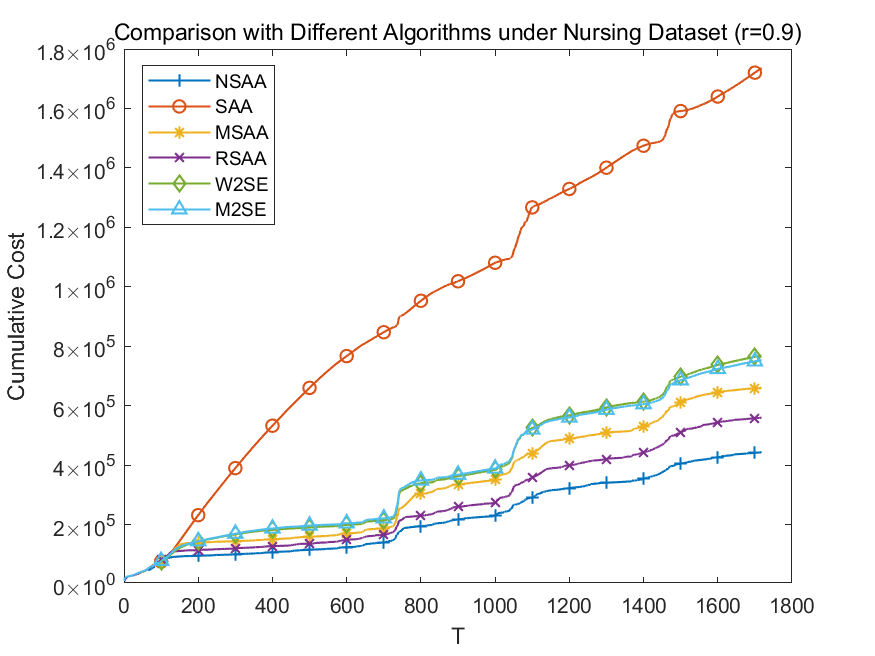}
\includegraphics[width =0.48\textwidth]{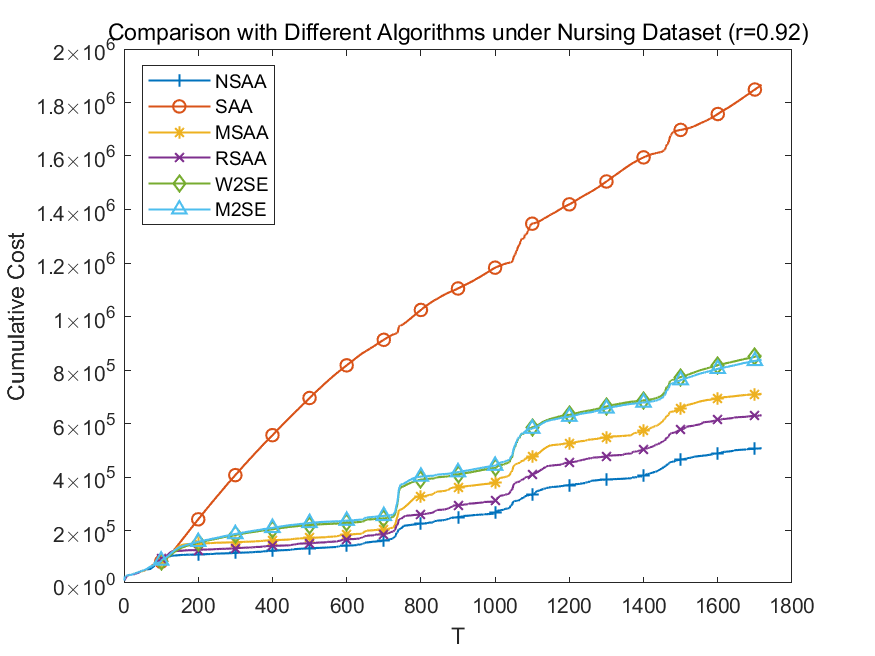}
\includegraphics[width =0.48\textwidth]{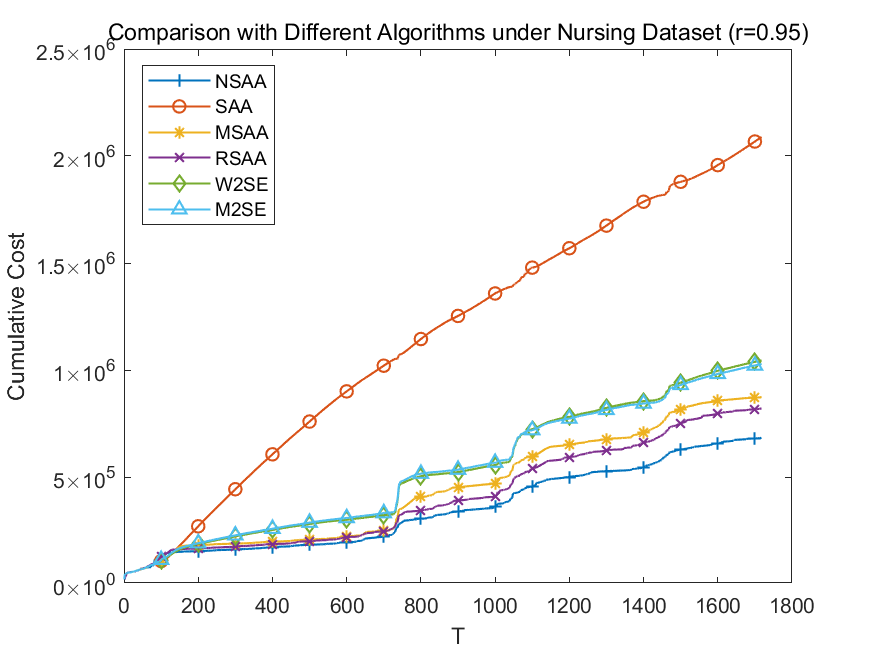}
\includegraphics[width =0.48\textwidth]{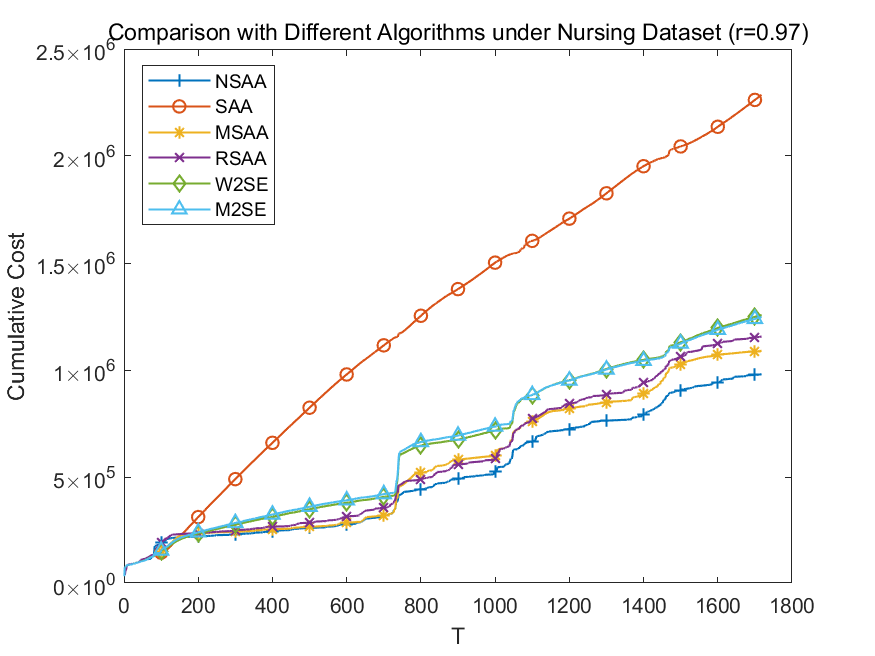}
\end{center}
\caption{{Performance of NSAA under Nursing Dataset and Large Ratio $r$}}\label{fig:robust2}
\end{figure}

In Figure~\ref{fig:robust1}, we present the results under critical ratios $r = 0.5,0.6,0.7,0.8$. In Figure~\ref{fig:robust2}, we present the results under large critical ratios $r = 0.9,0.92,0.95,0.97$.

As shown in Figures \ref{fig:robust1} and \ref{fig:robust2}, the NSAA algorithm achieves the best costs compared with other algorithms. These experiments show the robustness of the NSAA algorithm in parameter $b$.

\end{document}